\documentclass[12pt]{amsart}
\numberwithin{equation}{section}
\def\a{\alpha}

\def\b{\beta}

\newcommand\g{\gamma}

\newcommand\n{\mathbf n}

\newcommand\tosea{{}^t\hskip-2pt}
\newcommand\tose{{}^t\hskip-1pt}
\newcommand\T{{}^t\hskip-2pt}

\newcommand\w{\wedge}
\newcommand\x{\mathbf x}
\newcommand\y{\mathbf y}
\newcommand\z{\mathbf z}

\newcommand\C{\mathbf C}

\newcommand\I{\mathcal I}

\newcommand\R{\mathbf R}

\newtheorem*{introtheorem}{Classification Theorem}
\newtheorem{theorem}{Theorem}
\newtheorem{lemma}[theorem]{Lemma}
\newtheorem{corollary}[theorem]{Corollary}
\newtheorem{proposition}[theorem]{Proposition}

\theoremstyle{definition}
\newtheorem{definition}[theorem]{Definition}
\newtheorem{example}[theorem]{Example}

\theoremstyle{remark}
\newtheorem{remark}[theorem]{Remark}

\begin{document}

\title[Isoparametric hypersurfaces]{Isoparametric hypersurfaces with four principal curvatures}
\author[Cecil]{Thomas E. Cecil}
\thanks{The first author was partially supported by NSF Grant No. DMS-0071390}
\address{Department of Mathematics and Computer Science \\ College of the Holy Cross \\ Worcester, Massachusetts 01610-2395}
\email{cecil@mathcs.holycross.edu}
\author[Chi]{Quo-Shin Chi}
\thanks{The second author was partially supported by NSF Grant No. DMS-0103838}
\address{Department of Mathematics \\ Campus Box 1146 \\ Washington University \\ St. Louis, Missouri 63130}
\email{chi@math.wustl.edu}
\author[Jensen]{Gary R. Jensen}
\address{Department of Mathematics \\ Campus Box 1146 \\ Washington University \\ St. Louis, Missouri 63130}
\email{gary@math.wustl.edu}

\date{\today}


\begin{abstract}  Let $M$ be an isoparametric hypersurface in the sphere
$S^n$ with four distinct principal curvatures.  
M\"{u}nzner showed that the four principal curvatures 
can have at most two distinct multiplicities $m_1, m_2$,
and Stolz showed that the pair $(m_1,m_2)$ must either be 
$(2,2)$, $(4,5)$, or be equal to the multiplicities of an
isoparametric hypersurface of FKM-type, constructed by
Ferus, Karcher and M\"{u}nzner from orthogonal
representations of Clifford algebras.  In this paper, we
prove that if the multiplicities satisfy $m_2 \geq 3m_1 - 1$,
then the isoparametric hypersurface $M$ must be of
FKM-type.  Together with known results of Takagi
for the case $m_1 = 1$, and Ozeki and Takeuchi for $m_1 = 2$,
this handles all possible pairs of multiplicities except
for 10 cases, for which the classification problem remains open.
\end{abstract}

\maketitle

\section{Introduction}
A hypersurface $M$ in a real space-form $\tilde{M}^n (c)$ of 
constant sectional curvature $c$ is said
to be {\it isoparametric} if it has constant principal curvatures.
An isoparametric hypersurface $M$ in ${\mathbf R}^n$ can have at most
two distinct principal curvatures, and $M$ must be an open 
subset of a hyperplane, hypersphere or a spherical cylinder
$S^k \times {\mathbf R}^{n-k-1}$.  This was shown by Levi-Civita~\cite{Lev}
for $n=3$ and by B. Segre~\cite{Seg} for arbitrary $n$. Similarly,
E. Cartan \cite{Car1} proved that an isoparametric hypersurface $M$ 
in hyperbolic space $H^n$ can have
at most two distinct principal curvatures, and $M$ must be either
totally umbilic or else an open subset of a standard product
$S^k \times H^{n-k-1}$ in $H^n$ (see also \cite[pp.237-238]{CR}).
However, Cartan~\cite{Car1}-\cite{Car4}
showed in a series of four papers
written in the late 1930's that the situation is much more
interesting for isoparametric hypersurfaces in $S^n$. 
Cartan proved several general results and found examples with three and
four distinct principal curvatures, as well as those with one or two.
However, despite the beauty of Cartan's theory,
it was relatively unnoticed for thirty years,
until it was revived in the 1970's by 
Nomizu \cite{Nom1}-\cite{Nom2} and M\"{u}nzner~\cite{Mu}.

Cartan showed that isoparametric hypersurfaces
come as a family of parallel hypersurfaces, i.e., if
${\mathbf x}:M \rightarrow S^n$ is an isoparametric hypersurface, then so is 
any parallel hypersurface ${\mathbf x}_t$ at oriented
distance $t$ from the original
hypersurface ${\mathbf x}$.  However,
if $\lambda = \cot t$ is a principal curvature of $M$, then ${\mathbf x}_t$
is not an immersion, since it is constant on the
leaves of the principal foliation $T_{\lambda}$, and
${\mathbf x}_t$ factors through an immersion of the
space of leaves $M/T_{\lambda}$ into $S^n$.  In that case,
${\mathbf x}_t$ is a {\it focal submanifold} of codimension $m+1$ 
in $S^n$, where
$m$ is the multiplicity of $\lambda$.  M\"{u}nzner~\cite{Mu} showed that a
parallel family of isoparametric hypersurfaces in $S^n$ always
consists of the level sets in $S^n$ of a homogeneous polynomial $F$
defined on ${\mathbf R}^{n+1}$ satisfying certain differential
equations which are listed at the beginning of Section 2. 
He showed that the level sets of $F$ on $S^n$ are connected, and
thus
any connected isoparametric hypersurface can be extended to a unique
compact, connected isoparametric hypersurface.  M\"{u}nzner also showed that
regardless of the number of distinct principal curvatures
of $M$, there are only two distinct focal submanifolds in a
parallel family of isoparametric hypersurfaces, and each isoparametric
hypersurface in the family separates the sphere into two ball bundles
over the two focal submanifolds.  From this topological
information, M\"{u}nzner was able to prove his fundamental result that
the number $g$ of distinct principal curvatures of an isoparametric
hypersurface in $S^n$ must be $1,2,3,4$ or $6$.  As one would expect,
classification results on isoparametric hypersurfaces have 
been dependent on the number of distinct principal curvatures.

Cartan classified isoparametric hypersurfaces with 
$g \leq 3$ principal curvatures.  If $g=1$,
then $M$ is umbilic and it must be a great or small
sphere.  If $g=2$, then $M$ must be a standard product of two
spheres
\begin{displaymath}
S^k(r) \times S^{n-k-1}(s) \subset S^n, \quad r^2+s^2=1.
\end{displaymath}
In the case $g=3$, Cartan~\cite{Car2}
showed that all the principal curvatures must have the same multiplicity
$m=1,2,4$ or 8, and the isoparametric hypersurface must be a tube of
constant radius over a standard Veronese embedding of a projective
plane ${\mathbf F}P^2$ into $S^{3m+1}$, where ${\mathbf F}$ is the division 
algebra
${\mathbf R}$, ${\mathbf C}$, ${\mathbf H}$ (quaternions),
${\mathbf O}$ (Cayley numbers) for $m=1,2,4,8,$ respectively.  Thus, up to
congruence, there is only one such family for each value of $m$.

The classification of 
isoparametric hypersurfaces
with four or six principal curvatures has stood as one 
of the outstanding problems
in submanifold geometry for some time, and it was listed as Problem~34 on 
Yau's \cite{Yau} list of important open problems in geometry
in 1992. In this paper, we will provide a partial solution to this
classification problem in the case $g=4$, but first we will
describe the known results in the two cases.

In the case $g=6$, there exists
one homogeneous family
with six principal curvatures of multiplicity one in $S^7$,
and one homogeneous family with six principal curvatures of
multiplicity two in $S^{13}$ (see Miyaoka~\cite{Mi} for a description).  
These are the only known examples.  
M\"{u}nzner showed
that for $g=6$, all of the principal curvatures
must have the same multiplicity $m$, and then
Abresch \cite{Ab} showed  that $m$ must be 1 or 2. In the case
$m=1$, Dorfmeister and Neher~\cite{DN} showed 
in 1985 that an isoparametric 
hypersurface must be homogeneous, but it remains an
open question whether this is true in the case $m=2$. 

For $g=4$, there is  a much larger and more diverse
collection of known examples. 
Cartan produced examples of isoparametric hypersurfaces with four
principal curvatures
in $S^5$ and $S^9$.  These examples are homogeneous, and have the property
that all of the principal curvatures have the same multiplicity.
Cartan asked if all isoparametric hypersurfaces must be homogeneous,
and if there exists
an isoparametric hypersurface whose principal curvatures do not all have 
the same multiplicity.  Nomizu~\cite{Nom1} generalized Cartan's example
in $S^5$ to produce
a collection of isoparametric hypersurfaces
whose  principal curvatures have two distinct multiplicities $(1,k)$,
for any positive integer $k$, thereby answering Cartan's second
question in the affirmative.  At approximately the same
time as Nomizu's work, Takagi and Takahashi~\cite{TT} used the work
of Hsiang and Lawson~\cite{HL} on submanifolds of cohomogeneity two
to determine all homogeneous isoparametric hypersurfaces
of the sphere.  Takagi and Takahashi showed that every homogeneous
isoparametric hypersurface is a principal orbit of
the isotropy representation of a rank two symmetric space, and they
presented a complete list of examples.  This list included
some examples with 6 principal curvatures, as well as those with
$1,2,3$ or 4 distinct principal curvatures.  In a separate paper,
Takagi \cite{Ta} proved that in the case $g=4$,
if one of the principal curvatures of $M$
has multiplicity one, then $M$ must be homogeneous.

In a two-part paper, Ozeki and Takeuchi~\cite{OT} produced two infinite
series of inhomogeneous isoparametric hypersurfaces with 
multiplicities $(3,4k)$ and $(7,8k)$, for any positive integer $k$.
They also classified isoparametric hypersurfaces for which
one principal curvature has multiplicity two, proving
that they must be homogeneous.  In the process,
Ozeki and Takeuchi developed a formulation of the
Cartan-M\"{u}nzner polynomial $F$ in terms of the
second fundamental forms of the focal submanifolds which
is very useful in our work.

Next Ferus, Karcher and M\"{u}nzner~\cite{FKM} used representations of Clifford algebras to construct
for any positive integer
$m_1$ an infinite series of isoparametric hypersurfaces
with four principal curvatures having multiplicities
$(m_1,m_2)$, where $m_2$ is nondecreasing and
unbounded in each series. In fact, 
$m_2 = k \delta (m_1) - m_1 - 1$, where $\delta (m_1)$ is the
positive integer such that the Clifford algebra $C_{m_1-1}$
has an irreducible representation on ${\mathbf R}^{\delta (m_1)}$
(see~\cite{ABS}), 
and $k$ is any positive integer for which $m_2$ is positive.
Isoparametric hypersurfaces obtained by this construction
of Ferus, Karcher and M\"{u}nzner are said to be of
{\it FKM-type}.
The FKM-series with multiplicities $(3,4k)$ and $(7,8k)$ are
precisely those constructed by Ozeki and Takeuchi.
For isoparametric hypersurfaces of FKM-type, 
one of the focal submanifolds is always a Clifford-Stiefel
manifold (see Pinkall-Thorbergsson~\cite{PT}).

The set of FKM-type isoparametric hypersurfaces contains
all known examples 
with $g=4$ with the exception of two
homogeneous examples, with multiplicities $(m_1,m_2)$ equal to $(2,2)$ and
$(4,5)$ (see \cite[part II, p.27]{OT} for more detail on these two
exceptions).  Over the years, many restrictions on the
multiplicities were found by M\"{u}nzner~\cite{Mu},
Abresch~\cite{Ab}, Grove and Halperin~\cite{GH},
Tang~\cite{Tang} and Fang~\cite{Fang1}.  This series of papers
culminated in the recent work of Stolz~\cite{St}, who showed 
that the multiplicities of an isoparametric hypersurface with $g=4$
must be the same as those in the known examples of
Ferus, Karcher and M\"{u}nzner or the two homogeneous exceptions.
This certainly adds weight to the conjecture that the known
examples are actually the only isoparametric hypersurfaces
with $g=4$.  In this paper, we prove that this conjecture is true,
if the two multiplicities satisfy $m_2 \geq 3m_1 -1$.  Specifically,
we prove (see Theorem~\ref{th11.1}):  

\begin{introtheorem}
Let $M$ be an isoparametric hypersurface in the sphere $S^n$ with four
distinct principal curvatures, whose multiplicities $m_1$, $m_2$ satisfy
$m_2 \geq 3m_1 -1$. Then $M$ is of FKM-type.
\end{introtheorem}

Taken together with the classifications of Takagi for
the case $m_1 = 1$ and Ozeki and Takeuchi for $m_1 = 2$,
this handles all possible pairs $(m_1,m_2)$ of multiplicities,
with the exception of $(4,5)$ and 9 pairs of multiplicities,
$(3,4)$, $(4,7)$, $(5,10)$,
$(6,9)$, $(7,8)$, $(7,16)$, $(8,15)$,
$(9,22)$, $(10,21)$,
corresponding to isoparametric hypersurfaces of FKM-type.
For these 10 pairs, the classification problem for
isoparametric hypersurfaces remains open.  

The first part of this work (through \S9) gives necessary and sufficient
conditions in terms of a natural second order moving frame
for an isoparametric hypersurface to be of FKM-type.  The second part
shows that these conditions are satisfied if $m_2 \geq 3m_1 -1$.
It is entirely possible that many of the remaining open
cases can be resolved using our characterization of FKM-type,
but we have not been able to improve on our estimate
of $m_2 \geq 3m_1 -1$ at this point.

Next we will provide a detailed outline of the paper. For
more information on isoparametric hypersurfaces and the
extensive theory of isoparametric submanifolds
of codimension greater than one in the sphere which was
introduced by Carter and West~\cite{CW1} and Terng~\cite{Te1},
the reader is referred to the excellent survey
article by Thorbergsson~\cite{Th1}, who proved that all isoparametric submanifolds of codimension greater than
one in the sphere are homogeneous~\cite{Th2}.  

We think of an isoparametric hypersurface as an immersion $\tilde\x:M^{n-1} \to S^n$.  About any point of $M$ there is a neighborhood $U$ on which there is defined an orthonormal frame field $\tilde\x,\tilde e_0, e_a, e_p, e_\alpha, e_\mu$ for which $\tilde e_0$ is normal to the hypersurface and the other sets of vectors are principal directions for the four respective principal curvatures of $\tilde \x$.  The index range of $a,p$ has length $m$, and that of $\alpha, \mu$ has length $N$, where $m=m_1$ and $N = m_2$ are the multiplicities for our isoparametric hypersurface.  The dual coframe on $U$ is the set of $1$-forms $\theta^a, \theta^p, \theta^\alpha, \theta^\mu$ defined on $U$ by the equation (sum on repeated indices)
\begin{equation*}
d\tilde\x = \theta^a e_a + \theta^p e_p + \theta^\alpha e_\alpha + \theta^\mu e_\mu
\end{equation*}
The curvature surfaces are the integral submanifolds of the distribution obtained by setting any three sets of these forms equal to zero.  The Levi-Civita connection forms of a curvature surface are given, essentially, by the forms $\theta^a_b = de_a \cdot e_b$, $\theta^p_q = de_q \cdot e_p$, etc. 
The second fundamental tensors of the focal submanifolds are given in terms of our frame field by the four sets of tensors $F^\mu_{\alpha a}$, $F^\mu_{\alpha p}$, $F^\mu_{pa}$ and $F^\alpha_{pa}$ defined in~\eqref{eq:4:13} in which the coframe field $\omega^a, \omega^p, \omega^\alpha, \omega^\mu$ is defined in~\eqref{eq:4:10} as constant multiples of $\theta^a, \theta^p, \theta^\alpha, \theta^\mu$, respectively. 
We derive the identities imposed on these tensors and their derivatives by the Maurer-Cartan structure equations of the orthogonal group $O(n+1)$, the isometry group of $S^n$.

If our isoparametric hypersurface is of FKM-type, then a simple calculation shows that the following equations hold for an appropriate choice of the Darboux frame field.
\begin{align}
F^\mu_{\alpha \,a+m} &= F^\mu_{\alpha a}  \label{eq:1} \\
F^\alpha_{b+m\,a} + F^\alpha_{a+m\,b} &= 0 \label{eq:2} \\
F^\mu_{b+m\,a} + F^\mu_{a+m\,b} &= 0  \label{eq:3} \\
\theta^a_b - \theta^{a+m}_{b+m} &= L^a_{bc}(\omega^c + \omega^{c+m}), \quad L^a_{bc} = - L^b_{ac} = -L^a_{cb} \label{eq:4}
\end{align}
where $a,b,c = 1,\dots,m$ and $a+m, b+m$ run through the range of the indices $p,q$.  The matrices of the operators of the Clifford system in terms of our frame field have as entries certain constants and the functions $F^\mu_{\alpha a}$, $F^\mu_{\alpha p}$, $F^\mu_{pa}$, $F^\alpha_{pa}$ and $L^a_{bc}$.  Thus, using these matrices, we can define these operators for an arbitrary isoparametric hypersurface.  If equations~\eqref{eq:1}-\eqref{eq:4} hold for the isoparametric hypersurface, then by an elementary, but extremely long, calculation we show that these operators form a Clifford system whose FKM construction produces the given isoparametric hypersurface.  This calculation is contained in the proof of Theorem~\ref{th:g:1}.

In Proposition~\ref{pr:8:1} we prove that~\eqref{eq:1} implies~\eqref{eq:2}-\eqref{eq:4} on $U$ provided that $\tilde\x$ satisfies the \textit{spanning property} (Definition~\ref{spanning}), which is:  

\noindent
(a).  There exists a vector $x_\alpha e_\alpha$ such that 
\begin{equation*}
\{F^\mu_{\alpha a} x_\alpha y_\mu e_a: (y_\mu) \in \R^N\} = \mbox{span}\,\{e_1,\dots,e_m\}
\end{equation*}
(b).  There exists a vector $y_\mu e_\mu$ such that
\begin{equation*}
\{F^\mu_{\alpha a} x_\alpha y_\mu e_a : (x_\alpha) \in \R^N \} = \mbox{span}\, \{e_1,\dots,e_m\}
\end{equation*}

Combining these results, we see that if an isoparametric hypersurface satisfies the spanning property and~\eqref{eq:1} on $U$, then it is of FKM-type.  The next step is to see when~\eqref{eq:1} will be true.

The parallel hypersurface at an oriented distance $t$ from $\tilde\x$ is given by $\x=\cos t\, \tilde\x + \sin t\, \tilde e_0$.  Its unit normal vector is $e_0 = -\sin t\, \tilde\x + \cos t\, \tilde e_0$ and its principal directions are still given by the remaining vectors in the frame field.  
At some value of $t$ the rank of $\x$ is less than $n-1$, in which case the image of $\x$ is a focal submanifold of the isoparametric family.  Any multiple of $\pi/4$ added to this value of $t$ again gives a focal submanifold.  From M\"unzner's result that there are only two focal submanifolds, it follows that as $t$ changes by a multiple of $\pi/2$, we return to the same focal submanifold.  If $\x$ is a focal submanifold, then we may assume that $e_0, e_a$ is a normal frame field along $\x$ and the vectors $e_p$, $e_\alpha$, $e_\mu$ are the principal vectors for the second fundamental form $II_{e_0}$, of principal curvatures $0$, $1$ and $-1$, respectively.  Moving a distance $t=\pi/2$ from $\x$ along the geodesic in the direction of $e_0$, we arrive at $e_0$, which must then be a position vector on the same focal submanifold.  At $e_0$, the normal frame field is $\x, e_p$, and the principal vectors, of principal curvatures $0$, $1$ and $-1$ are $e_a$, $e_\alpha$ and $e_\mu$, respectively.

There is a simple relationship between the four sets of tensors at $e_0$, denoted with the same letters barred, and these tensors at $\x$.  For our purposes, the most important is
\begin{equation*}
\bar F^\mu_{\alpha a} = F^\mu_{\alpha \,a+m}
\end{equation*}
Use these tensors to define real bihomogeneous polynomials
\begin{equation*}
p_a(x,y) = F^\mu_{\alpha a}x_\alpha y_\mu, \quad \bar p_a(x,y) = \bar F^\mu_{\alpha a} x_\alpha y_\mu
\end{equation*}
In Proposition~\ref{pr:6:10} we prove that if $\x$ satisfies the spanning property on $U$ and if at each point of $U$ the $\bar p_a$ are contained in the ideal $I$ generated by $p_1, \dots, p_m$ in the polynomial ring $\R[x_\alpha, y_\mu]$, then the frame field can be chosen so that~\eqref{eq:1} holds on $U$.

The key to linking the set of polynomials $\bar p_a$ with the set of polynomials $p_a$ comes from a formula for the isoparametric function derived by Ozeki and Takeuchi~\cite{OT} (recorded in~\eqref{eq9.1} below).  In Proposition~\ref{pro9.2} (see also Proposition~\ref{pro9.3}) we use this formula to prove 
that the zero locus of $p_1,\dots,p_m$ in $\R P^{N-1} \times \R P^{N-1}$ is identical with that of $\bar p_1,\dots, \bar p_m$.  

Algebraic geometers have developed a substantial body of information about
the relationship between two polynomial ideals whose zero varieties coincide.  Let $I$ be the ideal generated by $p_1,\dots,p_m$ in the polynomial ring $\R[x_\alpha, y_\mu]$ and let $I^\C$ be the ideal they generate in the polynomial ring $\C[x_\alpha, y_\mu]$.  Define the affine bi-cones
\begin{align*}
V_I &= \{(x,y) \in \R^N \times \R^N : p_a(x,y)=0,\,\, a=1,\dots,m \} \\
V_I^\C &= \{(x,y) \in \C^N \times \C^N : p_a(x,y)=0,\,\, a=1,\dots,m \}
\end{align*}
Let $J_m$ be the complex subvariety of $V_I^\C$ where the Jacobian matrix of $p_1,\dots, p_m$ is of rank less than $m$.  In our Classification Theorem~\ref{th11.1} we prove the following.  Fix a point in $U$.  If the codimension of $J_m$ is greater than $1$ in $V_I^\C$, then, at the point,
\begin{enumerate}
\item[(I)] $p_1,\dots,p_m$ form a regular sequence in $\C [x_\alpha, y_\mu]$
\item[(II)] $\dim_\R V_I = \dim_\C V_I^\C$
\item[(III)] $I^\C$ is a prime ideal of codimension $m$
\item[(IV)]  The spanning property holds for $\x$.
\end{enumerate}
It follows then by Serre's criterion (see Proposition~\ref{pro10.2}) that the ideal $I$ is reduced (see Definition~\ref{reducedness}), which is precisely the condition which allows us to conclude that the $\bar p_a \in I$.

The final step in our argument is then provided by
Proposition~\ref{pro11.1} which states that for $m\geq 2$, if $N \geq
3m-1$ then $\mbox{codim}\,(J_m) \geq 2$ at every point of $U$.  The
proof of this estimate requires a detailed analysis of the second
fundamental forms $II_{e_a}$ of $\x$. In the case $m=1$, we give a
simpler proof that $M$ is of FKM-type, thereby providing another proof
of Takagi's result.

We would like to thank N. Mohan Kumar for substantial help with the algebraic geometry and
John Little for his comments on previous versions
of this paper.


\section{Second order frames}\label{section2}
An immersed connected oriented hypersurface $\tilde \x: M^{n-1} \to S^n$ is called \textit{isoparametric} if $\tilde \x$ has constant principal curvatures.  Such a hypersurface always occurs as part of a family, the level surfaces of an \textit{isoparametric function} $f$, which is a smooth function on $S^n$ such that $|\nabla f|^2 = a(f)$ and $\Delta f = b(f)$, for some smooth functions $a,b : \R\to \R$. 

Denote the principal curvatures of $\tilde \x$ by $k_i$, with multiplicity $m_i$, for $i = 1,\dots,g$, and assume that $k_1> \dots > k_g$.  M\"unzner~\cite[part I]{Mu} showed that the multiplicities satisfy $m_i = m_{i+2}$ (subscripts mod $g$).  He then showed that the isoparametric function $f$ must be the restriction to $S^n$ of a homogeneous polynomial $F:\R^{n+1} \to \R$ of degree $g$ satisfying the differential equations
\begin{equation*}
\aligned
|\mbox{grad}\,F|^2 &= g^2 r^{2g-2}, \quad r = |\x| \\
\Delta F &= \frac{m_2-m_1}2 g^2 r^{g-2}
\endaligned
\end{equation*}
where $m_1$ and $m_2$ are the two (possibly equal) multiplicities.
The polynomial $F$ is called the \textit{Cartan-M\"unzner polynomial}
of the family of isoparametric hypersurfaces, and $F$ takes values
between $-1$ and $1$ on the sphere $S^n$.  For $-1 <t<1$, the level
set $F^{-1}(t)$ is one of the isoparametric hypersurfaces in the
family.  The level sets $M_+ = F^{-1}(1)$ and $M_- = F^{-1}(-1)$ are
the two focal submanifolds of the family, having codimensions $m_1+1$
and $m_2+1$ in $S^n$, respectively. 

We now develop the local geometry of isoparametric hypersurfaces using
the method of moving frames in the sphere.  In the process, we will
reprove some of the results obtained by M\"unzner, although this is
not our primary goal. 

We assume now that $g=4$, even though many of the results in Sections
$2$ -- $4$ have analogues for arbitrary values of $g$.  Let $e_0$ be
the unit normal vector field along $\tilde \x$ defining the
orientation of $M$.  Any point of $M$ has an open neighborhood $U$ on
which there exists a Darboux frame field $\tilde \x, e_i, \tilde e_0 :
U \to SO(n+1)$, $1 \leq i \leq n-1$, for which each vector $e_i$ is a
principal direction.  We adopt the index ranges 
\begin{equation}\label{eq:2:0}
\aligned
&i,j,k \in \{1,\dots,n-1\} \\
&a,b,c \in \{1,\dots,m_1\}, \quad p,q,r \in \{m_1+1,\dots, m_1+m_3\} \\
&\alpha, \beta, \g \in \{m_1+m_3+1,\dots,m_1+m_2+m_3\} \\
&\mu, \nu, \sigma \in \{m_1+m_2+m_3+1,\dots, n-1\}
\endaligned
\end{equation}
Arrange the frame so that the $e_a$ span the principal space for
$k_1$, the $e_\alpha$ span the principal space for $k_2$, the $e_p$
span the principal space for $k_3$, and the $e_\mu$ span the principal
space for $k_4$.  We shall call such a Darboux frame field 
\begin{equation}\label{eq:2:0a}
\tilde \x, e_a, e_p, e_\alpha, e_\mu, \tilde e_0
\end{equation}
on $U$ a \textit{second order frame field} along $\tilde \x$, (a first
order Darboux frame field is one for which $e_0$ is normal and the
remaining vectors are tangent, but not necessarily principal
directions). 
For such a frame field
\begin{equation}\label{eq:2:1}
d\tilde\x = \theta^i e_i \mbox{ and } de_i = \theta^j_i e_j - \theta^i \tilde \x + \theta^0_i \tilde e_0
\end{equation}
where $\theta^i$, $\theta^0_i = -\theta^i_0$, $\theta^i_j = - \theta^j_i$ are $1$-forms on $U$ and
$\theta^1,\dots,\theta^{n-1}$ is an orthonormal coframe field on $U$ with respect to the metric induced by $\tilde \x$ on $M$.  Notice that $\theta^0 = d\tilde \x \cdot \tilde e_0 =0$.  We use the summation convention unless the contrary is stated explicitly.
These $1$-forms satisfy the Maurer-Cartan structure equations of $SO(n+1)$
\begin{equation}\label{eq:2:1a}
\aligned
d\theta^i &= -\theta^i_j \w \theta^j \\ 
d\theta^0_i &= -\theta^0_j \w \theta^j_i \\ 
d\theta^i_j &= \theta^i \w \theta^j - \theta^i_0 \w \theta^0_j - \theta^i_k \w \theta^k_j
\endaligned
\end{equation} 
We also have
\begin{equation}\label{eq:2:2}
d\tilde e_0 = \theta^i_0 e_i
\end{equation}
where the $1$-forms $\theta^i_0 = -\theta^0_i$ are linear combinations of the coframe forms, namely
\begin{equation}\label{eq:2:3}
\theta^0_i =  h_{ij} \theta^j
\end{equation}
where these coefficient functions on $U$ satisfy $h_{ij}=h_{ji}$ as a consequence of taking the exterior derivative of the equation $\theta^0=0$.  The second fundamental form of $\tilde \x$ is
\begin{equation}\label{eq:2:4}
\widetilde{II} = -d\tilde \x \cdot d\tilde e_0 =  h_{ij}\theta^i\theta^j
\end{equation}
Having chosen the $e_i$ to be principal vectors, we know that the symmetric matrix $h_{ij}$ is a diagonal matrix.  In fact, we have
\begin{equation}\label{eq:2:5}
\theta^0_a = k_1 \theta^a,\quad  \theta^0_p = k_3 \theta^p,\quad \theta^0_\alpha = k_2 \theta^\alpha,\quad \theta^0_\mu = k_4 \theta^\mu
\end{equation}
Set $\theta^i_j = \sum h^i_{jk} \theta^k$, where the smooth function coefficients satisfy $h^i_{jk} = - h^j_{ik}$, for all $i,j,k = 1,\dots,n-1$.
Take the exterior differential of equations~\eqref{eq:2:5}, using the structure equations of $SO(n+1)$, to find
\begin{equation}\label{eq:2:6}
\aligned
\theta^p_a &= h^p_{a\alpha}
\theta^\alpha + h^p_{a\mu} \theta^\mu, \mbox{ since $h^p_{ab}=0=h^p_{aq}$} \\
\theta^\alpha_a &= h^\alpha_{ap} \theta^p + h^\alpha_{a\mu} \theta^\mu, \mbox{ since $h^\alpha_{ab} = -h^a_{\alpha b} =0=h^\alpha_{a\beta}$} \\
\theta^\mu_a &= h^\mu_{ap} \theta^p + h^\mu_{a\alpha} \theta^\alpha, \mbox{ since $h^\mu_{ab} = -h^a_{\mu b} = 0 = h^\mu_{a\nu}$} \\
\theta^\alpha_p &= h^\alpha_{pa} \theta^a + h^\alpha_{p\mu} \theta^\mu, \mbox{ since $h^\alpha_{pq} = -h^p_{\alpha q} = 0 = h^\alpha_{p\beta}$} \\
\theta^\mu_p &= h^\mu_{pa} \theta^a + h^\mu_{p\alpha} \theta^\alpha, \mbox{ since $h^\mu_{pq} = -h^p_{\mu q} = 0 = h^\mu_{p\nu}$} \\
\theta^\mu_\alpha &= h^\mu_{\alpha a} \theta^a + h^\mu_{\alpha p} \theta^p, \mbox{ since $h^\mu_{\alpha\beta} = -h^\alpha_{\mu\beta} = 0 = h^\mu_{\alpha\nu}$}
\endaligned
\end{equation}
and the coefficient functions further satisfy
\begin{equation}\label{eq:2:7}
\aligned
(k_3-k_1) h^p_{a\alpha} &= (k_2 - k_1) h^\alpha_{ap} = (k_2-k_3) h^\alpha_{pa} \\
(k_3-k_1)h^p_{a\mu} &= (k_4-k_1) h^\mu_{ap} = (k_4-k_3) h^\mu_{pa} \\
(k_2-k_1) h^\alpha_{a\mu} &= (k_4-k_1) h^\mu_{a\alpha} = (k_4-k_2) h^\mu_{\alpha a} \\
(k_2-k_3) h^\alpha_{p\mu} &= (k_4 - k_3) h^\mu_{p\alpha} = (k_4-k_2) h^\mu_{\alpha p} 
\endaligned
\end{equation}
At a point of $M$ the set of principal vectors for a
principal curvature $k_i$ is a subspace of dimension $m_i$, defined by
the equations $\theta^j=0$, for all $j$ not in the range of the given
principal curvature.  This
$m_i$-plane distribution on $M$ is called a
\textit{curvature distribution} on $M$.

\begin{lemma}\label{le:2:1} The curvature distributions are completely
integrable.  Their integral submanifolds are called
\textit{curvature
surfaces}.  A curvature surface corresponding to $k_j$ is totally 
geodesic in $M$ and its induced metric has constant sectional
curvature $1+k_j^2$.
\end{lemma}

\begin{proof} This is a simple application of
the structure equations and the first
three equations in~\eqref{eq:2:6}. 
\end{proof}

Additional conditions are imposed by the structure equations on the
coefficients upon the exterior differentiation of
equations~\eqref{eq:2:6}.

\section{Parallel hypersurfaces}\label{section3}
Let $\tilde \x, e_a, e_p, e_\alpha, e_\mu, \tilde e_0$ be a second order frame field~\eqref{eq:2:0a} along $\tilde \x$ on $U$.
We may arrange to have $k_1>k_2>k_3>k_4$.  It will be convenient to set $k_i = \cot s_i$, for $i=1,\dots,4$, where $0< s_1<s_2<s_3<s_4 < \pi$.  For any fixed real number $t$, let
\begin{equation}\label{eq:3:1}
\x = \cos t\, \tilde\x + \sin t \,\tilde e_0
\end{equation}
From~\eqref{eq:2:1}, \eqref{eq:2:2} and \eqref{eq:2:5} we have
\begin{equation}\label{eq:3:2}
\aligned
d \x = &(\cos t - \sin t \cot s_1)\theta^a e_a + (\cos t - \sin t \cot s_3) \theta^p e_p \\
&+ (\cos t - \sin t \cot s_2) \theta^\alpha e_\alpha + (\cos t - \sin t \cot s_4) \theta^\mu e_\mu
\endaligned
\end{equation}
We conclude that $\x$ is an immersion of $M$ except when $t \equiv s_i \mod \pi$, for some $i=1,2,3,4$.  Suppose $t$ is not one of these exceptional values.  Then the unit normal vector field along $\x$ preserving the orientation of $M$ is
\begin{equation}\label{eq:3:4}
e_0 = -\sin t\, \tilde \x + \cos t \,\tilde e_0
\end{equation}
and again from~\eqref{eq:2:1}, \eqref{eq:2:2} and \eqref{eq:2:5} we have
\begin{equation}\label{eq:3:5}
\aligned
d e_0 = &-(\sin t + \cos t \cot s_1) \theta^a e_a -(\sin t + \cos t \cot s_3) \theta^p e_p \\
&- (\sin t + \cos t \cot s_2) \theta^\alpha e_\alpha - (\sin t + \cos t \cot s_4) \theta^\mu e_\mu
\endaligned
\end{equation}
Since $(\sin t + \cos t \cot s)/(\cos t - \sin t \cot s) = \cot(s-t)$, for any $s$ and $t$, we find that the second fundamental form of $\x$ is
\begin{equation}\label{eq:3:6}
\aligned
II &= -d \x \cdot d e_0  \\ 
&= \cot(s_1-t)\, \omega^a \omega^a + \cot(s_3-t)\, \omega^p \omega^p \\ 
&+ \cot(s_2-t)\, \omega^\alpha \omega^\alpha + \cot(s_4-t)\, \omega^\mu \omega^\mu
\endaligned
\end{equation}
We conclude that the principal curvatures of $\x$ are constant, equal to $\cot(s_i-t)$ with multiplicity $m_i$, for $i=1,2,3,4$, and that
\begin{equation}\label{eq:3:7}
\x, e_a, e_p, e_\alpha, e_\mu, e_0
\end{equation}
is a second order frame field along $\x$ on $U$.

\section{Focal submanifolds}\label{section4}
We consider now what happens when $t$ is one of the exceptional values.  To be specific, suppose that $t=s_1$.  Then $\x$ is defined in~\eqref{eq:3:1} and $e_0$ is defined in~\eqref{eq:3:4} with $t = s_1$.  For the frame field~\eqref{eq:3:7} along $\x$ on $U$, equation~\eqref{eq:3:2} becomes
\begin{equation}\label{eq:4:1a}
d\x = \omega^p e_p + \omega^\alpha e_\alpha + \omega^\mu e_\mu
\end{equation}
whose rank is $n-1-m_1$ at every point of $M$ and where
\begin{equation}\label{eq:4:1}
\aligned
\omega^p &= \frac{\sin(s_3-s_1)}{\sin s_3} \theta^p,& \omega^\alpha &= \frac{\sin(s_2-s_1)}{\sin s_2} \theta^\alpha, & 
\omega^\mu &= \frac{\sin(s_4 - s_1)}{\sin s_4} \theta^\mu
\endaligned
\end{equation}
Therefore, the image $\tilde\x(M)$ is a submanifold of codimension $m_1+1$ in $S^n$.  It is called the \textit{focal submanifold} for the principal curvature $\cot s_1$.  In the same way, there are focal submanifolds for each of the principal curvatures.  For a point $\mathbf v \in \x(M)$, the set $L = \x^{-1}\{\mathbf v \}$ is a curvature surface of $\x$ for the principal curvature $\cot s_1$.  Restricted to this curvature surface, the forms $\theta^a$ give a coframe field on it.

If $e_0$ is defined by~\eqref{eq:3:4}, then~\eqref{eq:4:1a} shows that $\x, e_p,e_\alpha, e_\mu, e_a, e_0$ is a Darboux frame field along $\x$, with $e_p, e_\alpha, e_\mu$ tangent and $e_0, e_a$ normal vectors.  Take a point $p$ in the curvature surface $L$ and let $N$ denote the normal space to $\x$ at $p$.  Let $S^{m_1}$ denote the unit sphere in $N$.  The next lemma shows that $e_0(L)$ covers an open neighborhood of $e_0(p)$ in this sphere.

\begin{lemma}\label{le:4:1}  The rank of $e_0 : L \to S^{m_1}$ is $m_1$ at every point of the curvature surface $L$.  Therefore, $e_0(L)$ covers an open neighborhood of $e_0(p)$ in $S^{m_1}$.
\end{lemma}

\begin{proof}  Consider the frame field $e_0, e_a, \x, e_p, e_\alpha, e_\mu$ along $e_0$ on $L$.  Since $\theta^p$, $ \theta^\alpha$ and $\theta^\mu$ are all zero pulled back to $L$, it follows from~\eqref{eq:2:6} that $\theta^p_0$, $\theta^\alpha_0$ and $\theta^\mu_0$ are also zero pulled back to $L$.  Therefore, restricted to $L$, and using~\eqref{eq:2:5}, in which now $k_1 = \cot s_1$, we have
\begin{equation}\label{eq:4:2}
d e_0 = -\sin s_1\,\theta^a e_a + \cos s_1\, \theta^a_0 e_a = -\csc s_1\, \theta^a e_a
\end{equation}
which has rank equal to $m_1$ at every point of $L$.
\end{proof}

We can now calculate the second fundamental form of the submanifold $\x$ at the point $\x(p) = \mathbf v$ with respect to any unit normal vector there.

\begin{lemma}\label{le:4:2}  At any point of $M$ and with respect to any unit normal vector at the point, the principal curvatures of the focal submanifold $\x$ are 
\begin{equation}\label{eq:4:2a}
\cot (s_2-s_1), \quad\cot (s_3-s_1), \quad\cot (s_4-s_1)
\end{equation}
with multiplicities $m_2, m_3, m_4$, respectively.
\end{lemma}

\begin{proof}  From~\eqref{eq:3:5} we have for $t = s_1$
\begin{equation}\label{eq:4:3}
\aligned
d e_0 = &-\frac 1{\sin s_1} \theta^a e_a - \frac{\cos(s_3-s_1)}{\sin s_3} \theta^p e_p \\
&- \frac{\cos(s_2-s_1)}{\sin s_2} \theta^\alpha e_\alpha - \frac{\cos(s_4-s_1)}{\sin s_4} \theta^\mu e_\mu
\endaligned
\end{equation}
Combining this with~\eqref{eq:4:1} we have for the second fundamental form at $p$ with respect to the normal vector $e_0$
\begin{equation*}\label{eq:4:4}
\aligned
II_{e_0} &= -d\x \cdot d e_0 \\
&= \cot(s_3-s_1) \omega^p \omega^p + \cot(s_2 - s_1) \omega^\alpha \omega^\alpha + \cot(s_4 - s_1) \omega^\mu \omega^\mu
\endaligned
\end{equation*}
where $\omega^p, \omega^\alpha, \omega^\mu$, defined in~\eqref{eq:4:2},
form an orthonormal coframe with respect to the metric induced by $\x$ on the orthogonal complement to the curvature surfaces of the principal curvature $\cot s_1$.  By Lemma~\ref{le:4:1} we know that $e_0(L)$ covers some open subset of the unit sphere in the normal space to $\x$ at $p$.  Since the characteristic polynomial of $II_\n$ is an analytic function of $\n$ in the unit sphere of the normal space, it follows that the eigenvalues of $II_\n$ must be given by~\eqref{eq:4:2a} for every unit normal vector at $p$. (See~\cite[Proof of Corollary 2.2 on p. 249]{CR}).
\end{proof}

M\"unzner~\cite[Part I]{Mu} proved Lemma~\ref{le:4:2} and used it to prove the following important consequence (see also~\cite[p. 249]{CR}).

\begin{corollary}\label{co:4:1} The angles $s_i = s_1 + (i-1)\pi/4$, for $i = 2,3,4$ and the multiplicities satisfy $m_1 = m_3$ and $m_2=m_4$.  To simplify the notation we set $m_1=m_3=m$ and $m_2 = m_4 = N$.
\end{corollary}
Given these facts, our index conventions~\eqref{eq:2:0} become
\begin{equation}\label{eq:4:convention}
\aligned
i,j,k &\in \{1,\dots,n-1\},\quad a,b,c \in \{1,\dots,m\} \\
p,q,r &\in \{m+1,\dots, 2m\}, \quad 
\alpha, \beta, \g \in \{2m+1,\dots,2m+N\} \\
\mu, \nu, \sigma &\in \{2m+N+1,\dots n-1\}
\endaligned
\end{equation}
so that $2m+2N = n-1$, and $n$ must be odd.  
Combining Lemma~\ref{le:4:2} and Corollary~\ref{co:4:1} yields the following.

\begin{corollary}\label{co:4:2}  At any point of $M$ and with respect to any unit normal vector of $\x$ at the point, the principal curvatures of $\x$ are
\begin{equation}\label{eq:4:5a}
1,\quad 0, \quad -1
\end{equation}
with multiplicities $N$, $m$ and $N$, respectively.
\end{corollary}

In the light of Corollary~\ref{co:4:1}, the principal curvatures $k_i
= \cot s_i$ of $\tilde \x$ satisfy
\begin{equation}\label{eq:4:6}
k_2 = \frac{k_1-1}{k_1+1}, \quad k_3 = -\frac1{k_1}, \quad k_4 = \frac{1+k_1}{1-k_1}
\end{equation}
We will have occasion to use the following differences of these principal curvatures.
\begin{equation}\label{eq:4:7}
\aligned
k_2-k_1 &= -\frac{1+k_1^2}{1+k_1}, &\quad k_3 - k_1 &= -\frac{1+k_1^2}{k_1} \\ 
k_4 - k_1 &= \frac{1+k_1^2}{1-k_1}, &\quad
k_3-k_2 &= -\frac{1+k_1^2}{k_1(1+k_1)} \\
k_4 - k_2 &= 2\frac{1+k_1^2}{1-k_1^2},  &\quad k_4 - k_3 &= \frac{1+k_1^2}{k_1(1-k_1)}
\endaligned
\end{equation}
We use equations~\eqref{eq:4:7} to rewrite equations~\eqref{eq:2:7} as
\begin{equation}\label{eq:4:7a}
\aligned
h^p_{a\alpha} &= -\frac1{1+k_1} h^\alpha_{pa}, &\quad h^\alpha_{ap} &= -\frac1{k_1} h^\alpha_{pa} \\
h^p_{a\mu} &= \frac 1{k_1-1} h^\mu_{pa}, &\quad h^\mu_{ap} &= \frac 1{k_1} h^\mu_{pa} \\
h^\alpha_{a\mu} &= \frac2{k_1-1} h^\mu_{\alpha a}, &\quad h^\mu_{a\alpha} &= \frac 2{1+k_1} h^\mu_{\alpha a} \\
h^\alpha_{p\mu} &= \frac{2k_1}{1-k_1} h^\mu_{\alpha p}, &\quad h^\mu_{p\alpha} &= \frac{2k_1}{1+k_1} h^\mu_{\alpha p}
\endaligned
\end{equation}
Now, with $s_i = s_1 + (i-1)\pi/4$, equation~\eqref{eq:4:1a} takes the form
\begin{equation}\label{eq:4:8}
d\x = \frac1{\sin s_3} \theta^p e_p + \frac1{\sqrt 2 \sin s_2} \theta^\alpha e_\alpha + \frac 1{\sqrt 2 \sin s_4} \theta^\mu e_\mu
\end{equation}
and with $t=s_1$ equation~\eqref{eq:3:5} becomes
\begin{equation}\label{eq:4:9}
de_0 = -\frac1 {\sin s_1} \theta^a e_a - \frac 1{\sqrt2 \sin s_2} \theta^\alpha e_\alpha + \frac 1 { \sqrt 2 \sin s_4} \theta^\mu e_\mu
\end{equation}
If we define a new coframe field on $U\subset M$ by
\begin{equation}\label{eq:4:10}
\aligned
\omega^a &= -\frac 1{\sin s_1} \theta^a, &\quad \omega^p &= \frac1{k_1 \sin s_1} \theta ^p \\  \omega^\alpha &= \frac 1{(1+k_1) \sin s_1} \theta^\alpha, &\quad \omega^\mu &= \frac 1{(k_1 -1)\sin s_1} \theta^\mu
\endaligned
\end{equation}
then, because
\begin{equation}\label{eq:4:11}
\sin s_2 = \frac{1+ k_1}{\sqrt 2} \sin s_1, \,\, \sin s_3 = k_1 \sin s_1, \,\, \sin s_4 = \frac{k_1-1}{\sqrt 2} \sin s_1
\end{equation}
equations~\eqref{eq:4:8} and~\eqref{eq:4:9} become
\begin{equation}\label{eq:4:12}
d \x = \omega^p e_p + \omega^\alpha e_\alpha + \omega^\mu e_\mu, \quad 
de_0 = \omega^a e_a - \omega^\alpha e_\alpha + \omega^\mu e_\mu 
\end{equation}
One conclusion we can draw from~\eqref{eq:4:12} is that
\begin{equation}\label{eq:4:12a}
\x, e_0, e_a, e_p, e_\alpha, e_\mu
\end{equation}
is a Darboux frame field along $\x$ on $U$, with $e_0, e_a$ normal vectors and $e_p, e_\alpha, e_\mu$ tangent vectors spanning the principal spaces of curvature $0$, $1$ and $-1$, respectively of $II_{e_0}$.  We shall call this a \textit{second order frame field} along the focal submanifold $\x$ on $U$.  
For each point of $U$, define linear subspaces of $\R^{n+1}$ by
\begin{equation}\label{eq:6:4}
V_+ = \mbox{span}\{e_\alpha\}, \quad V_- = \mbox{span}\{e_\mu\}, \quad V_0 = \mbox{span}\{e_p\}
\end{equation}
These are the $+1$, $-1$ and $0$ principal curvature spaces, respectively, for the normal vector $e_0$ at this point.  
If we express the Maurer-Cartan forms~\eqref{eq:2:6} in terms of our coframe field~\eqref{eq:4:10} as
\begin{equation}\label{eq:4:13}
\aligned
\theta^p_a &= F^\alpha_{pa} \omega^\alpha - F^\mu_{pa} \omega^\mu, &\quad \theta^\alpha_a &= F^\alpha_{pa}\omega^p - 2F^\mu_{\alpha a}\omega^\mu \\
\theta^\alpha_p &= F^\alpha_{pa}\omega^a - 2F^\mu_{\alpha p}\omega^\mu, &\quad \theta^\mu_a &= -F^\mu_{pa}\omega^p -2 F^\mu_{\alpha a} \omega^\alpha \\
\theta^\mu_p &= F^\mu_{pa}\omega^a + 2 F^\mu_{\alpha p} \omega^\alpha, &\quad \theta^\mu_\alpha &= F^\mu_{\alpha a} \omega^a + F^\mu_{\alpha p}\omega^p 
\endaligned
\end{equation}
then comparison with~\eqref{eq:2:6}, using~\eqref{eq:4:7a} and~\eqref{eq:4:10}, gives
\begin{equation}\label{eq:4:14}
\aligned
F^\alpha_{pa} &= -h^\alpha_{pa} \sin s_1, &\quad   F^\mu_{pa} &= - h^\mu_{pa} \sin s_1 \\
F^\mu_{\alpha a} &= - h^\mu_{\alpha a} \sin s_1, &\quad F^\mu_{\alpha p} &= h^\mu_{\alpha p} \cos  s_1
\endaligned
\end{equation}
Notice that the distribution obtained by setting any three sets of $\{\omega^a\}$, $\{\omega^p\}$, $\{\omega^\alpha \}$ and $\{\omega^\mu\}$ equal to zero is completely integrable and its integral submanifolds are the respective curvature surfaces.

Equations~\eqref{eq:2:1} become, for the Darboux frame field~\eqref{eq:4:12a},
\begin{equation}\label{eq:4:15}
\aligned
d\x &= \omega^p e_p + \omega^\alpha e_\alpha + \omega^\mu e_\mu \\
de_0 &= \omega^a e_a - \omega^\alpha e_\alpha + \omega^\mu e_\mu \\
de_a &= -\omega^a e_0 + \theta^b_a e_b + \theta^q_a e_q + \theta^\a_a e_\alpha + \theta^\mu_a e_\mu  \\
de_p &= -\omega^p \x + \theta^b_p e_b + \theta^q_p e_q + \theta^\a_p e_\alpha + \theta^\mu_p e_\mu \\
de_\alpha &= -\omega^\alpha \x + \omega^\alpha e_0 +\theta^a_\alpha e_a + \theta^q_\alpha e_q + \theta^\b_\alpha e_\b + \theta^\mu_\alpha e_\mu \\
de_\mu &= -\omega^\mu \x - \omega^\mu e_0 + \theta^a_\mu e_a + \theta^q_\mu e_q + \theta^\a_\mu e_\alpha + \theta^\nu_\mu e_\nu
\endaligned
\end{equation}
The Cartan-M\"unzner polynomial $F:\R^{n+1} \to \R$ defining the isoparametric function $f = F|_{S^n}: S^n \to [-1,1]$ has $\pm1$ as the only two singular values, and focal points at a distance $\pi/2$ along a normal geodesic from each other lie on the same focal submanifold.  If our second order Darboux frame field~\eqref{eq:4:12a} is along the focal submanifold
\begin{equation*}
\x: U \subset M \to M_+ = f^{-1}\{1\} \subset S^n
\end{equation*}
then the tube~\eqref{eq:3:1} with $t = \pi/2$ shows that the image of $\bar x = e_0: U \to M_+$ is the same focal submanifold. If we let $\bar e_0 = \x$, then by~\eqref{eq:4:12}
\begin{equation}\label{eq:4:21}
\aligned
d\bar \x &= de_0 = \omega^a e_a - \omega^\alpha e_\alpha + \omega^\mu e_\mu  \\
d \bar e_0 &= d\x = \omega^p e_p + \omega^\alpha e_\alpha + \omega^\mu e_\mu 
\endaligned
\end{equation}
which shows that $e_a, e_\alpha, e_\mu$ are tangent to $M_+$ at $\bar \x = e_0$, while $\bar e_0, e_p$ are normal to $M_+$ at $\bar \x$.
The second fundamental form at $\bar \x$ with respect to $\bar e_0$ is
\begin{equation*}
\overline{II}_{\bar e_0} = -d\bar \x \cdot d\bar e_0 = -de_0 \cdot d\x = II_{e_0} = \sum \omega^\alpha \omega^\alpha - \sum \omega^\mu \omega^\mu 
\end{equation*}
which implies that $V_+$ is the $+1$ eigenspace and $V_-$ is the $-1$ eigenspace of $\overline{II}_{\bar e_0}$ at $\bar \x$.  Therefore, the principal curvature spaces of $\bar e_0$ at $\bar \x$ are
\begin{equation}\label{eq:4:22}
\bar V_+ = V_+, \quad \bar V_- = V_-, \quad \bar V_0 = \mbox{span}\{e_a\}
\end{equation}
It follows that a second order Darboux frame field along $\bar \x$ on $U$ is 
\begin{equation}\label{eq:4:20}
\bar \x = e_0,\,\, \bar e_0 = \x,\,\, \bar e_a = e_{a+m},\,\, \bar e_{a+m} = e_a, \bar e_\alpha = e_{\alpha},\,\, \bar e_\mu = e_\mu
\end{equation}
From~\eqref{eq:4:21} we see that
\begin{equation}\label{eq:4:23}
\aligned
\bar \omega^a &= \omega^{a+m}, & \bar \omega^{a+m} &= \omega^a, & \bar \omega^\alpha &= -\omega^\alpha, & \bar \omega^\mu &=  \omega^\mu
\endaligned
\end{equation}
is the coframe field dual to~\eqref{eq:4:20}.  

Of the forms in~\eqref{eq:4:13} for the frame field~\eqref{eq:4:20} and its coframe field~\eqref{eq:4:23}, we consider
\begin{equation*}
\aligned
d \bar e_\alpha \cdot \bar e_\mu &= \bar \theta^\mu_\alpha = \bar F^\mu_{\alpha a} \bar \omega^a + \bar F^\mu_{\alpha\, a+m} \bar \omega^{a+m} \\
= de_\alpha \cdot de_\mu &= \theta^\mu_\alpha = F^\mu_{\alpha a} \omega^a + F^\mu_{\alpha\,a+m} \omega^{a+m}
\endaligned
\end{equation*}
to conclude that
\begin{equation}\label{eq:4:26}
\bar F^\mu_{\alpha a} = F^\mu_{\alpha \,a+m}, \quad \bar F^\mu_{\alpha\,a+m} = F^\mu_{\alpha a} 
\end{equation}
Therefore, if $v = \sum ( x_\alpha  e_\alpha +  y_\mu  e_\mu) \in V_+ \oplus V_-$,
then
\begin{equation}\label{eq:4:27}
\bar p_a(v) = \sum_{\alpha, \beta} \bar F^\mu_{\alpha a} x_\alpha  y_\mu = \sum_{\alpha, \beta} F^\mu_{\alpha\, a+m}x_\alpha  y_\mu = p_{a+m}(v)
\end{equation}
where the polynomials $\bar p_a$ and $p_{a+m}$ are defined by these equations.

\section{Consequences of the structure equations}\label{section5}
We continue working with a second order frame field~\eqref{eq:4:12a} along the focal submanifold $\x$ defined in~\eqref{eq:3:1} with $t=s_1$.
Equations~\eqref{eq:4:14} show that differentiating equations~\eqref{eq:2:6} is equivalent to differentiating equations~\eqref{eq:4:13}, which we now proceed to do.  In preparation for this we first take the exterior differential of the coframe field~\eqref{eq:4:10} to obtain
\begin{equation}\label{eq:5:1}
\aligned
d\omega^a &= - \theta^a_b\w \omega^b - F^\alpha_{pa}\omega^p\w\omega^\alpha - F^\mu_{pa}\omega^p\w\omega^\mu -4 F^\mu_{\alpha a}\omega^\alpha \w\omega^\mu \\
d\omega^p &= - \theta^p_q\w\omega^q + F^\alpha_{pa}\omega^a\w\omega^\alpha + F^\mu_{pa}\omega^a\w\omega^\mu +4 F^\mu_{\alpha p}\omega^\alpha\w\omega^\mu \\
d\omega^\alpha &= - \theta^\alpha_\beta\w\omega^\beta - F^\alpha_{pa}\omega^a\w\omega^p + F^\mu_{\alpha a}\omega^a\w\omega^\mu - F^\mu_{\alpha p}\omega^p\w\omega^\mu \\
d\omega^\mu &= - \theta^\mu_\nu \w\omega^\nu - F^\mu_{pa}\omega^a\w\omega^p - F^\mu_{\alpha a}\omega^a\w\omega^\alpha + F^\mu_{\alpha p} \omega^p\w\omega^\alpha.
\endaligned
\end{equation}
We define the \textit{covariant derivatives} of the tensors $F^\alpha_{pa}$, $F^\mu_{pa}$, $F^\mu_{\alpha a}$ and $F^\mu_{\alpha p}$, respectively, to be the $1$-forms
\begin{equation}\label{eq:5:2}
\aligned
F^\alpha_{pai}\omega^i &= dF^\alpha_{pa} - F^\alpha_{qa}\theta^q_p - F^\alpha_{pb}\theta^b_a + F^\beta_{pa} \theta^\alpha_\beta \\
F^\mu_{pai} \omega^i &= dF^\mu_{pa} - F^\mu_{qa}\theta^q_p - F^\mu_{pb}\theta^b_a + F^\nu_{pa}\theta^\mu_\nu \\
F^\mu_{\alpha a i} \omega^i &= dF^\mu_{\alpha a} - F^\mu_{\beta a} \theta^\beta_\alpha - F^\mu_{\alpha b} \theta^b_a + F^\nu_{\alpha a} \theta^\mu_\nu \\
F^\mu_{\alpha p i} \omega^i &= dF^\mu_{\alpha p} - F^\mu_{\beta p}\theta^\beta_\alpha - F^\mu_{\alpha q} \theta^q_p + F^\nu_{\alpha p} \theta^\mu_\nu
\endaligned
\end{equation}
Any other second order frame field along $\x$ is given in terms of~\eqref{eq:4:12} by
\begin{equation}\label{eq:5:3}
\x, e_0, \hat e_a, \hat e_p, \hat e_\alpha, \hat e_\mu
\end{equation}
where
\begin{equation}\label{eq:5:4}
\hat e_a = A^b_a e_b, \,\, \hat e_p = A^q_p e_q, \,\, \hat e_\alpha = A^\beta_\alpha e_\beta, \,\, \hat e_\mu = A^\nu_\mu e_\nu
\end{equation}
with $(A^b_a), (A^q_p) : U \to O(m)$ and $(A^\beta_\alpha), (A^\nu_\mu): U \to O(N)$ smooth maps.  If the coefficients with respect to this new frame field are denoted by the same letters covered by a hat, then the transformation rules are tensorial.  For example,
\begin{equation}\label{eq:5:5} 
\hat F^\alpha_{pa} = A^\alpha_\beta F^\beta_{qb}A^q_p A^b_a, \,\, \hat F^\alpha_{pab} = A^\alpha_\beta F^\beta_{qcd}A^q_p A^c_a A^d_b
\end{equation}
and so forth.  If we take the exterior differential of the equations~\eqref{eq:4:13} and use~\eqref{eq:5:1} and~\eqref{eq:5:2} together with the Maurer-Cartan structure equations~\eqref{eq:2:1a} we obtain the following sets of equations (compare~\cite[I, p.\ 536 and II, p.\ 45]{OT}).
\begin{equation}\label{eq:5:6}
\aligned
F^\alpha_{pa}F^\alpha_{qb} + F^\alpha_{pb}F^\alpha_{qa} - (F^\mu_{pa}F^\mu_{qb} + F^\mu_{pb}F^\mu_{qa}) &= 0 \\
F^\alpha_{pa}F^\beta_{pb} + F^\alpha_{pb}F^\beta_{pa} + 2(F^\mu_{\alpha a}F^\mu_{\beta b} + F^\mu_{\alpha b}F^\mu_{\beta a}) &= \delta_{\alpha\beta}\delta_{ab} \\
F^\alpha_{pa}F^\beta_{qa} + F^\alpha_{qa}F^\beta_{pa} + 2 (F^\mu_{\alpha p}F^\mu_{\beta q} + F^\mu_{\alpha q}F^\mu_{\beta p}) &=  \delta_{pq}\delta_{\alpha\beta} \\
F^\mu_{pa}F^\nu_{pb} + F^\mu_{pb}F^\nu_{pa} + 2(F^\mu_{\alpha a}F^\nu_{\alpha b} + F^\mu_{\alpha b}F^\nu_{\alpha a}) &=  \delta_{ab}\delta_{\mu\nu} \\
F^\mu_{pa}F^\nu_{qa} + F^\mu_{qa}F^\nu_{pa} + 2(F^\mu_{\alpha p}F^\nu_{\alpha q} + F^\mu_{\alpha q}F^\nu_{\alpha p}) &=  \delta_{pq}\delta_{\mu\nu}  \\
F^\mu_{\alpha a}F^\nu_{\beta a} + F^\mu_{\beta a}F^\nu_{\alpha a} 
- (F^\mu_{\alpha p} F^\nu_{\beta p} + F^\mu_{\beta p}F^\nu_{\alpha p}) & = 0
\endaligned
\end{equation}

\begin{equation}\label{eq:5:7}
\aligned
F^\alpha_{pab} &=  -F^\mu_{pa}F^\mu_{\alpha b} - 2F^\mu_{pb}F^\mu_{\alpha a} \\
F^\alpha_{paq} &= F^\mu_{pa}F^\mu_{\alpha q} + 2F^\mu_{\alpha p}F^\mu_{qa} \\
F^\alpha_{pa\beta} &= 2F^\mu_{\alpha p}F^\mu_{\beta a} -2 F^\mu_{\beta p}F^\mu_{\alpha a} 
\endaligned
\end{equation}

\begin{equation}\label{eq:5:8}
\aligned
F^\mu_{pab} &=  F^\alpha_{pa}F^\mu_{\alpha b} + 2 F^\alpha_{pb}F^\mu_{\alpha a}\\
F^\mu_{paq} &= - F^\alpha_{pa}F^\mu_{\alpha q} -2 F^\mu_{\alpha p} F^\alpha_{qa} \\
F^\mu_{pa\nu} &= 2 F^\mu_{\alpha p}F^\nu_{\alpha a} - 2F^\mu_{\alpha a}F^\nu_{\alpha p}
\endaligned 
\end{equation}

\begin{equation}\label{eq:5:9}
\aligned
F^\mu_{\alpha ab} &= - \frac 12 F^\mu_{pa}F^\alpha_{pb} + \frac 12 F^\mu_{pb}F^\alpha_{pa} \\
F^\mu_{\alpha a\beta} &= F^\mu_{\alpha p} F^\beta_{pa} + 2 F^\mu_{\beta p}F^\alpha_{pa} \\
F^\mu_{\alpha a \nu} &= F^\mu_{\alpha p} F^\nu_{pa} + 2 F^\mu_{pa}F^\nu_{\alpha p} 
\endaligned 
\end{equation}

\begin{equation}\label{eq:5:10}
\aligned
F^\mu_{\alpha pq} &= \frac 12 F^\mu_{pa}F^\alpha_{qa} - \frac12 F^\mu_{qa}F^\alpha_{pa} \\
F^\mu_{\alpha p \beta} &= - F^\mu_{\alpha a} F^\beta_{pa} -2F^\mu_{\beta a}F^\alpha_{pa} \\
F^\mu_{\alpha p\nu} &= - F^\mu_{\alpha a}F^\nu_{pa} - 2F^\mu_{pa}F^\nu_{\alpha a}
\endaligned 
\end{equation}

\begin{equation}\label{eq:5:11}
F^\alpha_{pa\mu} = -F^\mu_{pa\alpha} = -2 F^\mu_{\alpha ap} = -2F^\mu_{\alpha pa}
\end{equation}

\section{Second fundamental forms of a focal submanifold}\label{section6}
Consider the focal submanifold $\x$ of~\eqref{eq:3:1} with $t=s_1$ with a second order frame field~\eqref{eq:4:12a} along it on $U$.  For each point of $\x$, Corollary~\ref{co:4:2} tells us the principal curvatures of 
the second fundamental forms $II_{e_a}$ of $\x$.  In order to derive the consequence of this knowledge, we begin by finding the expression of $II_{e_a}$ of $\x$ in terms of the orthonormal coframe field $\omega^p, \omega^\alpha, \omega^\mu$ and from that obtain the matrices of the corresponding shape operators with respect to the orthonormal tangent frame field $e_p, e_\alpha, e_\mu$.  For our frame, equations~\eqref{eq:2:1} have become, in part,
\begin{equation}\label{eq:6:1}
\aligned
d\x &= \omega^p e_p + \omega^\alpha e_\alpha + \omega^\mu e_\mu \\
de_a &= (k_1 e_0 - \x)\theta^a + \theta^b_a e_b + \theta^p_a e_p + \theta^\alpha_a e_\alpha + \theta^\mu_a e_\mu
\endaligned
\end{equation}
The shape operator $S_a$ is the symmetric operator on the tangent space at $\x$ given by
\begin{equation}\label{eq:6:2}
II_{e_a} = -de_a \cdot d\x =  d\x \circ S_a \cdot d\x
\end{equation}
That is, $S_a$ is the tangential component of $-de_a$.  Combining the second equation in~\eqref{eq:6:1} with~\eqref{eq:4:13}, we find
\begin{equation*}\label{eq:6:3}
S_a  = (2F^\mu_{\alpha a} e_\mu -F^\alpha_{pa}e_p)\omega^\alpha + (2F^\mu_{\alpha a}e_\alpha + F^\mu_{pa}e_p) \omega^\mu + (-F^\alpha_{pa}e_\alpha + F^\mu_{pa}e_\mu)\omega^p
\end{equation*}
Recall the curvature spaces $V_0, V_+, V_-$ defined in~\eqref{eq:6:4}.  Define linear operators
\begin{equation}\label{eq:6:5}
\aligned
A_a &= 2F^\mu_{\alpha a}e_\alpha \omega^\mu : V_- \to V_+ \\ 
B_a &= -F^\alpha_{pa}e_\alpha \omega^p: V_0 \to V_+ \\ 
C_a &= F^\mu_{pa}e_\mu \omega^p : V_0 \to V_-
\endaligned
\end{equation}
and their transposes
\begin{equation}\label{eq:6:6}
\aligned
\tose A_a &= 2F^\mu_{\alpha a}e_\mu \omega^\alpha : V_+ \to V_- \\ 
\tose B_a &= -F^\alpha_{pa}e_p \omega^\alpha: V_+ \to V_0 \\ 
\tose C_a &= F^\mu_{pa}e_p \omega^\mu : V_- \to V_0
\endaligned
\end{equation}
With respect to the orthogonal direct sum decomposition $V_+ \oplus V_- \oplus V_0$ of the tangent space to $\x$ at the point, the operator $S_a$ has the block form
\begin{equation}\label{eq:6:7}
S_a = \begin{pmatrix} 0 & A_a & B_a \\ \tose A_a & 0 & C_a \\ \tose B_a & \tose C_a & 0 \end{pmatrix}
\end{equation}
Restriction of the second fundamental forms $II_{e_0}$ and $II_{e_a}$
to $V_+ \oplus V_-$ defines quadratic forms
\begin{equation}\label{eq:6:polys}
\aligned
p_0(x,y) &= II_{e_0}((x,y),(x,y)) = \sum_\a x_\a^2 - \sum_\mu
y_\mu^2 \\
p_a(x,y) &= \frac14 II_{e_a}((x,y),(x,y)) = F^\mu_{\a a}x_\a y_\mu
\endaligned
\end{equation}
where $x = x_\a e_\a \in V_+$ and $y = y_\mu e_\mu \in V_-$.

By Corollary~\ref{co:4:2}, the minimal polynomial of $S_a$ is $x(x^2-1)$, and therefore $S_a^3 = S_a$ at every point of $U$.  Actually, something stronger is true.  If $t = (t^a) \in S^{m-1}$ is any unit vector, then $t^a e_a = \n$ is a normal vector to $\x$ at the point, and if its shape operator is denoted $S$, then $S^3 = S$ as well.  On the other hand,
\begin{equation}\label{eq:6:8}
S = t^a S_a
\end{equation}
so that $S^3=S$ implies that
\begin{equation}\label{eq:6:9}
\aligned
\sum_a t^aS_a &= \sum_a (t^a)^3 S_a^3 + \sum_{a\neq b} (t^a)^2 t^b(S_aS_aS_b + S_aS_bS_a + S_bS_aS_a) \\
&\quad + \sum_{a\neq b \neq c} t^at^bt^c(\sum_{\sigma \in \mathcal S_3} S_{\sigma(a)} S_{\sigma(b)} S_{\sigma(c)})
\endaligned
\end{equation}
where $\mathcal S_3$ is the symmetric group on three letters.
Multiplying $\sum (t^a)^2 = 1$ by $t^b$ gives 
\begin{equation}\label{eq:6:10}
t^b = (t^b)^3 + \sum_{a\neq b} t^b(t^a)^2
\end{equation}
which substituted into~\eqref{eq:6:9} gives
\begin{equation*}\label{eq:6:11}
\aligned
0 &= \sum_a (S_a^3-S_a)(t^a)^3 + \sum_{a\neq b}(t^a)^2 t^b(S_aS_aS_b + S_aS_bS_a + S_bS_aS_a - S_b) \\
&\quad + \sum_{a\neq b \neq c} t^at^bt^c \sum_{\sigma \in \mathcal S_3} S_{\sigma(a)} S_{\sigma(b)} S_{\sigma(c)}
\endaligned
\end{equation*}
A homogeneous polynomial on $\R^m$ which is identically zero on the unit sphere $S^{m-1}\subset \R^m$ must be identically zero, that is, its coefficients must all be zero.  Therefore (compare~\cite[I, p.\ 534]{OT})
\begin{equation}\label{eq:6:12}
\aligned
S_a &= S_a^3, \mbox{ for all $a$} \\
S_b &= S_aS_aS_b + S_aS_bS_a + S_bS_aS_a, \mbox{ for all $a \neq b$} \\
0 &= \sum_{\sigma \in \mathcal S_3} S_{\sigma(a)} S_{\sigma(b)} S_{\sigma(c)}, \mbox{ for all $a \neq b \neq c$}
\endaligned
\end{equation}

\begin{proposition}\label{pr:8:4}  If $m<N$, then the operators $A_a$ in~\eqref{eq:6:5} must be linearly independent on an open dense subset of $U$.
\end{proposition}

\begin{proof}  Linear independence being an open condition implies that these operators are linearly independent on an open subset of $U$ (possibly empty).  It remains to show that they are linearly independent on a dense subset of $U$.  Suppose, to the contrary, that the operators $A_a$ are linearly dependent at every point of some open subset $V\subset U$.  This means that there exists a smooth unit vector
\begin{equation*}
u = (u^a): V \to \R^m
\end{equation*}
such that
\begin{equation} \label{eq:8:39a}
u^a F^\mu_{\alpha a} = 0
\end{equation}
for all $\mu$ and $\alpha$, at every point of $V$.  Then multiplying the second equation in~\eqref{eq:5:6} by $u^au^b$, summing on $a$ and $b$ and using~\eqref{eq:8:39a} gives
\begin{equation*}  
F^\alpha_{pa}u^aF^\beta_{pb}u^b = \delta_{\alpha \beta}
\end{equation*}
Therefore, 
\begin{equation*}
\{ F^\alpha_{pa}u^a e_p : \alpha = 2m+1, 2m+2,\dots, 2m+N\}
\end{equation*}
is an orthonormal set of $N$ vectors in the $m$-dimensional subspace $V_0$ defined in~\eqref{eq:6:4}, which contradicts the assumption that $m < N$.  We conclude that every open subset of $U$ contains a point at which the $A_a$ are linearly independent.  That is, the  $A_a$ are linearly independent on a dense subset of $U$.
\end{proof}
We need a condition which is stronger than the linear independence of the $A_a$.

\begin{definition}[Spanning Property]\label{spanning} The focal submanifold $\x$ satisfies the \textit{spanning property} at a point of $M$ if

(a).  There exists a vector $X = x_\alpha e_\alpha \in V_+$ such that the set of vectors $\{F^\mu_{\alpha a} x_\alpha e_\mu: a = 1,\dots,m\}$ in $V_-$ are linearly independent; and

(b).  There exists a vector $Y = y_\mu e_\mu \in V_-$ such that the set of vectors $\{F^\mu_{\alpha a} y_\mu e_\alpha : a = 1,\dots , m\}$ in $V_+$ are linearly independent.
\end{definition}

\begin{remark}\label{re:6:8} Observe that (a) is equivalent to

(a').  There exists $X =x_\alpha e_\alpha\in V_+$ such that
$\{F^\mu_{\alpha a} x_\alpha y_\mu e_a : Y = y_\mu e_\mu \in V_-\} =
\mbox{span}\{e_1,\dots, e_m\}$. 

and (b) is equivalent to

(b').  There exists $Y = y_\mu e_\mu\in V_-$ such that
$\{F^\mu_{\alpha a} x_\alpha y_\mu e_a : X = x_\alpha e_\alpha \in V_+
\} = \mbox{span}\{e_1,\dots,e_m \}$. 
\end{remark}

\begin{remark}\label{re:6:9}  If $\x$ satisfies the spanning property
at a point of $M$, then it satisfies it on some open neighborhood of
the point.  In fact, let $X = x_\alpha e_\alpha$ and $Y = y_\mu e_\mu$
be continuous vector fields whose values at the point give the vectors
in (a) and (b), respectively.  The vectors $F^\mu_{\alpha a} x_\alpha
e_\mu$, $a = 1,\dots,m$ are linearly independent at the point implies
that some $m\times m$ minor of the $N \times m$ matrix $(F^\mu_{\alpha
a} x_\alpha)$ is nonzero at the point.  Being a continuous function in
$M$, this minor remains nonzero on some open neighborhood of the
point, and thus (a) holds on this neighborhood.  A similar argument
shows that (b) holds on some neighborhood of the point, and thus (a)
and (b) hold on the intersection. 
\end{remark}

Let $\x, e_0, e_a, e_p, e_\alpha, e_\mu$ be a second order frame field~\eqref{eq:4:12a} along $\x$ on $U$, where $\x(U) \subset M_+$ is a focal submanifold.  Let the same letters with bars denote the second order frame field~\eqref{eq:4:20} along $\bar \x = e_0$ on $U$.  At each point of $U$ define bihomogeneous polynomials $p_a$ and $\bar p_a$ in $\R[x_\alpha, y_\mu]$ by
\begin{equation}\label{eq:6:13}
p_a(x,y) = F^\mu_{\alpha a} x_\alpha y_\mu, \quad \bar p_a(x,y) = \bar F^\mu_{\alpha a} x_\alpha y_\mu
\end{equation}
where $F^\mu_{\alpha a}$ and $\bar F^\mu_{\alpha a}$ are defined in~\eqref{eq:4:13} for the respective frame fields.

\begin{proposition}\label{pr:6:10} If at each point of $U$ there exist polynomials $f_{ab}$ in the polynomial ring $\R[x_\alpha, y_\mu]$ such that
\begin{equation}\label{eq:6:14}
\bar p_a = f_{ab} p_b
\end{equation}
and if the spanning condition holds for $\x$ on $U$, then there exists a second order frame field $\x, e_0, \hat e_a, \hat e_p, \hat e_\alpha, \hat e_\mu$ along $\x$ on $U$ with respect to which
\begin{equation}\label{eq:6:15}
\hat F^\mu_{\alpha\, a+m} = \hat F^\mu_{\alpha a}
\end{equation}
for all $a,\alpha,\mu$, at each point of $U$.
\end{proposition}

\begin{proof}  If we let
$p_{a+m}(x,y) = F^\mu_{\alpha\, a+m} x_\alpha y_\mu$, then by~\eqref{eq:4:27}, $p_{a+m} = \bar p_a$ and therefore~\eqref{eq:6:14} implies that at each point of $U$
\begin{equation}\label{eq:6:16}
p_{a+m} = f_{ab} p_b
\end{equation}
If we expand the right side of this equation in terms of the bi\-homo\-ge\-neous components of the $f_{ab}$ and collect all terms of the same bi-degrees, then all terms must cancel except those of bi-degree $(1,1)$, since $p_{a+m}$ has bi-degree $(1,1)$.  This results in an expression for $p_{a+m}$ as a linear combination of the $p_b$ with constant coefficients, since each $p_b$ has bi-degree $(1,1)$.  Hence, we may assume that
the $f_{ab}$ in~\eqref{eq:6:16} are constant polynomials.
Now~\eqref{eq:6:16} implies that
\begin{equation}\label{eq:6:17}
F^\mu_{\alpha\, a+m} = f_{ab} F^\mu_{\alpha b}
\end{equation}
for all $\alpha, \mu$ at each point of $U$.  We claim that the functions $f_{ab}: U \to \R$ are smooth.  In fact, if we let $A_{a+m} = 2F^\mu_{\alpha\,a+m} e_\alpha \omega^\mu:V_- \to V_+$ and let $A_a$ be the operators defined in~\eqref{eq:6:5}, then~\eqref{eq:6:17} implies that $A_{a+m} = f_{ab}A_b$.  The spanning property implies that the operators $A_b$ are linearly independent in $\text{End}(V_-,V_+)$, and therefore at each point of $U$ an inner product can be defined on this space of endomorphisms, depending smoothly on the point of $U$, such that $\{A_b\}$ is an orthonormal set.  Then $f_{ab} = \langle A_{a+m}, A_b\rangle:U \to \R$ is smooth.

Fix $\alpha = \alpha_0$ and for each $\mu$ define vectors in $\R^m$
\begin{equation*}
W_\mu = \begin{pmatrix} F^\mu_{\alpha_0 1} \\ \vdots \\ F^\mu_{\alpha_0 m} \end{pmatrix}, \quad 
V_\mu = \begin{pmatrix} F^\mu_{\alpha_0\, m+1} \\ \vdots \\ F^\mu_{\alpha_0 \, m+m} \end{pmatrix}
\end{equation*}
If we define the $m\times m$ matrix $B = (f_{ab})$, then by~\eqref{eq:6:17}, we have 
\begin{equation}\label{eq:6:18}
V_\mu = BW_\mu
\end{equation}
for each $\mu$.
The sixth equation in~\eqref{eq:5:6} says that for any $\mu$ and $\nu$
\begin{equation*}
V_\mu \cdot V_\nu = W_\mu \cdot W_\nu
\end{equation*}
Combining these equations, we have
\begin{equation}\label{eq:6:19}
W_\mu \cdot W_\nu = BW_\mu \cdot BW_\nu
\end{equation}
for all $\mu, \nu$.  It follows that $B$ is orthogonal, provided that the set $\{W_\mu\}$ spans $\R^m$.  By the spanning property, this is true for some choice of $\alpha_0$.  Therefore, assuming we have made that choice, we have a smooth map
\begin{equation*}
B = (f_{ab}):U \to O(m)
\end{equation*}
Alter the second order frame field along $\x$ by
\begin{equation*}
\hat e_{a+m} = e_{b+m} f_{ba}
\end{equation*}
leaving the other vectors in the frame unchanged.  If we let $\hat F^\mu_{\alpha a}$, etc. be the coefficients with respect to this new frame field, then by~\eqref{eq:5:5}, we have $\hat F^\mu_{\alpha a} = F^\mu_{\alpha a}$ and, also using~\eqref{eq:6:17}, we have
\begin{equation*}
\hat F^\mu_{\alpha\, a+m} = F^\mu_{\alpha\, b+m} f_{ba}  =  f_{bc} F^\mu_{\alpha c}f_{ba} = \hat F^\mu_{\alpha c} \delta_{ca} = \hat F^\mu_{\alpha a}
\end{equation*}
which proves~\eqref{eq:6:15}.
\end{proof}

\section{The Ferus-Karcher-M\"unzner construction}\label{section7}
Let $P_0,P_1,\dots,P_m$ be a Clifford system on $\R^{2l}$.  Recall that this means that these are symmetric operators on $\R^{2l}$ satisfying 
\begin{equation}\label{eq:7:1}
P_iP_j+P_jP_i = 2\delta_{ij}I, \quad i,j = 0,1,\dots,m
\end{equation}
It follows that each operator $P_i$ is also orthogonal.
For this section we modify the index conventions~\eqref{eq:4:convention} by
\begin{equation}\label{eq:7:2}
i,j,k \in \{0,\dots,m\}
\end{equation}
and now $N = l-m-1$ and $n+1 = 2l$.
If $A \in SO(m+1)$, and if we let 
\begin{equation}\label{eq:7:3}
Q_i = A^j_i P_j
\end{equation}
then $Q_0, Q_1,\dots, Q_m$ is also a Clifford system on $\R^{2l}$.  Since $Q_0^2 = I$, the eigenvalues of $Q_0$ must be $\pm1$.  If $E_\pm$ are the eigenspaces of $Q_0$, then $\R^{2l} = E_+ \oplus E_-$ is an orthogonal direct sum and $E_\pm$ each has dimension $l$, because for any $a$, the operator $Q_a$ interchanges $E_+$ and $E_-$.

Because $P_0,\dots,P_m$ are linearly independent, 
\begin{equation}\label{eq:7:4}
M_+ = \{\x\in S^{2l-1} \subset \R^{2l} : P_i\x\cdot \x = 0, \quad i = 0,\dots,m \} 
\end{equation}
is a submanifold of $S^{2l-1}$ of codimension $m+1$.
If $\x \in M_+$, then $Q_0\x,\dots,Q_m\x$ is an orthonormal set of unit normal vectors to $M_+$ in $S^{2l-1}$.  Therefore, this is a global frame field for the normal bundle of $M_+$ and the unit normal bundle of $M_+$ is isomorphic to the trivial bundle 
\begin{equation}\label{eq:7:5}
M = M_+ \times S^m
\end{equation}
Consider the principal bundle
\begin{equation}\label{eq:7:6}
\aligned
SO(m+1) &\to S^m \\
A &\mapsto A_0
\endaligned
\end{equation}
where for any $A\in SO(m+1)$ we let $A_i$ denote the $i^{\mbox{th}}$ column of $A$.
For a section $A$ of~\eqref{eq:7:6}, denote its pull-back to $S^m$ of the Maurer-Cartan form of $SO(m+1)$ by
\begin{equation}\label{eq:7:7}
A^{-1}dA = \tau = (\tau^i_j)
\end{equation}
an $\mathfrak o(m+1)$-valued form on $S^m$.  Then $dA_i = A_j \nu^j_i$, and thus, for the Clifford systems
\begin{equation}\label{eq:7:8}
Q_i = A^j_i P_j
\end{equation}
depending on $A\in SO(m+1)$, we have
\begin{equation}\label{eq:7:9}
dQ_i = Q_j\tau^j_i
\end{equation}
for each $i$.  Observe that $\tau^1_0,\dots,\tau^m_0$ is a local coframe field in $S^m$. 
For each $(x,A_0) \in M=M_+\times S^m$, there is an orthogonal direct sum
\begin{equation}\label{eq:7:10}
\R^{2l} = \text{span}\{\x \} \oplus M_+^\perp(\x) \oplus T_0(\x,A_0) \oplus T_1(\x,A_0) \oplus T_{-1}(\x,A_0),
\end{equation}
where
\begin{equation}\label{eq:7:11}
\aligned
M_+^\perp(\x) &= \mbox{span}\{Q_0\x,\dots,Q_m\x\} = \mbox{span}\{P_0\x,\dots,P_m\x\} \\
T_0(\x,A_0) &= \mbox{span}\{Q_aQ_0\x : \mbox{ for all $a$}\} \\
T_+(\x,A_0) &= E_- \cap T_\x M_+ = \{ X \in E_- : X\cdot Q_i\x = 0 \mbox{ for all $i$}\} \\
 &= \{X\in \R^{2l} : Q_0X= -X \mbox{ and } X\cdot P_i\x=0, \mbox{ for all $i$} \} \\
T_-(\x,A_0) &= E_+ \cap T_\x M_+ = \{X\in E_+ : X\cdot Q_i \x =0, \mbox{ for all $i$} \} \\
 &= \{X \in \R^{2l} : Q_0 X = X \mbox{ and } X\cdot P_i\x=0, \mbox{for all $i$} \}
\endaligned
\end{equation}
Then $\mbox{dim}\,M_+^\perp(\x)=m+1$, $\mbox{dim}\,T_0(\x,A_0)=m$, $\mbox{dim}\,T_+(\x,A_0) = N$ and $\mbox{dim}\,T_-(\x,A_0) = N$, where $N = l-(m+1)$.  Notice that
\begin{equation}\label{eq:7:12}
Q_0 : T_0(\x,A_0) \to M_+^\perp(\x)
\end{equation}
because $Q_0Q_aQ_0\x = -Q_a\x \in M_+^\perp$, for any $a$.

For any point in $M = M_+ \times S^m$, there is an open neighborhood about it of the form $U\times V$, where $U\subset M_+$ and $V\subset S^m$, such that the section $A$ of~\eqref{eq:7:6} is defined on $V$ and such that there exist smooth orthonormal bases $e_\alpha$ of $T_+(\x,A_0)$ and $e_\mu$ of $T_-(\x,A_0)$ on $U\times V$.  This means that at each point of $U\times V$
\begin{equation}\label{eq:7:13}
\aligned
Q_0 e_\alpha &= -e_\alpha \mbox{ and } e_\alpha \cdot Q_i\x =0 \\
Q_0 e_\mu &= e_\mu \mbox{ and } e_\mu \cdot Q_i\x =0
\endaligned
\end{equation}
Compose $\x:M_+ \to S^{2l-1}$ with the projection $M = M_+\times S^m \to M_+$ so that we may regard it as a mapping $\x: M\to S^{2l-1}$.  Then
\begin{equation}\label{eq:7:14}
\x,\,\, e_i = Q_i\x, \,\, e_p = Q_{p-m}Q_0\x, \,\, e_\alpha, \,\, e_\mu
\end{equation}
is a Darboux frame field along $\x$ on $U\times V$, where the $e_i$ are normal vectors and the rest are tangent to $\x$.

\begin{lemma}\label{le:7:A}  For any $\x \in M_+$
\begin{equation}\label{eq:7:14a}
Q_iQ_jQ_k\x \cdot \x =0
\end{equation}
for all $i,j,k$ and
\begin{equation}\label{eq:7:14aa}
L^a_{bc} = Q_aQ_bQ_c e_0 \cdot \x
\end{equation}
is skew-symmetric in $a,b,c$.
\end{lemma}

\begin{proof}  If $i,j,k$ are distinct, then
\begin{equation*}
Q_iQ_jQ_k \x \cdot \x = \x \cdot Q_kQ_jQ_i \x = -\x \cdot Q_iQ_jQ_k \x
\end{equation*}
which implies~\eqref{eq:7:14a}.  If the indices are not distinct, then the product is a single $\pm Q_i$ and $Q_i\x\cdot \x = 0$ by definition of $M_+$.

If any two of $a,b,c$ are the same, then the product $Q_aQ_bQ_c$ is a single operator $\pm Q_a$, for some $a$, and we know that $Q_ae_0 \cdot \x =0$.  If $a,b,c$ are distinct, then $Q_aQ_bQ_c$ changes sign if any two indices are switched.  Therefore, $L^a_{bc}$ is skew-symmetric in $a,b,c$.
\end{proof}

\begin{lemma}\label{le:7:B}  For the Darboux frame field~\eqref{eq:7:14} along $\x$,
\begin{align}
Q_i\x \cdot \x &= 0, \text{ for all $i$} \label{eq:7:B1} \\
Q_ie_j \cdot e_k &= 0, \text{ for all $i,j,k$} \label{eq:7:B2} \\
Q_i e_p \cdot e_q &= 0, \text{ for all $i,p,q$} \label{eq:7:B3} \\
Q_a e_\alpha \cdot e_\beta &= 0, \text{ for all $a,\alpha,\beta$} \label{eq:7:B4} \\
Q_a e_\mu \cdot e_\nu &= 0, \text{ for all $a, \mu, \nu$} \label{eq:7:B5}
\end{align}
at each point of $U\times V$.
\end{lemma}

\begin{proof} The first equation follows from the definition of $M_+$.  For the second equation
\begin{equation*}
Q_i e_j \cdot e_k = Q_iQ_j \x \cdot Q_k \x = Q_kQ_iQ_j \x \cdot \x =0
\end{equation*}
by Lemma~\ref{le:7:A}.  For the third equation
\begin{equation*}
Q_ae_p\cdot e_q = Q_aQ_{p-m}Q_0\x \cdot Q_{q-m}Q_0 \x =-Q_{q-m}Q_aQ_{p-m}\x \cdot \x =0
\end{equation*}
by Lemma~\ref{le:7:A} and
\begin{equation*}
Q_0 e_p \cdot e_q = Q_0Q_{p-m}Q_0\x \cdot Q_{q-m}Q_0 \x = -\x \cdot Q_{p-m}Q_{q-m} Q_0 \x =0
\end{equation*}
by Lemma~\ref{le:7:A}.  Equations $4$ and $5$ follow from the observation made above that $Q_a$ interchanges $E_-$ and $E_+$.
\end{proof}

\begin{lemma}\label{le:7:1}  For the frame field~\eqref{eq:7:14},
\begin{align}
d\x &= \omega^p e_p + \omega^\alpha e_\alpha + \omega^\mu e_\mu \label{eq:7:15} \\
de_0 &= \omega^a e_a -\omega^\alpha e_\alpha + \omega^\mu e_\mu \label{eq:7:15a}
\end{align}
where $\omega^p, \omega^\alpha, \omega^\mu$ are linearly independent one forms on $U$ with coefficients being functions on $U\times V$, and
\begin{equation}\label{eq:7:16}
\omega^a = \tau^a_0 - \omega^{a+m}
\end{equation}
A smooth coframe field on $U\times V$ is given by $\omega^a, \omega^p, \omega^\alpha, \omega^\mu$.
\end{lemma}

\begin{proof}  The expression~\eqref{eq:7:15} for $d\x$ follows from the fact that $\x:U \to \R^{2l}$ is an immersion and then $\omega^A = d\x \cdot e_A$, for $A = m+1,\dots,2l-1$.  Combining this with~\eqref{eq:7:9}, we have
\begin{equation}\label{eq:7:17}
\aligned
de_0 &= dQ_0\,\x + Q_0d\x \\
&= \tau^a_0 Q_a\x + \omega^{a+m}Q_0Q_aQ_0\x + \omega^\alpha Q_0e_\alpha + \omega^\mu Q_0 e_\mu \\
&= (\tau^a_0 - \omega^{a+m}) e_a - \omega^\alpha e_\alpha + \omega^\mu e_\mu
\endaligned
\end{equation}
which proves~\eqref{eq:7:15a}.
\end{proof}

For $t \in \R$, the tube of radius $t$ about $M_+$ is given by the immersion
\begin{equation}\label{eq:7:18}
\tilde \x : M \to S^{2l-1}, \quad
\tilde \x = \cos t\, \x + \sin t\, e_0 
\end{equation}
A unit normal vector field along $\tilde \x$ is
\begin{equation}\label{eq:7:19}
\tilde e_0 = -\sin t\, \x + \cos t\, e_0
\end{equation}
and a Darboux frame field along $\tilde \x$ is given by 
\begin{equation}\label{eq:7:19a}
\tilde\x,e_a,e_p,e_\alpha,e_\mu, \tilde e_0
\end{equation}
From~\eqref{eq:7:15} we compute
\begin{equation}\label{eq:7:20}
\aligned
d\tilde \x &= \sin t\, \omega^a e_a + \cos t\, \omega^p e_p \\
&\quad + (\cos t - \sin t) \omega^\alpha e_\alpha + (\cos t + \sin t )\omega^\mu e_\mu \\
d\tilde e_0 &= \cos t \, \omega^a e_a -\sin t\, \omega^p e_p \\
&\quad - (\cos t + \sin t) \omega^\alpha e_\alpha + (\cos t - \sin t ) \omega^\mu e_\mu
\endaligned
\end{equation}
which shows that
\begin{equation}\label{eq:7:21}
\aligned
\theta^a &= \sin t\, \omega^a, &\quad \theta^p &= \cos t\, \omega^p \\ 
\theta^\alpha &= (\cos t - \sin t) \omega^\alpha, &\quad  \theta^\mu &= (\cos t + \sin t) \omega^\mu
\endaligned
\end{equation}
is an orthonormal coframe field in $M$ for the metric $d\tilde \x \cdot d\tilde \x$ induced by $\tilde \x$.  The second fundamental form of $\tilde \x$ is then
\begin{equation*}\label{eq:7:22}
\aligned
&II_{\tilde e_0} = -d\tilde \x \cdot d\tilde e_0 \\
&= -\cot t\,\theta^a \theta^a + \tan t\,  \theta^p  \theta^p + \frac{\cot t + 1}{\cot t -1}  \theta^\alpha  \theta^\alpha - \frac{\cot t -1}{\cot t +1}  \theta^\mu  \theta^\mu \\
&= \cot(-t)\theta^a\theta^a + \cot(\frac \pi 2 - t) \theta^p \theta^p + \cot(\frac \pi 4 -t) \theta^\alpha \theta^\alpha + \cot(\frac{3\pi}4 -t) \theta^\mu \theta^\mu
\endaligned
\end{equation*}
from which we conclude that the principal curvatures are the constants $\cot(-t)$ and $\cot(\pi/2 - t)$, each with multiplicity $m$ and the constants $\cot(\pi/4 - t)$ and $\cot(3\pi/4 - t)$, each with multiplicity $N$.  In addition, the Darboux frame field~\eqref{eq:7:19a} along $\tilde \x$ is of second order.  Therefore, the $\tilde \x$ for $t\in \R$ is an isoparametric family of hypersurfaces in $S^{2l-1}$ and $\x$ is a focal submanifold.  This is the \textit{Ferus-Karcher-M\"unzner construction}, (FKM construction)~\cite{FKM}, of an isoparametric hypersurface from a given Clifford system.

We next calculate equations~\eqref{eq:4:13} for the FKM construction for a given Clifford system.
\begin{lemma}\label{le:7:2}  For the Darboux frame field~\eqref{eq:7:14} along $\x$, the coefficients of the forms $\theta^A_B = de_A \cdot e_B$ in~\eqref{eq:4:13} are given by
\begin{equation}\label{eq:7:23}
\aligned
F^\alpha_{pa} &= Q_{p-m}Q_a \x \cdot e_\alpha, &\quad F^\mu_{pa} &= Q_{p-m} Q_a \x \cdot e_\mu \\
F^\mu_{\alpha a} &= -\frac12 Q_a e_\mu \cdot e_\alpha, &\quad F^\mu_{\alpha p} &= -\frac12 Q_{p-m} e_\mu \cdot e_\alpha
\endaligned
\end{equation}
\end{lemma}

\begin{proof}  These coefficients are determined by $\theta^p_a$, $\theta^\alpha_a$ and $\theta^\alpha_p$.  From~\eqref{eq:7:9} and~\eqref{eq:7:15} we have
\begin{equation}\label{eq:7:24}
\aligned
de_a &= dQ_a \,\x + Q_a d\x \\
&=-\tau^a_0 e_0 + \tau^b_a e_b + \omega^{b+m} Q_a e_{b+m} + \omega^\alpha Q_a e_\alpha + \omega^\mu Q_a e_\mu
\endaligned
\end{equation}
and from~\eqref{eq:7:9} and~\eqref{eq:7:15a} we have
\begin{equation}\label{eq:7:24a}
\aligned
de_{a+m} &=  dQ_a\,e_0 + Q_a de_0 \\
&= -\tau^a_0 \x + \tau^b_a e_{b+m} + \omega^b Q_ae_b - \omega^\alpha Q_a e_\alpha + \omega^\mu Q_a e_\mu
\endaligned
\end{equation}
Using Lemma~\ref{le:7:B} and~\eqref{eq:4:13} we have
\begin{equation}\label{eq:7:25}
\aligned
&F^\alpha_{b+m\, a} \omega^\alpha - F^\mu_{b+m\, a} \omega^\mu = \theta^{b+m}_a = de_a \cdot e_{b+m} \\
&=  \omega^\alpha Q_a e_\alpha \cdot e_{b+m} + \omega^\mu Q_a e_\mu \cdot e_{b+m}
\endaligned
\end{equation}
which implies that
\begin{equation*}
\aligned
F^\alpha_{b+m\, a} &= Q_a e_\alpha \cdot e_{b+m} = Q_ae_\alpha \cdot Q_bQ_0 \x \\
&= Q_0 e_\alpha \cdot Q_aQ_b\x = -e_\alpha \cdot Q_aQ_b\x = Q_bQ_a \x \cdot e_\alpha
\endaligned
\end{equation*}
which is the first formula in~\eqref{eq:7:23}, and similarly
\begin{equation*}
-F^\mu_{b+m\, a} = Q_a e_\mu \cdot Q_bQ_0 \x = Q_a e_\mu \cdot Q_b\x = e_\mu \cdot Q_aQ_b \x
\end{equation*}
which gives the second formula in~\eqref{eq:7:23}.
In the same way,
\begin{equation}\label{eq:7:28}
\aligned
&F^\alpha_{b+m\,a} \omega^{b+m} -2 F^\mu_{\alpha a} \omega^\mu =\theta^\alpha_a = de_a \cdot e_\alpha \\
&= \omega^{b+m} Q_ae_{b+m} \cdot e_\alpha + \omega^\mu Q_a e_\mu \cdot e_\alpha
\endaligned
\end{equation}
which implies that  $F^\mu_{\alpha a} = Q_ae_\mu \cdot e_\alpha$, which is the third formula in~\eqref{eq:7:23}.
Next,
\begin{equation}\label{eq:7:29}
\aligned
&F^\alpha_{a+m\, b} \omega^b -2 F^\mu_{\alpha a+m} \omega^\mu = \theta^\alpha_{a+m} = de_{a+m} \cdot e_\alpha \\
&=  \omega^b Q_ae_b \cdot e_\alpha + \omega^\mu Q_a e_\mu \cdot e_\alpha
\endaligned
\end{equation}
which implies that $-2F^\mu_{\alpha a+m} = Q_a e_\mu \cdot e_\alpha$, which is the fourth formula in~\eqref{eq:7:23}.
\end{proof}

\begin{corollary}\label{co:7:1}
With respect to a Darboux frame~\eqref{eq:7:14} along an FKM construction $\x: M \to S^{2l-1}$, the coefficients~\eqref{eq:7:23} satisfy the equations
\begin{equation}\label{eq:7:30}
\aligned
F^\mu_{\alpha\, a+m} &= F^\mu_{\alpha a} \\
F^\alpha_{a+m\, b} &= -F^\alpha_{b+m\, a} \\
F^\mu_{a+m\,b} &= - F^\mu_{b+m\, a}
\endaligned
\end{equation}
\end{corollary}

\begin{proof}  From~\eqref{eq:7:23}, 
\begin{equation}\label{eq:7:31}
\aligned
F^\mu_{\alpha\,a+m} &= -\frac12 Q_a e_\mu \cdot e_\alpha = F^\mu_{\alpha a} \\
F^\alpha_{a+m\, b} + F^\alpha_{b+m\,a} &= (Q_aQ_b+ Q_bQ_a)\x \cdot e_\alpha = 0 \\
F^\mu_{a+m\, b} + F^\mu_{b+m\,a} &= (Q_aQ_b + Q_bQ_a)\x \cdot e_\mu = 0
\endaligned
\end{equation}
\end{proof}

\begin{proposition}\label{pr:7:1} For the Darboux frame field~\eqref{eq:7:14}, at any point of $U\times V \subset M$, the operators $Q_0, Q_a$ are given by
\begin{equation}\label{eq:7:32}
\aligned
Q_0 \x &= e_0 &\quad Q_0 e_0 &= \x &\quad Q_0 e_a &= -e_{a+m} \\
Q_0 e_{a+m} &= -e_a &\quad Q_0 e_\alpha &= -e_\alpha &\quad Q_0 e_\mu &= e_\mu
\endaligned
\end{equation}
and for each $a$
\begin{equation}\label{eq:7:33}
\aligned
Q_a\x &= e_{a} \\
Q_a e_0 &= e_{a+m} \\
Q_a e_b &= \delta_{ab}\x - L^c_{ab}e_{c+m} + F^\alpha_{a+m\,b}e_\alpha + F^\mu_{a+m\,b}e_\mu \\
Q_a e_{b+m} &= \delta_{ab}e_0 + L^c_{ab}e_c + F^\alpha_{b+m\,a}e_\alpha - F^\mu_{b+m\,a}e_\mu \\
Q_a e_\alpha &= F^\alpha_{a+m\, b}e_b + F^\alpha_{b+m\,a}e_{b+m} -2 F^\mu_{\alpha\, a}e_\mu \\
Q_a e_\mu &= F^\mu_{a+m\, b}e_b - F^\mu_{b+m\, a} e_{b+m} -2 F^\mu_{\alpha\, a}e_\alpha 
\endaligned
\end{equation}
where the coefficients are defined in~\eqref{eq:7:14aa} and~\eqref{eq:7:23}.
\end{proposition}

\begin{proof} The expansion~\eqref{eq:7:32} of $Q_0$ can be verified by inspection.  Also easy are the calculations $Q_b \x = e_b$ and $Q_b e_0 = Q_bQ_0 \x = e_{b+m}$.  To calculate $Q_b$ on the remaining basis vectors, we use the fact that the basis is orthonormal.  In the following calculations we use~\eqref{eq:7:1}, \eqref{eq:7:14a}, \eqref{eq:7:14}, \eqref{eq:7:14aa} and~\eqref{eq:7:23}.
\begin{align*}
Q_b e_a\cdot \x &= Q_bQ_a \x \cdot \x = \delta_{ab} \\
Q_b e_a \cdot e_0 &= Q_b Q_a \x \cdot Q_0 \x = Q_0 Q_b Q_a \x \cdot \x =0 \\
Q_b e_a \cdot e_c &= Q_b Q_a \x \cdot Q_c \x = Q_cQ_bQ_a \x \cdot \x =0 \\
Q_b e_a \cdot e_{c+m} &= Q_b Q_a \x \cdot Q_c Q_0 \x = Q_aQ_bQ_c Q_0 \x \cdot \x = L^a_{bc} \\
Q_b e_a \cdot e_\alpha &= Q_b Q_a \x \cdot e_\alpha = F^\alpha_{b+m\, a} \\
Q_b e_a \cdot e_\mu &= Q_b Q_a \x \cdot e_\mu = F^\mu_{b+m\, a}
\end{align*}
give the expansion of $Q_b e_a$.
\begin{align*}
Q_b e_{c+m} \cdot \x &= Q_b Q_c Q_0 \x \cdot \x =0 \\
Q_b e_{c+m} \cdot e_0 &= Q_bQ_cQ_0\x \cdot Q_0 \x = \delta_{bc} \\
Q_b e_{c+m} \cdot e_a &= Q_bQ_cQ_0\x \cdot Q_a \x = Q_aQ_bQ_cQ_0 \x \cdot \x = L^a_{bc} \\
Q_be_{c+m} \cdot e_{a+m} &= Q_bQ_cQ_0 \x \cdot Q_aQ_0 \x = -Q_aQ_bQ_c \x \cdot \x =0 \\
Q_b e_{c+m} \cdot e_\alpha &= Q_bQ_cQ_0 \x \cdot e_\alpha = -Q_bQ_c\x \cdot e_\alpha = F^\alpha_{c+m\,b} \\
Q_b e_{c+m} \cdot e_\mu &= Q_bQ_cQ_0 \x \cdot e_\mu = Q_bQ_c\x \cdot e_\mu = -F^\mu_{c+m\,b}
\end{align*}
give the expansion of $Q_b e_{c+m}$.  Using also~\eqref{eq:7:13}, we find
\begin{align*}
Q_b e_\alpha \cdot \x &= e_\alpha \cdot Q_b \x =0 \\
Q_b e_\alpha \cdot e_0 &= Q_b e_\alpha \cdot Q_0 \x = Q_b e_\alpha \cdot \x =0 \\
Q_b e_\alpha \cdot e_a &= Q_b e_\alpha \cdot Q_a \x = e_\alpha \cdot Q_b Q_a \x = F^\alpha_{b+m\,a} \\
Q_b e_\alpha \cdot e_\beta &= 0 \\
Q_b e_\alpha \cdot e_\mu &= Q_b e_\alpha \cdot e_\mu = -2F^\mu_{\alpha b}
\end{align*}
give the expansion of $Q_b e_\alpha$.
\begin{align*}
Q_b e_\mu \cdot \x &= e_\mu \cdot Q_b \x =0 \\
Q_b e_\mu \cdot e_0 &= Q_b e_\mu \cdot Q_0 \x = - Q_b e_\mu \cdot \x = 0 \\
Q_b e_\mu \cdot e_a &= Q_b e_\mu \cdot Q_a \x = e_\mu \cdot Q_b Q_a \x = F^\mu_{b+m\,a} \\
Q_b e_\mu \cdot e_{a+m} &= Q_b e_\mu \cdot Q_a Q_0 \x = e_\mu \cdot Q_b Q_a \x = -F^\mu_{a+m\, b} \\
Q_b e_\mu \cdot e_\alpha &= -2 F^\mu_{\alpha b} \\
Q_b e_\mu \cdot e_\nu &= 0
\end{align*}
give the expansion of $Q_b e_\mu$.
\end{proof}

\begin{lemma}\label{le:7:3}  For the Darboux frame field~\eqref{eq:7:14} along $\x$,
\begin{equation}\label{eq:7:34}
\aligned
\theta^b_a &= \tau^b_a + L^b_{ac}\omega^{c+m} + F^\alpha_{a+m\, b} \omega^\alpha + F^\mu_{a+m\, b}\omega^\mu \\
\theta^{b+m}_{a+m} &= \tau^b_a + L^c_{ab}\omega^c + F^\alpha_{a+m\,b} \omega^\alpha + F^\mu_{a+m\, b} \omega^\mu
\endaligned
\end{equation}
and therefore
\begin{equation}\label{eq:7:35}
\theta^b_a - \theta^{b+m}_{a+m} = L^b_{ac}(\omega^c+\omega^{c+m})
\end{equation}
\end{lemma}

\begin{proof}  Using~\eqref{eq:7:9} and~\eqref{eq:7:15}, we find
\begin{equation}\label{eq:7:36}
\aligned
\theta^b_a &= de_a \cdot e_b = d(Q_a\x)\cdot e_b \\
&= (\tau^i_a Q_i\x + \omega^p Q_a e_p + \omega^\alpha Q_a e_\alpha + \omega^\mu Q_a e_\mu) \cdot e_b \\
&= \tau^b_a + \omega^{c+m} Q_ae_{c+m} \cdot e_b + \omega^\alpha Q_ae_\alpha \cdot e_b + \omega^\mu Q_a e_\mu \cdot e_b
\endaligned
\end{equation}
which combined with~\eqref{eq:7:33} gives the first formula in~\eqref{eq:7:34}.  The second formula is derived in the same way.
\end{proof}

\section{Necessary conditions to be FKM}\label{section8}

Let $\tilde \x, e_a, e_p, e_\alpha, e_\mu, \tilde e_0$ be a second order frame field~\eqref{eq:2:0a} in $U \subset M$ along an isoparametric hypersurface $\tilde \x: M \to S^n$.  We continue using the index conventions in~\eqref{eq:4:convention}.
Let $\x = \cos s_1\, \tilde \x + \sin s_1\,\tilde e_0$ be a focal submanifold and let $e_0 = -\sin s_1\,\tilde \x + \cos s_1\, e_0$ so that $\x, e_0, e_a, e_p, e_\alpha, e_\mu$ is a Darboux frame field~\eqref{eq:4:12a} along $\x$ on $U$.  Let $\omega^a, \omega^p, \omega^\alpha, \omega^\mu$ be its coframe field~\eqref{eq:4:10} on $U$.  We look for conditions on this Darboux frame field which imply that $\x$ comes from an FKM construction.

\begin{proposition}\label{pr:8:1}  Suppose that $\x$ satisfies the spanning condition (definition~\ref{spanning}) on $U$.  If
\begin{equation}\label{eq:8:1}
F^\mu_{\alpha \, a+m} = F^\mu_{\alpha a}
\end{equation}
on $U$, then
\begin{align}
&F^\alpha_{b+m\, a} + F^\alpha_{a+m\, b} = 0 \label{eq:8:2} \\
&F^\mu_{b+m\, a} + F^\mu_{a+m\, b} = 0  \label{eq:8:3} \\
&\theta^a_b - \theta^{a+m}_{b+m} = L^a_{bc}(\omega^c + \omega^{c+m}), \text{ where $L^a_{bc} = - L^b_{ac} = - L^a_{cb}$} \label{eq:8:4}
\end{align}
on $U$.
\end{proposition}

\begin{remark}\label{re:8:1}  By Corollary~\ref{co:7:1} and
  Lemma~\ref{le:7:3}, equations~\eqref{eq:8:1} - \eqref{eq:8:4} hold
  for the Darboux frame field~\eqref{eq:7:14} defined along an FKM
  $\x$. 
\end{remark}

\begin{proof}  The summation convention is not used in this proof.    If we subtract the fourth equation in~\eqref{eq:5:2}, with $p = a+m$, from the third equation in~\eqref{eq:5:2}, we obtain
\begin{equation}\label{eq:8:5}
\sum_i (F^\mu_{\alpha a i} - F^\mu_{\alpha\, a+m\, i}) \omega^i = \sum_b F^\mu_{\alpha b} ( \theta^{b+m}_{a+m} - \theta^b_a)
\end{equation}
Putting \eqref{eq:8:1} into the second equation of~\eqref{eq:5:9} gives
\begin{equation}\label{eq:8:6}
F^\mu_{\alpha a \b} = \sum_b(F^\mu_{\alpha b} F^\b_{b+m\, a} + 2F^\mu_{\b b} F^\alpha_{b+m\, a})
\end{equation}
and putting~\eqref{eq:8:1} into the second equation of~\eqref{eq:5:10} gives
\begin{equation}\label{eq:8:7}
F^\mu_{\alpha\, a+m\, \b} = -\sum_b(F^\mu_{\alpha b} F^\b_{a+m\, b} +
2F^\mu_{\b b} F^\alpha_{a+m\, b}) 
\end{equation}
Subtracting \eqref{eq:8:7} from \eqref{eq:8:6} we get
\begin{equation}\label{eq:8:8}
F^\mu_{\alpha a \b} - F^\mu_{\alpha\,a+m\, \b} = \sum_b \left(
F^\mu_{\alpha b}(F^\b_{b+m\, a} + F^\b_{a+m\, b}) + 2F^\mu_{\b
  b}(F^\alpha_{b+m\,a} + F^\alpha_{a+m\, b} )\right) 
\end{equation}
Likewise, using the third equation in~\eqref{eq:5:9} and
in~\eqref{eq:5:10}, gives 
\begin{equation}\label{eq:8:9}
F^\mu_{\alpha a \nu} - F^\mu_{\alpha\, a+m\, \nu} = \sum_b \left(
F^\mu_{\alpha b}(F^\nu_{b+m\, a} + F^\nu_{a+m\, b}) + 2F^\nu_{\alpha
  b} ( F^\mu_{b+m\, a} + F^\mu_{a+m\,b})\right) 
\end{equation}
Expressing $\theta^a_b - \theta^{a+m}_{b+m}$ in terms of our coframe field, we have
\begin{equation}\label{eq:8:10}
\theta^a_b - \theta^{a+m}_{b+m} = \sum_c (L^a_{bc} \omega^c + L^a_{b\,c+m} \omega^{c+m}) + \sum_\alpha L^a_{b\a} \omega^\alpha + \sum_\mu L^a_{b\mu} \omega^\mu
\end{equation}
where the coefficients are smooth functions on $U$, each skew symmetric in $a,b$.  

By the Spanning Property, as expressed in (a') of Remark~\ref{re:6:8}, we
may assume the basis of $V_+$ chosen so that for some $\alpha$, the
set of vectors 
\[
\{ \sum_a F^\mu_{\a a} e_a:  \mbox{ all $\mu$}\}
\]
spans $V_0$.
Fix this choice of $\alpha$.  Substitute~\eqref{eq:8:10}
into~\eqref{eq:8:5} and compare the coefficients of $\omega^\alpha$ on
each side to obtain 
\begin{equation}\label{eq:8:11}
F^\mu_{\alpha a \alpha} - F^\mu_{\alpha\, a+m\, \alpha} = \sum_b
F^\mu_{\alpha b} L^a_{b\a} 
\end{equation}
Compare this to~\eqref{eq:8:8}, in which we set $\b = \alpha$, to obtain
\begin{equation}\label{eq:8:12}
\sum_b F^\mu_{\alpha b} \left( 3(F^\alpha_{b+m\, a} + F^\alpha_{a+m\,
  b}) - L^a_{b\a} \right) = 0
\end{equation}
for all $a$ and $\mu$.  By the Spanning Property, then, the vectors 
\[
\sum_b(3(F^\a_{b+m\,a} + F^\a_{a+m\,b}) - L^a_{b\a})e_b
\]
for each $a$ and $\mu$, are orthogonal to every vector in
$V_0$. Therefore, 
\begin{equation}\label{eq:8:13}
3(F^\alpha_{b+m\, a} + F^\alpha_{a+m\, b}) = L^a_{b\a}
\end{equation}
The left side of this equation is symmetric in $a,b$, while the
right side is skew-symmetric in $a,b$.  Therefore, for our choice
of $\alpha$,~\eqref{eq:8:2}
holds and
\begin{equation}\label{eq:8:14}
L^a_{b\a} = 0
\end{equation}
for all $a,b$.    Now, \eqref{eq:8:8} becomes, for our choice of $\a$ and
for any $\beta$,
\begin{equation}\label{eq:8:14a}
F^\mu_{\a a \b} - F^\mu_{\a\,a+m\,\b} = \sum_b F^\mu_{\a b}
(F^\b_{b+m\,a} + F^\b_{a+m\,b})
\end{equation}
Substitute \eqref{eq:8:10} into \eqref{eq:8:5} and compare the
coefficient of $\omega^\b$ with~\eqref{eq:8:14a} to obtain
\[
\sum_b F^\mu_{\a b}(F^\b_{b+m\,a} + F^\b_{a+m\,b} - L^a_{b\b}) = 0
\]
for all $a$, $\b$, and  $\mu$.  Again, the Spanning Property then implies that
\[
F^\b_{b+m\, a} + F^\b_{a+m\, b} = L^a_{b\b}
\]
for all $a$, $b$, and $\b$.  Hence, as before, each side of this
equation must be zero. 
Therefore,~\eqref{eq:8:2} and~\eqref{eq:8:14} hold for all $a$, $b$, and $\a$.

We can prove~\eqref{eq:8:3} and
\begin{equation}\label{eq:8:15}
L^a_{b\mu} = 0
\end{equation}
for all $a$, $b$ and $\mu$
in a similar way, by first fixing an appropriate $\mu$ and comparing
coefficients of 
$\omega^\mu$ in~\eqref{eq:8:5} after substitution of~\eqref{eq:8:10}
into it.  In this case (b') of the spanning property is used. 

With \eqref{eq:8:2} and~\eqref{eq:8:3} now true, we see
that~\eqref{eq:8:8} and~\eqref{eq:8:9} become 
\begin{equation}\label{eq:8:16}
F^\mu_{\alpha a \b} = F^\mu_{\alpha\,a+m\,\b}, \quad F^\mu_{\alpha a
  \nu} = F^\mu_{\alpha\,a+m\,\nu} 
\end{equation}
and~\eqref{eq:8:14} and~\eqref{eq:8:15} substituted into~\eqref{eq:8:10} give
\begin{equation}\label{eq:8:17}
\theta^a_b - \theta^{a+m}_{b+m} = \sum_c(L^a_{bc}\omega^c + L^a_{b\,c+m} \omega^{c+m})
\end{equation}
Substitute this into~\eqref{eq:8:5} and compare coefficients of
$\omega^c$ and $\omega^{c+m}$ to get
\begin{equation}\label{eq:8:20}
\aligned
\sum_b F^\mu_{\alpha b} L^a_{bc} &= F^\mu_{\alpha ac} - F^\mu_{\alpha\,a+m\, c} \\
\sum_b F^\mu_{\alpha b} L^a_{b\,c+m} &= F^\mu_{\alpha a\,c+m} - F^\mu_{\alpha\,a+m\,c+m}
\endaligned
\end{equation}
Subtracting gives
\begin{equation}\label{eq:8:21}
\sum_c F^\mu_{\alpha c}(L^a_{cb} - L^a_{c\,b+m}) = F^\mu_{\alpha a b} - F^\mu_{\alpha\,a+m \, b} - F^\mu_{\alpha a\,b+m} + F^\mu_{\alpha\,a+m\,b+m}
\end{equation}
We want to show now that the right hand side of this equation is zero on $U$.  To that end, we begin with the first equation in~\eqref{eq:5:9}, which says
\begin{equation}\label{eq:8:22}
F^\mu_{\alpha a b} = -\frac12 \sum_c F^\mu_{c+m\, a} F^\alpha_{c+m\, b} + \frac12 \sum_c F^\mu_{c+m\, b} F^\alpha_{c+m\, a}
\end{equation}
and~\eqref{eq:5:11} says
\begin{equation}\label{eq:8:23}
F^\mu_{\alpha\,a+m\,b} = \frac12 F^\mu_{a+m\, b\a}, \quad F^\mu_{\alpha a\,b+m} = \frac12 F^\mu_{b+m\,a\a}
\end{equation}
and the first equation in~\eqref{eq:5:10} says
\begin{equation}\label{eq:8:24}
F^\mu_{\alpha\,a+m\,b+m} = \frac 12 \sum_c F^\mu_{a+m\,c} F^\alpha_{b+m\,c} - \frac12 \sum_c F^\mu_{b+m\,c} F^\alpha_{a+m\,c}
\end{equation}
Hence, using~\eqref{eq:8:2} and~\eqref{eq:8:3}, the right hand side of~\eqref{eq:8:21} is
\begin{equation*}\label{eq:8:25}
\aligned
&F^\mu_{\alpha a b} - F^\mu_{\alpha\,a+m \, b} - F^\mu_{\alpha a\,b+m} + F^\mu_{\alpha\,a+m\,b+m} = \\
&-\frac12\sum_c F^\mu_{c+m\, a} F^\alpha_{c+m\,b} + \frac12 \sum_c F^\mu_{c+m\,b}F^\alpha_{c+m\,a} - \frac12 F^\mu_{a+m\,b\a} -\frac12 F^\mu_{b+m\,a\a} \\
\quad &+\frac12 \sum_c F^\mu_{a+m\,c} F^\alpha_{b+m\,c} - \frac12 \sum_c F^\mu_{b+m\,c} F^\alpha_{a+m\,c} =
-\frac12(F^\mu_{a+m\,b\a} + F^\mu_{b+m\,a\a})\\
&+\frac12 \sum_c F^\mu_{a+m\,c}(F^\alpha_{c+m\,b} + F^\alpha_{b+m\,c}) + \frac12 \sum_c F^\mu_{c+m\,b}(F^\alpha_{c+m\,a} + F^\alpha_{a+m\,c})  \\
&= -\frac12(F^\mu_{a+m\,b\a} + F^\mu_{b+m\,a\a})
\endaligned
\end{equation*}
and so we want to show that this last term is zero on $U$ when~\eqref{eq:8:1},~\eqref{eq:8:2} and~\eqref{eq:8:3} hold.  By the second equation in~\eqref{eq:5:2},
\begin{equation*}
\sum_i F^\mu_{a+m\,bi} \omega^i = dF^\mu_{a+m\,b} - \sum_c F^\mu_{c+m\, b} \theta^{c+m}_{a+m} - \sum_c F^\mu_{a+m\, c} \theta^c_b + \sum_\nu F^\nu_{a+m\,b}\theta^\mu_\nu
\end{equation*}
and
\begin{equation*}
\sum_i F^\mu_{b+m\,ai}\omega^i = dF^\mu_{b+m\,a} - \sum_a F^\mu_{c+m\,a} \theta^{c+m}_{b+m} - \sum_c F^\mu_{b+m\,c} \theta^c_a + \sum_\nu F^\nu_{b+m\, a}\theta^\mu_\nu
\end{equation*}
Sum these two equations and use~\eqref{eq:8:2} and~\eqref{eq:8:3} to get
\begin{equation*}\label{eq:8:26}
\sum_i(F^\mu_{a+m\,bi} + F^\mu_{b+m\,ai})\omega^i = \sum_c\left( F^\mu_{b+m\,c}(\theta^{c+m}_{a+m} - \theta^c_a) + F^\mu_{a+m\,c}(\theta^{c+m}_{b+m} - \theta^c_b)\right)
\end{equation*}
By~\eqref{eq:8:17}, the right hand side of this equation is in the span of the set of $1$-forms $\{\omega^c, \omega^{c+m}: c = 1,\dots,m\}$, and therefore the coefficients of $\omega^\alpha$ and $\omega^\mu$ on the left hand side must vanish, to give
\begin{equation}\label{eq:8:27}
F^\mu_{a+m\,b\a} + F^\mu_{b+m\,a\a} = 0, \quad F^\mu_{a+m\,b\nu} + F^\mu_{b+m\,a\nu} =0
\end{equation}
and we have finally proved that the right hand side of~\eqref{eq:8:21}
is zero on $U$, and therefore 
\begin{equation}\label{eq:8:28}
\sum_b F^\mu_{\alpha b}(L^a_{bc} - L^a_{b\,c+m} )=0
\end{equation}
on $U$, for all $a$, $c$, $\alpha$, and $\mu$.  Multiplying this
equation by the 
$X = \sum x_\alpha e_\alpha$ of (a) of the spanning property, we conclude
that 
\begin{equation}\label{eq:8:29}
L^a_{bc} - L^a_{b\,c+m} = 0
\end{equation}
on $U$ for all $a,b,c$.  Substitution of this into~\eqref{eq:8:17} gives
\begin{equation}\label{eq:8:30}
\theta^a_b - \theta^{a+m}_{b+m} = \sum_c L^a_{bc}(\omega^c+ \omega^{c+m})
\end{equation}
To complete the proof of~\eqref{eq:8:4}, it remains to show that
\begin{equation}\label{eq:8:31}
L^a_{bc} + L^a_{cb}=0
\end{equation}
on $U$, for all $a,b,c$. 
By~\eqref{eq:5:2}, \eqref{eq:8:1} and~\eqref{eq:8:30}, and the known
skew-symmetry $L^a_{bc} = - 
L^b_{ac}$, we have
\begin{equation}\label{eq:8:31a}
\sum_i F^\mu_{\alpha\,a+m\,i} \omega^i = \sum_i F^\mu_{\alpha a i} \omega^i +
\sum_{b,c} F^\mu_{\alpha b} L^b_{ac}(\omega^c + \omega^{c+m})
\end{equation}
Comparing the coefficients of $\omega^c$, we have
\begin{equation}\label{eq:8:31b}
F^\mu_{\alpha\, a+m\, c} = F^\mu_{\alpha a c} + \sum_b F^\mu_{\alpha
  b} L^b_{ac}
\end{equation}
Interchanging $a$ and $c$ and then summing, we have
\begin{equation}\label{eq:8:31c}
F^\mu_{\alpha\, a+m\, c} + F^\mu_{\alpha\, c+m\, a} = F^\mu_{\alpha a
  c}+ F^\mu_{\alpha c a} + \sum_b F^\mu_{\alpha b}( L^b_{ac} +
  L^b_{ca})
\end{equation}
By the first equation in~\eqref{eq:5:9}, we have 
\begin{equation}\label{eq:8:31d}
F^\mu_{\alpha ac} + F^\mu_{\alpha ca} = 0
\end{equation}
Hence
\begin{equation}\label{eq:8:31e}
F^\mu_{\alpha\,a+m\,c} + F^\mu_{\alpha\,c+m\,a} = \sum_b F^\mu_{\alpha b}(L^b_{ac} + L^b_{ca})
\end{equation}
on $U$ for all $\alpha$ and $\mu$.  In~\eqref{eq:8:31a} compare the
coefficients of $\omega^{c+m}$ to get
\[
F^\mu_{\alpha\, a+m\, c+m} = F^\mu_{\alpha a\,c+m} + \sum_b
F^\mu_{\alpha b} L^b_{ac}
\]
Interchange $a$ and $c$ and sum, to get
\begin{equation}\label{eq:8:31f}
F^\mu_{\alpha\, a+m\,c+m} + F^\mu_{\alpha\, c+m\, a+m} = F^\mu_{\alpha
  a\,c+m} + F^\mu_{\alpha c\, a+m} + \sum_b F^\mu_{\alpha b}(L^b_{ac}
  + L^b_{ca})
\end{equation}
By the first equation in~\eqref{eq:5:10},
\[
F^\mu_{\alpha\,a+m\,c+m} + F^\mu_{\alpha \,c+m\, a+m} = 0
\]
and the last equation in~\eqref{eq:5:11} says that
\[
F^\mu_{\alpha a \,c+m} = F^\mu_{\alpha\,c+m\,a} \mbox{ and }
F^\mu_{\alpha c \,a+m} = F^\mu_{\alpha\, a+m\,c}
\]
Therefore, \eqref{eq:8:31f} is
\begin{equation}\label{eq:8:31g}
F^\mu_{\alpha\,c+m\,a} + F^\mu_{\alpha\, a+m\, c} = -\sum_b
F^\mu_{\alpha b} ( L^b_{ac} + L^b_{ca})
\end{equation}
Combining this with \eqref{eq:8:31e}, we conclude that
\begin{equation}\label{eq:8:32}
\sum_b F^\mu_{\alpha b}(L^b_{ac}+L^b_{ca}) = 0
\end{equation}
for all $a,c, \alpha, \mu$.  The spanning property then implies~\eqref{eq:8:31}.
\end{proof}

Resume use of the summation convention.

\begin{proposition}\label{pr:8:3} If equations~\eqref{eq:8:1}
  through~\eqref{eq:8:4} hold on $U$, then 
\begin{align}
F^\alpha_{c+m\, a} L^c_{bd} + F^\alpha_{c+m\, b}L^c_{ad} &=
2(F^\mu_{\alpha a}F^\mu_{d+m\, b} + F^\mu_{\alpha b}F^\mu_{d+m\, a})
\label{eq:8:33} \\ 
F^\mu_{c+m\, a} L^c_{bd} + F^\mu_{c+m\, b}L^c_{ad} &=
2(F^\alpha_{b+m\, d}F^\mu_{\alpha a} + F^\alpha_{a+m\, d}
F^\mu_{\alpha b}) \label{eq:8:34} \\ 
F^\mu_{\alpha\, b+m\,a} = L^c_{ba}F^\mu_{\alpha c} &- \frac12
F^\mu_{d+m\, b}F^\alpha_{d+m\,a} + \frac12 F^\mu_{d+m\,a}
F^\alpha_{d+m\,b} \label{eq:8:35} 
\end{align}
\end{proposition}

\begin{proof}  These identities come from
  differentiating~\eqref{eq:8:1} through~\eqref{eq:8:3}.  Using our
  definition of covariant derivative in~\eqref{eq:5:2}, we have 
\begin{equation}\label{eq:8:36}
\aligned
dF^\alpha_{b+m\, a} + F^\beta_{b+m\, a} \theta^\alpha_\beta -
F^\alpha_{c+m\,a} \theta^{c+m}_{b+m} - F^\alpha_{b+m\, c} \theta^c_a
&= F^\alpha_{b+m\, ai} \omega^i  \\ 
dF^\alpha_{a+m\,b} + F^\beta_{a+m\, b} \theta^\alpha_\beta -
F^\alpha_{c+m\,b} \theta^{c+m}_{a+m} - F^\alpha_{a+m\, c} \theta^c_b
&= F^\alpha_{a+m\, bi} \omega^i 
\endaligned
\end{equation}
Summing these two equations and using~\eqref{eq:8:2}
and~\eqref{eq:8:4}, we get 
\begin{equation}\label{eq:8:37}
(F^\alpha_{c+m\, a} L^c_{bd} + F^\alpha_{c+m\, b} L^c_{ad})(\omega^d +
\omega^{d+m}) = (F^\alpha_{b+m\, ai} + F^\alpha_{a+m\,bi})\omega^i 
\end{equation}
Equating the coefficients of $\omega^d$, we have
\begin{equation}\label{eq:8:38}
F^\alpha_{c+m\,a}L^c_{bd} + F^\alpha_{c+m\,b}L^c_{ad} =
F^\alpha_{b+m\,ad} + F^\alpha_{a+m\, bd} 
\end{equation}
From~\eqref{eq:5:7} we see that the right side of~\eqref{eq:8:38} is
\begin{equation}\label{eq:8:39}
\aligned
-F^\mu_{b+m\,a} F^\mu_{\alpha d} - &2F^\mu_{b+m\,d}F^\mu_{\alpha a} - F^\mu_{a+m\,b} F^\mu_{\alpha d} - 2F^\mu_{a+m\,d}F^\mu_{\alpha b} \\
&= 2F^\mu_{d+m\,b}F^\mu_{\alpha a} + 2F^\mu_{d+m\,a}F^\mu_{\alpha b}
\endaligned
\end{equation}
where the last equality comes from using~\eqref{eq:8:3}.  Now~\eqref{eq:8:33} follows from~\eqref{eq:8:38} and~\eqref{eq:8:39}.  Equating the coefficients of $\omega^{d+m}$ in~\eqref{eq:8:37} leads again to~\eqref{eq:8:33}. Equating the other coefficients leads to the identities
\begin{equation}\label{eq:8:40}
F^\alpha_{b+m\,a\beta} + F^\alpha_{a+m\,b\beta} = 0 \mbox{ and }
F^\alpha_{b+m\,a\mu} + F^\alpha_{a+m\, b\mu} = 0
\end{equation}
We next find the consequences of taking the covariant derivative of
equation~\eqref{eq:8:3}.  Again by~\eqref{eq:5:2}, we have 
\begin{equation}\label{eq:8:41}
\aligned
dF^\mu_{b+m\,a} + F^\nu_{b+m\,a}\theta^\mu_\nu - F^\mu_{c+m\,a} \theta^{c+m}_{b+m} - F^\mu_{b+m\, c}\theta^c_a &= F^\mu_{b+m\,ai} \omega^i  \\
dF^\mu_{a+m\,b} + F^\nu_{a+m\,b}\theta^\mu_\nu - F^\mu_{c+m\,b} \theta^{c+m}_{a+m} - F^\mu_{a+m\, c}\theta^c_b &= F^\mu_{a+m\,bi} \omega^i
\endaligned
\end{equation}
Summing these equations and using~\eqref{eq:8:3} and~\eqref{eq:8:4}, we get
\begin{equation}\label{eq:8:42}
(F^\mu_{c+m\,a} L^c_{bd} + F^\mu_{c+m\,b} L^c_{ad})( \omega^d + \omega^{d+m}) = (F^\mu_{b+m\,ai} + F^\mu_{a+m\,bi})\omega^i
\end{equation}
Equating the coefficients of $\omega^d$ we have
\begin{equation}\label{eq:8:43}
F^\mu_{c+m\,a}L^c_{bd} + F^\mu_{c+m\,b}L^c_{ad} = F^\mu_{b+m\,ad} + F^\mu_{a+m\,bd}
\end{equation}
By~\eqref{eq:5:8}, the right side of~\eqref{eq:8:43} is
\begin{equation}\label{eq:8:44}
F^\alpha_{b+m\,a}F^\mu_{\alpha d} + 2F^\alpha_{b+m\,d}F^\mu_{\alpha a} + F^\alpha_{a+m\,b}F^\mu_{\alpha d} + 2F^\alpha_{a+m\,d}F^\mu_{\alpha b}
\end{equation}
Using~\eqref{eq:8:2} in~\eqref{eq:8:44}, we then arrive
at~\eqref{eq:8:34}.  Equating coefficients of $\omega^{d+m}$
in~\eqref{eq:8:42} also leads to~\eqref{eq:8:34}.  Equating
coefficients of $\omega^\alpha$ and of $\omega^\mu$ gives 
\begin{equation}\label{eq:8:45}
F^\mu_{b+m\,a\alpha} + F^\mu_{a+m\, b\alpha} = 0 \mbox{ and }
F^\mu_{b+m\,a\nu} + F^\mu_{a+m\, b\nu} = 0
\end{equation}
Finally, substitute the first equation
of~\eqref{eq:5:9} into~\eqref{eq:8:31b} to arrive at~\eqref{eq:8:35}.
\end{proof}

We define the \textit{covariant derivatives} of the $L^a_{bc}$ to be
the coefficients $L^a_{bci}$ of the $1$-form 
\begin{equation}\label{eq:8:52}
dL^a_{bc} + L^d_{bc}\theta^a_d - L^a_{dc}\theta^d_b - L^a_{bd} \theta^d_c = L^a_{bci}\omega^i
\end{equation}

\begin{remark}\label{re:8:21}
If the $L^a_{bc}$ are skew symmetric in all three indices, then the
functions $L^a_{bci}$ are skew symmetric in $a,b,c$. 
\end{remark}

\begin{proposition}\label{pr:8:5}  If equations~\eqref{eq:8:1} through~\eqref{eq:8:4} hold, then the $L^a_{bcd}$ are skew symmetric in all four indices, and
\begin{equation}\label{eq:8:53}
\aligned
L^a_{bcd} &= \frac12(\delta_{ad}\delta_{bc} - \delta_{ac}\delta_{bd}) + \frac12 (L^a_{ce}L^e_{bd} - L^a_{de}L^e_{bc})  \\
&\quad + F^\alpha_{c+m\, a} F^\alpha_{d+m\, b} - F^\mu_{c+m\, b} F^\mu_{d+m\, a} 
\endaligned
\end{equation}
\begin{equation}\label{eq:8:54}
\aligned
L^a_{bc\,d+m} &= \frac12(\delta_{ac}\delta_{bd} - \delta_{ad}\delta_{bc}) + L^a_{be}L^e_{dc} + \frac12(L^a_{ce}L^e_{bd} - L^a_{de}L^e_{bc})  \\
&\quad + F^\alpha_{d+m\, a}F^\alpha_{c+m\, b} - F^\mu_{d+m\, b} F^\mu_{c+m\, a} 
\endaligned
\end{equation}
\begin{align}
L^a_{bc\alpha} &= L^a_{be} F^\alpha_{e+m\, c} + 2(F^\mu_{\alpha a}F^\mu_{c+m\, b} - F^\mu_{\alpha b} F^\mu_{c+m\, a}) \label{eq:8:55} \\
L^a_{bc\mu} &= L^a_{be} F^\mu_{e+m\, c} + 2(F^\mu_{\alpha b} F^\alpha_{c+m\,a} - F^\mu_{\alpha a} F^\alpha_{c+m\,b}) \label{eq:8:56} 
\end{align}
\begin{equation}\label{eq:8:57}
\aligned
&2\delta_{ac}\delta_{bd} - \delta_{ad}\delta_{bc} - \delta_{ab}\delta_{dc} = \\
L^a_{be}L^c_{de} + L^a_{de}L^c_{be} &+ 2(F^\alpha_{b+m\, c} F^\alpha_{d+m\, a} + F^\alpha_{b+m\, a} F^\alpha_{d+m\, c})
\endaligned
\end{equation}
\begin{equation}\label{eq:8:58}
L^a_{bc\,d+m} + L^a_{bd\,c+m} = 0 
\end{equation}
\end{proposition}

\begin{proof} This proposition is a consequence of taking the exterior derivative of~\eqref{eq:8:4}. Notice that~\eqref{eq:8:58} follows directly from~\eqref{eq:8:54}. 

Using~\eqref{eq:4:13} and the structure equations~\eqref{eq:2:1a}, we find
\begin{equation*}
\aligned
&d(\theta^a_b - \theta^{a+m}_{b+m}) = \omega^a \w\omega^b - \omega^{a+m}\w\omega^{b+m} \\
& + (F^\alpha_{c+m\,a}F^\alpha_{d+m\,b} + F^\mu_{c+m\,a}F^\mu_{d+m\,b})(\omega^{c+m}\w \omega^{d+m} - \omega^c\w\omega^d) \\
&+ [L^a_{dc}\theta^d_b - L^d_{bc}\theta^a_d + L^a_{ed}L^e_{bc}(\omega^d+\omega^{d+m}) 
 + 2(F^\mu_{\alpha a}F^\mu_{c+m\,b} - F^\mu_{\alpha b} F^\mu_{c+m\,a}) \omega^\alpha \\
&+ 2(F^\mu_{\alpha b}F^\alpha_{c+m\,a} - F^\mu_{\alpha a} F^\alpha_{c+m\,b}) \omega^\mu ]\w (\omega^c + \omega^{c+m}) 
\endaligned
\end{equation*}
Using~\eqref{eq:5:1}, we find
\begin{equation*}
\aligned
d(L^a_{bc}(\omega^c + \omega^{c+m})) &= (dL^a_{bc} - L^a_{bd}\theta^d_c - L^a_{be}L^e_{dc}\omega^{d+m} \\
&\quad - L^a_{bd}F^\alpha_{d+m\,c}\omega^\alpha - L^a_{bd}F^\mu_{d+m\,c}\omega^\mu) \w (\omega^c + \omega^{c+m})
\endaligned
\end{equation*}
The exterior differential of~\eqref{eq:8:4} is obtained by equating the preceding two equations and using~\eqref{eq:8:52}, to get
\begin{equation}\label{eq:8:59}
\aligned
&L^a_{bei}\omega^i \w (\omega^e + \omega^{e+m}) = \\
&\quad [-\delta_{ae}\omega^b - \delta_{be}\omega^{a+m} \\
&+ (L^a_{dc}L^d_{be}+ F^\alpha_{e+m\,a}F^\alpha_{c+m\,b} + F^\mu_{e+m\,a}F^\mu_{c+m\,b})\omega^c \\
&+ (L^a_{bd}L^d_{ce} + L^a_{dc}L^d_{be} + F^\alpha_{c+m\,a}F^\alpha_{e+m\,b} + F^\mu_{c+m\,a}F^\mu_{e+m\,b})\omega^{c+m} \\
&\quad + (L^a_{bc}F^\alpha_{c+m\,e}+ 2F^\mu_{\alpha a}F^\mu_{e+m\,b} - 2F^\mu_{e+m\,a}F^\mu_{\alpha b})\omega^\alpha \\
&\quad + (L^a_{bc}F^\mu_{c+m\,e} +
2F^\alpha_{e+m\,a}F^\mu_{\alpha b} - 2F^\mu_{\alpha a}F^\mu_{e+m\,b})\omega^\mu]\w (\omega^e + \omega^{e+m})
\endaligned
\end{equation}
Equating the skew symmetrized coefficients of $\omega^c \w \omega^e$ in this equation, we have
\begin{equation}\label{eq:8:60}
\aligned
L^a_{bec}-L^a_{bce} &= L^a_{dc}L^d_{be} - L^a_{de}L^d_{bc} - \delta_{ae}\delta_{bc} + \delta_{ac}\delta_{be} \\
&\quad + F^\alpha_{e+m\,a}F^\alpha_{c+m\,b} - F^\alpha_{c+m\,a}F^\alpha_{e+m\,b} \\
&\quad + F^\mu_{e+m\,a}F^\mu_{c+m\,b} - F^\mu_{c+m\,a}F^\mu_{e+m\,b} 
\endaligned
\end{equation}
Rewrite~\eqref{eq:8:60} with $b$ and $e$ interchanged and add the result to~\eqref{eq:8:60}.  Using that $L^d_{be}$, $L^a_{bec}$, $F^\alpha_{e+m\,b}$ and $F^\mu_{e+m\,b}$ are all skew-symmetric in $b$ and $e$, we get from this sum
\begin{equation}\label{eq:8:61}
\aligned
L^a_{bce}+L^a_{ecb} &= L^a_{bd}L^c_{ed} + L^a_{ed}L^c_{bd} + \delta_{ae}\delta_{bc} + \delta_{ab}\delta_{ec} - 2\delta_{ac}\delta_{be} \\
&\quad + F^\alpha_{b+m\,c}F^\alpha_{e+m\,a} + F^\alpha_{b+m\,a}F^\alpha_{e+m\,c} \\
&\quad + F^\mu_{e+m\,a}F^\mu_{b+m\,c} + F^\mu_{b+m\,a}F^\mu_{e+m\,c}
\endaligned
\end{equation}
Equating the coefficients of $\omega^c \w \omega^{e+m}$ in~\eqref{eq:8:59}, we find
\begin{equation}\label{eq:8:62}
L^a_{bec} - L^a_{bc\,e+m} = L^a_{dc}L^d_{be} - L^a_{bd}L^d_{ec} - L^a_{de}L^d_{bc}
\end{equation}
Rewrite this equation with $b$ and $c$ interchanged and add the result to~\eqref{eq:8:62}.  From the skew symmetry of $L^a_{bc}$ and $L^a_{bcd}$ in $a,b,c$, it follows from this sum that
\begin{equation}\label{eq:8:63}
L^a_{bec} + L^a_{ceb} = 0
\end{equation}
from which we conclude that $L^a_{bcd}$ is skew symmetric in all four indices.  Putting~\eqref{eq:8:63} into~\eqref{eq:8:61}, interchanging $d$ and $e$ and using the first equation in~\eqref{eq:5:6}, we arrive at~\eqref{eq:8:57}.  Putting~\eqref{eq:8:63} into~\eqref{eq:8:60} and using the first equation of~\eqref{eq:5:6}, we get~\eqref{eq:8:53}.  Substitute~\eqref{eq:8:53} into~\eqref{eq:8:62} to obtain~\eqref{eq:8:54}.  Go back to~\eqref{eq:8:59} and equate coefficients of $\omega^\alpha \w \omega^c$ to obtain~\eqref{eq:8:55}, and equate coefficients of $\omega^\mu \w \omega^c$ to obtain~\eqref{eq:8:56}.
\end{proof}

\section{A sufficient condition to be FKM}\label{section9}  Let $\tilde \x, e_a, e_p, e_\alpha, e_\mu, \tilde e_0$ be a second order frame field~\eqref{eq:2:0a} in $U \subset M$ along an isoparametric hypersurface $\tilde \x: M \to S^n \subset \R^{n+1}$.  We continue using the index conventions in~\eqref{eq:4:convention}.
Let $\x = \cos s_1\, \tilde \x + \sin s_1\,\tilde e_0$ be a focal
submanifold and let $e_0 = -\sin s_1\,\tilde \x + \cos s_1\, \tilde e_0$ so
that  
\begin{equation}\label{eq:g:1}
\x, e_0, e_a, e_p, e_\alpha, e_\mu
\end{equation}
is a Darboux frame field~\eqref{eq:4:12a} along $\x$ on $U$.  Let 
\begin{equation}\label{eq:g:2}
\omega^a, \omega^p, \omega^\alpha, \omega^\mu
\end{equation}
be its coframe field~\eqref{eq:4:10} on $U$.

\begin{theorem}\label{th:g:1} If $\x$ satisfies the spanning condition (definition~\ref{spanning}) and condition~\eqref{eq:8:1}, $F^\mu_{\a\,a+m} = F^\mu_{\alpha a}$, on $U$, then it comes from an FKM construction.
\end{theorem}

\begin{proof} It is sufficient to prove the theorem locally, on some open neighborhood, because isoparametric hypersurfaces are algebraic.  For each point in $U$, the vectors of our Darboux frame field~\eqref{eq:g:1} form an orthonormal basis of $\R^{n+1}$.  Linear operators $Q_0,Q_a$ on $\R^{n+1}$, depending on the point in $U$, can thus be defined by~\eqref{eq:7:32} and~\eqref{eq:7:33}, which we recopy here for easier reference
\begin{equation}\label{eq:g:3}
\aligned
Q_0 \x &= e_0 &\quad Q_0 e_0 &= \x &\quad Q_0 e_a &= -e_{a+m} \\
Q_0 e_{a+m} &= -e_a &\quad Q_0 e_\alpha &= -e_\alpha &\quad Q_0 e_\mu &= e_\mu
\endaligned
\end{equation}
and for each $a$
\begin{equation}\label{eq:g:4}
\aligned
Q_a\x &= e_{a} \\
Q_a e_0 &= e_{a+m} \\
Q_a e_b &= \delta_{ab}\x - L^c_{ab}e_{c+m} + F^\alpha_{a+m\,b}e_\alpha + F^\mu_{a+m\,b}e_\mu \\
Q_a e_{b+m} &= \delta_{ab}e_0 + L^c_{ab}e_c + F^\alpha_{b+m\,a}e_\alpha - F^\mu_{b+m\,a}e_\mu \\
Q_a e_\alpha &= F^\alpha_{a+m\, b}e_b + F^\alpha_{b+m\,a}e_{b+m} -2 F^\mu_{\alpha\, a}e_\mu \\
Q_a e_\mu &= F^\mu_{a+m\, b}e_b - F^\mu_{b+m\, a} e_{b+m} -2 F^\mu_{\alpha\, a}e_\alpha 
\endaligned
\end{equation}
where the coefficients are defined now in~\eqref{eq:4:13} and~\eqref{eq:8:4}.  We first outline the quite elementary proof of the theorem, and then follow that with a proof of the details.  The first detail is:

\textbf{(I).}  At each point of $U$ these operators are symmetric, orthogonal and satisfy
\begin{equation}\label{eq:g:5}
Q_iQ_j+Q_jQ_i = 2\delta_{ij}I, \quad \mbox{ for $i,j = 0,1,\dots,m$}
\end{equation}
Given that, one next proves the second detail:

\textbf{(II).}  There exist a (constant) Clifford system $P_0,\dots,P_m$ on $\R^{n+1}$ and a smooth map 
\begin{equation}\label{eq:g:5a}
B:U \to SO(m+1)
\end{equation}
such that at every point of $U$,
\begin{equation}\label{eq:g:6}
Q_j = \sum_{i=0}^m B^i_j P_i, \quad \mbox{ for $j = 0,1,\dots,m$}
\end{equation}
It will then follow that $\x$ maps $U$ onto an open subset of the focal submanifold $M_+$ defined in~\eqref{eq:7:4} by this Clifford system, and that the Darboux frame field~\eqref{eq:7:14} coming from the FKM construction applied to $P_0,\dots,P_m$ coincides with our frame field~\eqref{eq:g:1}. Therefore, our $\x:U \to S^n$ coincides with the FKM construction applied to this Clifford system.

We turn now to the proof of detail (I).  
The verification that each $Q_i$ is symmetric can be done almost by inspection.  It is equally clear that $Q_0$ is orthogonal, since it sends the orthonormal basis~\eqref{eq:g:1} to an orthonormal basis.
The operator $Q_{a}$ sends the orthonormal basis~\eqref{eq:g:1} to the set of vectors given on the right hand side of~\eqref{eq:g:4}.  Among these vectors, $Q_a \x, Q_a e_0$ is an orthonormal pair orthogonal to the remaining vectors because $L^a_{bc}$ are skew symmetric in $a,b,c$ and $F^\alpha_{a+m\,b}$ and $F^\mu_{a+m\,b}$ are skew symmetric in $a$ and $b$.  

In the following verification that
\begin{equation*}
\{Q_a e_b, Q_a e_{b+m}, Q_a e_\alpha, Q_a e_\mu : b,a,\mu \}
\end{equation*}
is orthonormal, 
we do not use the Einstein summation convention as $a$ will always be a repeated index which is not summed.  We proceed through all the cases.
\begin{equation*} 
\aligned
Q_ae_b \cdot Q_ae_d &= \delta_{ab}\delta_{ad} + \sum_c L^b_{ac}L^d_{ac} \\
&+ \sum_\alpha F^\alpha_{a+m\,b}F^\alpha_{a+m\,d} + \sum_\mu F^\mu_{a+m\,b}F^\mu_{a+m\,d} = \delta_{bd}
\endaligned
\end{equation*}
by~\eqref{eq:8:57} with $c$ changed to $b$ and $b$ changed to $a$.
\begin{equation*} 
Q_a e_b\cdot Q_a e_{d+m} = \sum_\alpha F^\alpha_{a+m\,b}F^\alpha_{d+m\,a} - \sum_\mu F^\mu_{a+m\,b} F^\mu_{d+m\,a} = 0
\end{equation*}
by the first equation in~\eqref{eq:5:6}.
\begin{equation*} 
Q_a e_b\cdot Q_a e_\alpha = \sum_c L^b_{ac} F^\alpha_{c+m\,a} -2\sum_\mu F^\mu_{a+m\,b} F^\mu_{\alpha a} =0 
\end{equation*}
by~\eqref{eq:8:33} with $d$ changed to $a$.
\begin{equation*}  
Q_a e_b\cdot Q_a e_\mu = -\sum_c L^b_{ac} F^\mu_{c+m\,a} - 2\sum_\alpha F^\alpha_{a+m\,b} F^\mu_{\alpha a} = 0
\end{equation*}
by~\eqref{eq:8:34} with $d$ changed to $a$.
\begin{equation*} 
\aligned
Q_a e_{b+m} &\cdot Q_a e_{d+m} = \delta_{ab}\delta_{ad} + \sum_c L^c_{ab}L^c_{ad} \\ 
&+ \sum_\alpha F^\alpha_{b+m\,a} F^\alpha_{d+m\,a} 
+ \sum_\mu F^\mu_{b+m\,a} F^\mu_{d+m\,a} = \delta_{bd} 
\endaligned
\end{equation*}
by~\eqref{eq:8:57} with $c$ changed to $b$ and $b$ changed to $a$.
\begin{equation*}  
Q_a e_{b+m} \cdot Q_a e_\alpha = \sum_c L^c_{ab} F^\alpha_{a+m\,c} + 2\sum_\mu F^\mu_{b+m\,a} F^\mu_{\alpha a} = 0
\end{equation*}
by~\eqref{eq:8:33} with $d$ changed to $a$.
\begin{equation*}  
Q_a e_{b+m} \cdot Q_a e_\mu = \sum_c L^c_{ab} F^\mu_{a+m\,c} -2\sum_\alpha F^\alpha_{b+m\,a} F^\mu_{\alpha a} = 0
\end{equation*}
by~\eqref{eq:8:34} with $d$ changed to $a$.
\begin{equation*}  
Q_a e_\alpha \cdot Q_a e_\b = 2\sum_b F^\alpha_{b+m\,a} F^\b_{b+m\,a} + 4\sum_\mu F^\mu_{\alpha a} F^\mu_{\b a} = \delta_{\alpha\b} 
\end{equation*}
by the second equation in~\eqref{eq:5:6}.
\begin{equation*}  
Q_a e_\alpha \cdot Q_a e_\mu = \sum_b F^\alpha_{a+m\,b} F^\mu_{a+m\,b} - \sum_b F^\alpha_{b+m\,a} F^\mu_{b+m\,a} = 0
\end{equation*}
by~\eqref{eq:8:2} and~\eqref{eq:8:3}.
\begin{equation*}  
Q_a e_\mu \cdot Q_a e_\nu = 2\sum_b F^\mu_{b+m\,a} F^\nu_{b+m\,a} + 4\sum_\alpha F^\mu_{\alpha a}F^\nu_{\alpha a} = \delta_{\mu\nu}
\end{equation*}
by the fourth equation in~\eqref{eq:5:6}.  That completes the verification that each $Q_i$ is an orthogonal transformation.  

We proceed now to verify~\eqref{eq:g:5}.  For this we return to using the Einstein summation convention.
Clearly $Q_0^2 = I$.  To verify that $Q_0Q_a + Q_aQ_0 = 0$, for all $a$, we set $S = Q_0 Q_a + Q_a Q_0$ and evaluate it on the basis vectors.
\begin{equation*}
S \x= Q_0e_a  + Q_a e_0 = -e_{a+m} + e_{a+m} =0
\end{equation*}
\begin{equation*}
S e_0 = Q_0e_{a+m} + Q_a \x = -e_a + e_a =0
\end{equation*}
\begin{equation*}
\aligned
S e_b &= Q_0(\delta_{ab}\x+ L^b_{ac}e_{c+m} + F^\alpha_{a+m\,b} e_\alpha + F^\mu_{a+m\,b} e_\mu) 
+ Q_a (-e_{b+m}) \\
& = \delta_{ab}e_0 - L^b_{ac} e_c -F^\alpha_{a+m\,b} e_\alpha + F^\mu_{a+m\,b} e_\mu \\
&\quad -\delta_{ab}e_0 -L^c_{ab}e_c -F^\alpha_{b+m\,a}e_\alpha + F^\mu_{b+m\,a} e_\mu =0
\endaligned
\end{equation*}
\begin{equation*}
\aligned
S e_{b+m} &= Q_0(\delta_{ab} e_0 + L^c_{ab}e_c + F^\alpha_{b+m\,a}e_\alpha - F^\mu_{b+m\,a}e_\mu)  + Q_a(-e_b) \\
&= \delta_{ab} \x- L^c_{ab}e_{c+m} - F^\alpha_{b+m\,a}e_\alpha - F^\mu_{b+m\,a}e_\mu \\
&\quad - \delta_{ab}\x- L^b_{ac}e_{c+m} - F^\alpha_{a+m\,b}e_\alpha - F^\mu_{a+m\,b}e_\mu =0
\endaligned
\end{equation*}
\begin{equation*}
\aligned
S e_\alpha &= Q_0(F^\alpha_{a+m\,b}e_b + F^\alpha_{b+m\,a}e_{b+m} -2F^\mu_{\alpha a} e_\mu) + Q_a(- e_\alpha)\\
&= -F^\alpha_{a+m\,b}e_{b+m} -F^\alpha_{b+m\,a}e_b -2F^\mu_{\alpha a}e_\mu \\
&\quad -F^\alpha_{a+m\,b}e_b - F^\alpha_{b+m\,a}e_{b+m} + 2F^\mu_{\alpha a}e_\mu =0
\endaligned
\end{equation*}
\begin{equation*}
\aligned
S e_\mu &= Q_0(F^\mu_{a+m\,b}e_b - F^\mu_{b+m\,a}e_{b+m} -2F^\mu_{\alpha a} e_\alpha) + Q_ae_\mu\\
&= -F^\mu_{a+m\,b}e_{b+m} + F^\mu_{b+m\,a}e_b + 2F^\mu_{\alpha a} e_\alpha \\
&\quad + F^\mu_{a+m\,b}e_b - F^\mu_{b+m\,a}e_{b+m} - 2F^\mu_{\alpha a}e_\alpha =0
\endaligned
\end{equation*}
Therefore, $S=0$, which is what we wanted to prove.

Next we verify that $Q_aQ_d+Q_dQ_a = 2\delta_{ad}I$ for all $a$ and $d$.  For this verification we let $T = Q_aQ_d+Q_dQ_a$ and we evaluate it on the basis vectors.
\begin{equation*}
\aligned
T&\x = Q_a e_d+Q_d e_a = \delta_{ad}\x+ L^d_{ac} e_{c+m} + F^\alpha_{a+m\,d} e_\alpha + F^\mu_{a+m\,d} e_\mu \\
&+ \delta_{da} \x+ L^a_{dc} e_{c+m} + F^\alpha_{d+m\,a} e_\alpha + F^\mu_{d+m\,a} e_\mu = 2\delta_{ad} \x \\
T &e_0 = Q_a e_{d+m}+Q_d e_{a+m} = \delta_{ad} e_0 + L^c_{ad} e_c + F^\alpha_{d+m\,a} e_\alpha - F^\mu_{d+m\,a} e_\mu \\
&+ \delta_{da} e_0 + L^c_{da} e_c + F^\alpha_{a+m\,d} e_\alpha - F^\mu_{a+m\,d} e_\mu = 2\delta_{ad} e_0\\
T &e_b = Q_a (\delta_{db} \x+ L^b_{dc} e_{c+m} + F^\alpha_{d+m\,b} e_\alpha + F^\mu_{d+m\,b} e_\mu) \\
&+ Q_d (\delta_{ab} \x+ L^b_{ac} e_{c+m} + F^\alpha_{a+m\,b} e_\alpha + F^\mu_{a+m\,b}e_\mu) \\
&= (L^b_{da}+L^b_{ad})e_0 + (\delta_{bd}\delta_{ae} + \delta_{ba}\delta_{de} + L^b_{dc}L^e_{ac} + L^b_{ac}L^e_{dc} \\
&+ F^\alpha_{d+m\,b}F^\alpha_{a+m\,e} + F^\alpha_{a+m\,b}F^\alpha_{d+m\,e} + F^\mu_{d+m\,b}F^\mu_{a+m\,e} + F^\mu_{a+m\,b}F^\mu_{d+m\,e})e_e \\
&+(F^\alpha_{d+m\,b}F^\alpha_{e+m\,a} + F^\alpha_{a+m\,b}F^\alpha_{e+m\,d} - F^\mu_{d+m\,b}F^\mu_{e+m\,a} - F^\mu_{a+m\,b}F^\mu_{e+m\,d})e_{e+m} \\
&+ (L^b_{dc}F^\alpha_{c+m\,a} + L^b_{ac}F^\alpha_{c+m\,d} -2F^\mu_{d+m\,b}F^\mu_{\alpha a} -2 F^\mu_{a+m\,b} F^\mu_{\alpha d})e_\alpha \\
&- (L^b_{dc}F^\mu_{c+m\,a} + L^b_{ac}F^\mu_{c+m\,d} + 2F^\alpha_{d+m\,b}F^\mu_{\alpha a} + 2F^\alpha_{a+m\,b}F^\mu_{\alpha d})e_\mu \\
&= 2\delta_{ad}\delta_{be}e_e = 2\delta_{ad}e_b
\endaligned
\end{equation*}
where the coefficient of $e_e$ comes from~\eqref{eq:8:57}, the coefficient of $e_{e+m}$ is zero by the first equation of~\eqref{eq:5:6}, the coefficient of $e_\alpha$ is zero by~\eqref{eq:8:33} (with the roles of $b$ and $d$ reversed) and the coefficient of $e_\mu$ is zero by~\eqref{eq:8:34} (with the roles of $b$ and $d$ reversed).
\begin{equation*}
\aligned
T &e_{b+m} = (L^a_{db}+L^d_{ab}) \x \\
&+ (F^\alpha_{b+m\,d} F^\alpha_{a+m\,c} + F^\alpha_{b+m\,a} F^\alpha_{d+m\,c} - F^\mu_{b+m\,d}F^\mu_{a+m\,c} - F^\mu_{b+m\,a}F^\mu_{d+m\,c})e_c \\
&+ (\delta_{db}\delta_{ae} + \delta_{ab}\delta_{de} + L^c_{db}L^c_{ae} + L^c_{ab}L^c_{de} \\
&+ F^\alpha_{b+m\,d}F^\alpha_{e+m\,a} + F^\alpha_{b+m\,a}F^\alpha_{e+m\,d} + F^\mu_{b+m\,d}F^\mu_{e+m\,a} + F^\mu_{b+m\,a}F^\mu_{e+m\,d})e_{e+m} \\
&+ (L^c_{db}F^\alpha_{a+m\,c} + L^c_{ab}F^\alpha_{d+m\,c} + 2F^\mu_{b+m\,d}F^\mu_{\alpha a} + 2F^\mu_{b+m\,a} F^\mu_{\alpha d})e_\alpha \\
&+ (L^c_{db}F^\mu_{a+m\,c} + L^c_{ab}F^\mu_{d+m\,c} -2F^\alpha_{b+m\,d}F^\mu_{\alpha a} -2 F^\alpha_{b+m\,a}F^\mu_{\alpha d})e_\mu \\
&= 2\delta_{ad}\delta_{be}e_{e+m} = 2\delta_{ad}e_{b+m}
\endaligned
\end{equation*}
where the coefficient of $e_{e+m}$ comes from~\eqref{eq:8:57}, the coefficient of $e_c$ is zero by the first equation of~\eqref{eq:5:6}, the coefficient of $e_\alpha$ is zero by~\eqref{eq:8:33} (with the roles of $b$ and $d$ reversed) and the coefficient of $e_\mu$ is zero by~\eqref{eq:8:34} ( with the roles of $b$ and $d$ reversed).
\begin{equation*}
\aligned
T &e_\alpha = (F^\alpha_{d+m\,a} + F^\alpha_{a+m\,d}) \x + (F^\alpha_{a+m\,d} + F^\alpha_{d+m\,a}) e_0 \\
&+ (F^\alpha_{b+m\,d} L^c_{ab} + F^\alpha_{b+m\,a} L^c_{db} - 2F^\mu_{\alpha d} F^\mu_{a+m\,c} -2F^\mu_{\alpha a} F^\mu_{d+m\,c}) e_c \\
&+ (F^\alpha_{d+m\,b} L^b_{ac} + F^\alpha_{a+m\,b} L^b_{dc} + 2F^\mu_{\alpha d} F^\mu_{c+m\,a} + 2F^\mu_{\alpha a} F^\mu_{c+m\,d}) e_{c+m} \\
&+ 2(F^\alpha_{b+m\,d} F^\b_{b+m\,a} + F^\alpha_{b+m\,a} F^\b_{b+m\,d} + 2F^\mu_{\alpha d}F^\mu_{\b a} + 2F^\mu_{\alpha a} F^\mu_{\b d})e_\b \\
&+ (F^\alpha_{b+m\,d}F^\mu_{b+m\,a} + F^\alpha_{b+m\,a}F^\mu_{b+m\,d} - F^\alpha_{b+m\,d}F^\mu_{b+m\,a} - F^\alpha_{b+m\,a}F^\mu_{b+m\,d})e_\mu \\
&= 2\delta_{\alpha\b}\delta_{ad}e_\b = 2\delta_{ad}e_\alpha
\endaligned
\end{equation*}
where the coefficients of $\x$ and $e_0$ are clearly zero, the coefficients of $e_c$ and of $e_{c+m}$ are zero by~\eqref{eq:8:33} (in which the roles of $a,b,c,d$ are here played by $a,d,b,c$), the coefficient of $e_\b$ comes from the second equation of~\eqref{eq:5:6} and the coefficient of $e_\mu$ is clearly zero.
\begin{equation*}
\aligned
T &e_\mu = (F^\mu_{d+m\,a} + F^\mu_{a+m\, d})\x - (F^\mu_{a+m\,d} + F^\mu_{d+m\,a})e_0 \\
&- (2(F^\mu_{\alpha d}F^\alpha_{a+m\,b} + F^\mu_{\alpha a} F^\alpha_{d+m\, b})+ F^\mu_{c+m\,d}L^b_{ac} + F^\mu_{c+m\,a}L^b_{dc})e_b \\
&-(F^\mu_{d+m\,c}L^b_{ac} + 2F^\mu_{\alpha d}F^\alpha_{b+m\,a} + F^\mu_{a+m\,c} L^b_{dc} + 2F^\mu_{\alpha a} F^\alpha_{b+m\,d}) e_{b+m} \\
&+(F^\mu_{d+m\,b}F^\alpha_{a+m\,b} - F^\mu_{b+m\,d}F^\alpha_{b+m\,a} + F^\mu_{a+m\,b}F^\alpha_{d+m\,b} - F^\mu_{b+m\,a} F^\alpha_{b+m\,d})e_\alpha \\
&+ 2(F^\mu_{d+m\,b}F^\nu_{a+m\,b} + 2F^\mu_{\alpha d}F^\nu_{\alpha a} + F^\mu_{a+m\,b} F^\nu_{d+m\,b} + 2F^\mu_{\alpha a} F^\nu_{\alpha d})e_\nu \\
&= 2\delta_{ad}\delta_{\mu\nu} e_\nu = 2 \delta_{ad} e_\mu
\endaligned
\end{equation*}
where the coefficients of $\x$ and $e_0$ are clearly zero, the coefficients of $e_b$ and $e_{b+m}$ are zero by~\eqref{eq:8:34} (in which the roles of $b$ and $d$ are reversed), the coefficient of $e_\alpha$ is zero by~\eqref{eq:8:2} and~\eqref{eq:8:3}, and the coefficient of $e_\mu$ comes from the fifth equation of~\eqref{eq:5:6}.  This completes the proof of detail (I).  

In order to prove (II), we must find a Clifford system $P_0,\dots,P_m$ which is related to $Q_0,\dots,Q_m$ by~\eqref{eq:g:6}.  We do this by finding the map $B:U \to SO(m+1)$ of~\eqref{eq:g:5a}.  Let
\begin{equation}\label{eq:g:17}
\nu^a = \omega^a + \omega^{a+m}
\end{equation}
Use~\eqref{eq:5:1} together with~\eqref{eq:8:2}--\eqref{eq:8:4} to find
\begin{equation}\label{eq:g:18}
d\nu^a = - \nu^a_b \w \nu^b
\end{equation}
where
\begin{equation}\label{eq:g:19}
\nu^a_b = \theta^a_b + L^a_{cb} \omega^{c+m} + F^\alpha_{a+m\,b} \omega^\alpha + F^\mu_{a+m\,b} \omega^\mu = -\nu^b_a
\end{equation}
Set
\begin{equation}\label{eq:g:20}
\nu^0_b = -\nu^b_0 = -\nu^b = -(\omega^b+\omega^{b+m})
\end{equation}
We shall verify below that
\begin{equation}\label{eq:g:21}
dQ_j = \sum_{k=0}^m Q_k \nu^k_j, \mbox{ for $j = 0,\dots,m$}
\end{equation}
Differentiating this, we find that 
\begin{equation}\label{eq:g:22}
d\nu^i_j = -\sum_{k=0}^m \nu^i_k \w \nu^k_j, \mbox{ for $i,j =0,\dots m$}
\end{equation}
In fact, \eqref{eq:g:18} is the case $i=a$, $j=0$ and also implies the case $i=0$, $j=a$.  To verify the remaining cases in~\eqref{eq:g:22}, we take the exterior derivative of~\eqref{eq:g:21} when $j=a$, and then use~\eqref{eq:g:21} and~\eqref{eq:g:18} to find
\begin{equation}\label{eq:g:23}
\aligned
0 &= ddQ_a = dQ_0\w \nu^0_a + Q_0 d\nu^0_a + dQ_b \w \nu^b_a + Q_b d\nu^b_a \\
&= Q_b\nu^b \w \nu^0_a + Q_0 \nu^a_b \w \nu^b + (Q_0 \nu^0_b + Q_c \nu^c_b)\w \nu^b_a + Q_b d\nu^b_a \\
&= Q_b( d\nu^b_a + \nu^b_c \w \nu^c_a + \nu^b \w \nu^0_a) + Q_0(\nu^a_b \w \nu^b + \nu^0_b \w \nu^b_a)
\endaligned 
\end{equation}
which
implies~\eqref{eq:g:22} because the coefficient of $Q_0$ is zero and the $Q_b$ are linearly independent at each point of $U$, as can be seen from the fact that $Q_b \x=e_b$ are linearly independent at each point.
Define the $o(m+1)$-valued 1-form $\nu$ to be
\begin{equation}\label{eq:g:24}
\nu = \begin{pmatrix} 0 & \nu^0_b \\ \nu^a_0 & \nu^a_b \end{pmatrix}
\end{equation}
Then~\eqref{eq:g:18} and~\eqref{eq:g:22} imply that $d\nu =-\nu \w \nu$.  Therefore, on a simply connected subset of $U$, which we continue to call $U$, there exists a smooth map 
\begin{equation}\label{eq:g:25}
A:U\to SO(m+1)
\end{equation}
such that $A^{-1}dA = \nu$.  Denote the entries of $A$ by the functions $A^i_j$, $i,j = 0,\dots,m$, so that the entries of $dA = A \nu$ are given by
\begin{equation}\label{eq:g:26}
dA^i_j = \sum_{k=0}^m A^i_k \nu^k_j
\end{equation}
Let
\begin{equation}\label{eq:g:27}
P_i = \sum_{j=0}^m A^i_j Q_j, \mbox{ for $i=0,\dots,m$}
\end{equation}
which, at each point of $U$, is a set of symmetric, orthogonal transformations of $\R^{n+1}$ satisfying the conditions $P_iP_j+P_jP_i = 2\delta_{ij}I$, since $Q_iQ_j + Q_jQ_i =2\delta_{ij} I$ and $A\in SO(m+1)$.
Using~\eqref{eq:g:21} and~\eqref{eq:g:26}, we have
\begin{equation}\label{eq:g:28}
dP_i = \sum_{j=0}^m((dA^i_j)Q_j + A^i_j dQ_j) = \sum_{j,k=0}^m(A^i_k Q_j\nu^k_j + A^i_jQ_k\nu^k_j) =0
\end{equation}
since $\nu^i_j+\nu^j_i=0$.  Therefore, each $P_i$ is constant on $U$ and $P_0,\dots,P_m$ define a Clifford system on $\R^{n+1}$ and~\eqref{eq:g:6} holds with $B=A^{-1}$.

All that remains of the proof of detail (II) is to verify~\eqref{eq:g:21}, for which we need the Maurer-Cartan equations~\eqref{eq:4:15} for our Darboux frame field.
We first verify~\eqref{eq:g:21} for $j=0$, then for $j=a$, in both cases by evaluating each side on the basis vectors.
Differentiating equations in~\eqref{eq:g:3} and using~\eqref{eq:4:15}, we get
\begin{equation*}  
\aligned
(dQ_0&)\x = d(Q_0\x) - Q_0\,d\x\\
&= de_0 -Q_0(\omega^{a+m} e_{a+m} + \omega^\alpha e_\alpha + \omega^\mu e_\mu) \\
&= \omega^a e_a - \omega^\alpha e_\alpha +\omega^\mu e_\mu + \omega^{a+m} e_a + \omega^\alpha e_\alpha - \omega^\mu e_\mu \\ 
&= \nu^a e_a = \nu^aQ_a\x
\endaligned
\end{equation*}
\begin{equation*}  
\aligned
(dQ_0&)e_0 = d(Q_0 e_0) - Q_0\, de_0 \\
&= d\x- Q_0(\omega^a e_a - \omega^\alpha e_\alpha + \omega^\mu e_\mu) \\
&= \nu^a e_{a+m} = \nu^a Q_ae_0
\endaligned
\end{equation*}
\begin{equation*}  
\aligned
(dQ_0&)e_a = d(Q_0e_a)-Q_0de_a = -de_{a+m} - Q_0\, de_a \\
&= \nu^a\x+ (\omega^{b+m}_a - \omega^b_{a+m}) e_b + (\omega^b_a - \omega^{b+m}_{a+m})e_{b+m} \\
&\quad + (\omega^\alpha_a - \omega^\alpha_{a+m} ) e_\alpha - (\omega^\mu_{a+m} + \omega^\mu_a) e_\mu \\
&= \nu^b(\delta_{ab} \x+ L^c_{ab} e_{c+m} + F^\alpha_{b+m\,a} e_\alpha + F^\mu_{b+m\,a} e_\mu) \\
&= \nu^bQ_be_a 
\endaligned
\end{equation*}
\begin{equation*}  
\aligned
(dQ_0&)e_{a+m} = d(Q_0e_{a+m}) - Q_0\, de_{a+m} = -de_a - Q_0\,de_{a+m} \\
&= (\omega^a+\omega^{a+m})e_0 +(\omega^{b+m}_{a+m} - \omega^b_a)e_b + (\omega^b_{a+m} - \omega^{b+m}_a)e_{b+m}\\
&\quad + (\omega^\alpha_{a+m}- \omega^\alpha_a)e_\alpha - (\omega^\mu_{a+m} +\omega^\mu_a)e_\mu \\
&= \nu^b(\delta_{ab}e_0 - L^c_{ab}e_c + F^\alpha_{a+m\,b} e_\alpha -F^\mu_{a+m\,b}e_\mu) = \nu^b Q_be_{a+m} 
\endaligned
\end{equation*}
\begin{equation*}  
\aligned
(dQ_0&)e_\alpha = d(Q_0 e_\alpha) - Q_0\,de_\alpha = -de_\alpha -Q_0\,de_\alpha \\
&= (\omega^{a+m}_\alpha - \omega^a_\alpha)e_a + (\omega^a_\alpha - \omega^{a+m}_\alpha) e_{a+m}- 2\omega^\mu_\alpha e_\mu \\
&= \nu^b (F^\alpha_{b+m\,a}e_a - F^\alpha_{b+m\,a}e_{a+m} -2F^\mu_{\alpha b} e_\mu) = \nu^b Q_be_\alpha
\endaligned
\end{equation*}
\begin{equation*}  
\aligned
(dQ_0&)e_\mu = d(Q_0e_\mu) - Q_0\,de_\mu = de_\mu - Q_0\,de_\mu \\
&= (\theta^a_\mu+ \theta^{a+m}_\mu)(e_a+e_{a+m}) + 2\theta^\alpha_\mu e_\alpha \\
&= \nu^b(F^\mu_{b+m\,a}(e_a+ e_{a+m}) -2 F^\mu_{\alpha b}e_\alpha) = \nu^bQ_b e_\mu
\endaligned
\end{equation*}
That completes the verification of~\eqref{eq:g:21} for the case $j=0$.

We now verify the equations in~\eqref{eq:g:21} for the cases $j = a$ by applying each side to the basis vectors.
Using~\eqref{eq:g:4} and~\eqref{eq:4:15}, we have
\begin{equation*}  
\aligned
&(dQ_a)\x = d(Q_a \x) - Q_a d\x = de_a - Q_a(\omega^{b+m}e_{b+m} + \omega^\alpha e_\alpha + \omega^\mu e_\mu) \\
&= -\nu^a e_0 +(\theta^b_a - L^b_{ac}\omega^{c+m} -F^\alpha_{a+m\,b}\omega^\alpha - F^\mu_{a+m\,b} \omega^\mu) e_b\\
&+ (\theta^{b+m}_a - F^\alpha_{b+m\,a}\omega^\alpha - F^\mu_{b+m\,a}\omega^\mu)e_{b+m} \\
&+ (\theta^\alpha_a - F^\alpha_{b+m\,a}\omega^{b+m} + 2F^\mu_{\alpha a} \omega^\mu) e_\alpha
+ (\theta^\mu_a + F^\mu_{b+m\,a} \omega^{b+m} + 2 F^\mu_{\alpha a} \omega^\alpha) e_\mu \\
&= -\nu^a_0 e_0 + \nu^b_a e_b = (\nu^0_aQ_0 + \nu^b_a Q_b)\x
\endaligned
\end{equation*}
\begin{equation*}  
\aligned
(dQ_a&)e_0 = d(Q_ae_0) - Q_ade_0 = de_{a+m} -Q_a(\omega^be_b - \omega^\alpha e_\alpha + \omega^\mu e_\mu) \\
&= (-\omega^{a+m} - \omega^a) \x + (\theta^b_{a+m} + F^\alpha_{a+m\, b} \omega^\alpha - F^\mu_{a+m\,b} \omega^\mu)e_b \\
&+ (\theta^{b+m}_{a+m} - L^c_{ab} \omega^c + F^\alpha_{b+m\,a} \omega^\alpha + F^\mu_{b+m\,a}\omega^\mu) e_{b+m} \\
&+ (\theta^\alpha_{a+m} - F^\alpha_{a+m\,b} \omega^b + 2F^\mu_{\alpha a} \omega^\mu) e_\alpha \\
&+ (\theta^\mu_{a+m} - F^\mu_{a+m\,b} \omega^b - 2F^\mu_{\alpha a} \omega^\alpha) e_\mu \\
&= -\nu^a \x + (\theta^{b+m}_{a+m} - L^c_{ab}\omega^c + F^\alpha_{b+m\,a}\omega^\alpha + F^\mu_{b+m\,a}\omega^\mu)e_{b+m} \\
&= -\nu^a\x+ \nu^b_ae_{b+m} = (\nu^0_a Q_0 + \nu^b_a Q_b)e_0
\endaligned
\end{equation*}
where the coefficients of $e_b, e_\alpha$ and $e_\mu$ are zero by~\eqref{eq:4:13}, and~\eqref{eq:g:19} is used in the coefficient of $e_{b+m}$.

In order to verify~\eqref{eq:g:21} when both sides are applied to $e_b$, we must verify that
\begin{equation}\label{eq:g:37}
\aligned
&d(Q_a e_b) - Q_a de_b = (dQ_a)e_b = \nu^0_a Q_0 e_b + \nu^c_a Q_c e_b \\
&= \nu^b_a \x + (\delta_{bd}\nu^a + L^d_{bc} \nu^c_a) e_{d+m}
+ F^\a_{c+m\,b} \nu^c_a e_\alpha + F^\mu_{c+m\,b} \nu^c_a e_\mu
\endaligned
\end{equation}
Using~\eqref{eq:g:4} and~\eqref{eq:4:15}, and gathering together the coefficients of each basis vector, we get
\begin{equation}\label{eq:g:dQaeb}
\aligned
&d(Q_a e_b) - Q_a de_b = \\
&\quad (L^c_{ab} \omega^{c+m} - F^\alpha_{a+m\,b} \omega^\alpha - F^\mu_{a+m\,b} \omega^\mu - \theta^a_b) \x \\
&+ (F^\alpha_{a+m\,b} \omega^\alpha - F^\mu_{a+m\,b} \omega^\mu - \theta^{a+m}_b) e_0  \\
&+(-L^c_{ab} \theta^d_{c+m} + F^\alpha_{a+m\,b} \theta^d_\alpha + F^\mu_{a+m\,b} \theta^d_\mu \\
&\quad - L^d_{ac} \theta^{c+m}_b - F^\alpha_{a+m\,d} \theta^\alpha_b - F^\mu_{a+m\,d} \theta^\mu_b) e_d \\
&+ ( \delta_{ab} \omega^{c+m} - dL^c_{ab} - L^d_{ab} \theta^{c+m}_{d+m} + F^\alpha_{a+m\,b} \theta^{c+m}_\alpha + F^\mu_{a+m\,b} \theta^{c+m}_\mu \\
&\quad+ \delta_{ac} \omega^b + 
L^c_{ad} \theta^d_b -F^\alpha_{c+m\,a} \theta^\alpha_b + F^\mu_{c+m\,a} \theta^\mu_b) e_{c+m} \\
&+ (\delta_{ab} \omega^\alpha - L^c_{ab} \theta^\alpha_{c+m} + dF^\alpha_{a+m\,b}
+ F^\beta_{a+m\, b} \theta^\alpha_\beta \\
&\quad + F^\mu_{a+m\, b} \theta^\alpha_\mu - F^\alpha_{a+m\,c} \theta^c_b 
- F^\alpha_{c+m\,a} \theta^{c+m}_b + 2 F^\mu_{\alpha a} \theta^\mu_b ) e_\alpha \\
&+ (\delta_{ab} \omega^\mu - L^c_{ab} \theta^\mu_{c+m} + F^\alpha_{a+m\,b} \theta^\mu_\alpha + dF^\mu_{a+m\, b} \\
&\quad + F^\nu_{a+m\,b} \theta^\mu_\nu - F^\mu_{a+m\,c} \theta^c_b 
+ F^\mu_{c+m\,a} \theta^{c+m}_b + 2F^\mu_{\alpha a} \theta^\alpha_b) e_\mu
\endaligned
\end{equation}
The coefficient of $\x$ is
$\nu^b_a$
by~\eqref{eq:g:19}.  The coefficient of $e_0$ is $0$
Substituting~\eqref{eq:4:13} into the coefficient of $e_d$, we get
\begin{equation*}  
\aligned
&(-F^\alpha_{a+m\,b}F^\alpha_{c+m\,d} - F^\alpha_{a+m\,d}F^\alpha_{c+m\,b}
+ F^\mu_{a+m\,b} F^\mu_{c+m\,d} + F^\mu_{a+m\,d}F^\mu_{c+m\,b}
)\omega^{c+m} \\ 
&- (L^b_{ac}F^\alpha_{c+m\,d} + L^d_{ac}F^\alpha_{c+m\,b}
-2F^\mu_{a+m\,b}F^\mu_{\alpha d} - 2F^\mu_{a+m\,d}F^\mu_{\alpha
b})\omega^\alpha \\ 
&+ (L^b_{ac}F^\mu_{c+m\,d} + L^d_{ac}F^\mu_{c+m\,b}
+2F^\alpha_{a+m\,b}F^\mu_{\alpha d} + 2F^\alpha_{a+m\,d}F^\mu_{\alpha
b})\omega^\mu 
\endaligned
\end{equation*}
which is zero since the coefficient of $\omega^{c+m}$ is zero by the
first equation in~\eqref{eq:5:6}, the coefficient of $\omega^\alpha$
is zero by~\eqref{eq:8:33} and the coefficient of $\omega^\mu$ is zero
by~\eqref{eq:8:34}.  Thus, the coefficient of $e_d$ is zero, in
agreement with the right hand side of~\eqref{eq:g:37}. 

By~\eqref{eq:4:13}, \eqref{eq:8:52} and~\eqref{eq:g:17} with~\eqref{eq:g:19}, the coefficient of $e_{c+m}$ becomes
\begin{equation*}  
\aligned
&-L^c_{db}\nu^d_a + \nu^a \delta^c_b \\
&- (L^c_{abd} - L^e_{ab}L^c_{ed} + F^\alpha_{a+m\,b} F^\alpha_{c+m\, d} +
F^\mu_{a+m\,b}F^\mu_{c+m\, d} - \delta_{ac} \delta_{bd} + \delta_{bc}\delta_{ad}) \omega^d \\
&+( - L^c_{ab\,d+m} + L^c_{eb} L^e_{da} + L^e_{ab} L^c_{ed} -F^\alpha_{c+m\,a} F^\alpha_{d+m\,b} \\
&\qquad - F^\mu_{c+m\,a} F^\mu_{d+m\, b} + \delta_{ab} \delta_{cd} - \delta_{bc} \delta_{ad}) \omega^{d+m} \\
&+(-L^c_{ab\alpha} + L^c_{db}F^\alpha_{d+m\,a} -2F^\mu_{a+m\,b} F^\mu_{\alpha\, c+m} - 2 F^\mu_{c+m\,a} F^\mu_{\alpha b}) \omega^\alpha \\
&+( -L^c_{ab\mu} + L^c_{db} F^\mu_{d+m\,a} + 2F^\alpha_{a+m\,b} F^\mu_{\alpha\, c+m} + 2F^\alpha_{c+m\,a} F^\mu_{\alpha b}) \omega^\mu
\endaligned
\end{equation*}
We now verify that zero is the coefficient of each of $\omega^d$, $\omega^{d+m}$, $\omega^\alpha$, $\omega^\mu$.  

The coefficient of $\omega^d$ can be seen to be zero by taking~\eqref{eq:8:53} (with indices in the order $c,a,b,d$) and subtracting half of~\eqref{eq:8:57} (with indices as is).

The coefficient of $\omega^{d+m}$ can be seen to be zero by using~\eqref{eq:8:62}, then adding~\eqref{eq:8:53} (with indices in the order $c,a,b,d$), then adding half of~\eqref{eq:8:57} (with indices as is), and then using~\eqref{eq:8:57} again (with the roles of $d$ and $c$ reversed).

The coefficient of $\omega^\alpha$ can be seen to be zero from~\eqref{eq:8:55} and~\eqref{eq:8:33}.

The coefficient of $\omega^\mu$ is zero by~\eqref{eq:8:56}.

Hence, we have shown that the coefficient of $e_{c+m}$ in $(dQ_a)e_b$ is as given in~\eqref{eq:g:37}.

Using~\eqref{eq:5:2}, \eqref{eq:8:4} and~\eqref{eq:g:19}, we can rewrite the coefficient of $e_\alpha$ in $(dQ_a)e_b$ in~\eqref{eq:g:dQaeb} as
\begin{equation*}  
\aligned
&F^\alpha_{c+m\,b} \nu^c_a \\
&+(F^\alpha_{a+m\, bc} + L^b_{ad}F^\alpha_{d+m\,c} - F^\mu_{a+m\,b} F^\mu_{\alpha c} - F^\alpha_{d+m\,b}L^d_{ac}) \omega^c \\
&+ (F^\alpha_{a+m\,b\,c+m} - F^\mu_{a+m\,b} F^\mu_{\alpha\,c+m} - 2F^\mu_{\alpha a} F^\mu_{c+m\,b})\omega^{c+m} \\
&+(F^\alpha_{a+m\,b \beta} + \delta_{ab}\delta_{\alpha\beta} - F^\alpha_{c+m\,a} F^\beta_{c+m\,b} - 4 F^\mu_{\alpha a} F^\mu_{\beta b} - F^\alpha_{c+m\, b} F^\beta_{c+m\,a}) \omega^\beta \\
&+ (F^\alpha_{a+m\,b\mu} - 2L^b_{ac} F^\mu_{\alpha c} + F^\alpha_{c+m\,a} F^\mu_{c+m\,b} - F^\alpha_{c+m\,b} F^\mu_{c+m\,a}) \omega^\mu
\endaligned
\end{equation*}
The coefficient of $\omega^c$ is seen to be zero by using the first equation in~\eqref{eq:5:7} and then using~\eqref{eq:8:33}.  The coefficient of $\omega^{c+m}$ is zero by the second equation in~\eqref{eq:5:7}.  The coefficient of $\omega^\beta$ is seen to be zero by using the third equation in~\eqref{eq:5:7} and the using the second equation in~\eqref{eq:5:6}.  The coefficient of $\omega^\mu$ is seen to be zero by using~\eqref{eq:5:11} and then using~\eqref{eq:8:35} (with the roles of $a$ and $b$ interchanged).

Using~\eqref{eq:5:2} and~\eqref{eq:g:19}, we can rewrite the coefficient of $e_\mu$ in $(dQ_a)e_b$ in~\eqref{eq:g:dQaeb} as
\begin{equation*}   
\aligned
&F^\mu_{c+m\, b} \nu^c_a +\\
&(F^\mu_{a+m\,bc} + L^b_{ad}F^\mu_{d+m\,c} + F^\alpha_{a+m\,b} F^\mu_{\alpha c} - F^\mu_{d+m\,b} L^d_{ac})\omega^c \\
&+ (F^\mu_{a+m\,b\,c+m} + F^\alpha_{a+m\,b} F^\mu_{\alpha\,c+m} + 2F^\mu_{\alpha a} F^\alpha_{c+m\,b}) \omega^{c+m} \\
&+(F^\mu_{a+m\,b\alpha} + 2L^b_{ac} F^\mu_{\alpha\,c+m} + F^\mu_{c+m\,a} F^\alpha_{c+m\,b} - F^\mu_{c+m\,b} F^\alpha_{c+m\,a}) \omega^\alpha \\
&+(F^\mu_{a+m\,b\nu} + \delta_{ab} \delta_{\mu\nu} - F^\mu_{c+m\,a} F^\nu_{c+m\,b} - 4 F^\mu_{\alpha a} F^\nu_{\alpha b} - F^\mu_{c+m\,b} F^\nu_{c+m\,a}) \omega^\nu
\endaligned
\end{equation*}
The coefficient of $\omega^c$ is seen to be zero by using the first equation in~\eqref{eq:5:8} and then using~\eqref{eq:8:34}.  The coefficient of $\omega^{c+m}$ is zero by the second equation in~\eqref{eq:5:8}.  The coefficient of $\omega^\alpha$ is zero by~\eqref{eq:8:35}.  The coefficient of $\omega^\mu$ is seen to be zero by using the third equation in~\eqref{eq:5:8} and then using the fourth equation in~\eqref{eq:5:6}.  This completes the verification of~\eqref{eq:g:37}.

The next case is to verify~\eqref{eq:g:21} when both sides are applied to $e_{b+m}$.  We must verify that
\begin{equation}\label{eq:g:55a}
\aligned
&d(Q_a e_{b+m}) - Q_a d e_{b+m} = (dQ_a) e_{b+m} = \nu^0_a Q_0 e_{b+m} + \nu^c_a Q_c e_{b+m} \\
&=\nu^a e_b + \nu^b_a e_0 + L^c_{db} \nu^d_a e_c + F^\alpha_{b+m\,c} \nu^c_a e_\alpha - F^\mu_{b+m\,c} \nu^c_a e_\mu
\endaligned
\end{equation}
by~\eqref{eq:g:4}.
Using~\eqref{eq:g:4} to compute $Q_a e_{b+m}$ and~\eqref{eq:4:15} to compute $d e_{b+m}$, the left hand side becomes
\begin{equation}\label{eq:g:42}
\aligned
&(-F^\alpha_{b+m\,a} \omega^\alpha + F^\mu_{b+m\,a} \omega^\mu - \theta^a_{b+m}) \x\\
&+ (-L^c_{ab}\omega^c + F^\alpha_{b+m\, a} \omega^\alpha + F^\mu_{b+m\,a} \omega^\mu - \theta^{a+m}_{b+m})e_0 \\
&+ (\delta_{ab} \omega^c + dL^c_{ab} + L^d_{ab} \theta^c_d + F^\alpha_{b+m\,a} \theta^c_\alpha - F^\mu_{b+m\,a} \theta^c_\mu \\
&\qquad \quad+ \delta_{ac} \omega^{b+m} - L^c_{ad} \theta^{d+m}_{b+m} - F^\alpha_{a+m\,c} \theta^\alpha_{b+m} - F^\mu_{a+m\,c} \theta^\mu_{b+m})e_c \\
&+  ( L^c_{ab} \theta^{d+m}_c + F^\alpha_{b+m\,a} \theta^{d+m}_\alpha - F^\mu_{b+m\,a} \theta^{d+m}_\mu - L^c_{ad} \theta^c_{b+m} \\
&\qquad \quad- F^\alpha_{d+m\,a} \theta^\alpha_{b+m} + F^\mu_{d+m\,a} \theta^\mu_{b+m}) e_{d+m}\\
&+ ( -\delta_{ab} \omega^\alpha + L^c_{ab} \theta^\alpha_c + dF^\alpha_{b+m\,a} + F^\beta_{b+m\,a} \theta^\alpha_\beta - F^\mu_{b+m\,a} \theta^\alpha_\mu \\
&\qquad \quad- F^\alpha_{a+m\,c} \theta^c_{b+m} - F^\alpha_{c+m\,a} \theta^{c+m}_{b+m} + 2F^\mu_{\alpha a} \theta^\mu_{b+m})e_\alpha \\
&+ ( \delta_{ab} \omega^\mu + L^c_{ab} \theta^\mu_c + F^\alpha_{b+m\,a} \theta^\mu_\alpha - dF^\mu_{b+m\,a} - F^\nu_{b+m\,a} \theta^\mu_\nu \\
&\qquad \quad- F^\mu_{a+m\,c} \theta^c_{b+m} + F^\mu_{c+m\,a} \theta^{c+m}_{b+m} + 2F^\mu_{\alpha a} \theta^\alpha_{b+m})e_\mu
\endaligned
\end{equation}
We want to verify that this is equal to the right side of~\eqref{eq:g:55a}, where $\nu^c_a$ is given by~\eqref{eq:g:19}.
We do this by comparing the coefficients of the basis vectors $\x, e_0, e_c, e_{c+m}, e_\alpha, e_\mu$.

The coefficient of $\x$ is 0 by~\eqref{eq:4:13}.  

The coefficient of $e_0$ is 
\begin{equation*}
\aligned
&\theta^b_a + L^b_{ca} \omega^{c+m} + F^\alpha_{b+m\,a} \omega^\alpha + F^\mu_{b+m\,a} \omega^\mu - L^b_{ca}(\omega^c + \omega^{c+m}) +(\theta^{b+m}_{a+m} - \theta^b_a) \\
&= \nu^b_a
\endaligned
\end{equation*}
by~\eqref{eq:8:4}, \eqref{eq:g:17} and~\eqref{eq:g:19}.

The coefficient of $e_c$ in $(dQ_a)e_{b+m}$ in~\eqref{eq:g:42} is, using~\eqref{eq:4:13} and~\eqref{eq:8:52} and the skew symmetry of $L^a_{bcd}$ in all four indices,
\begin{equation*}  
\aligned
&L^c_{db}(\theta^d_a + L^d_{ea} \omega^{e+m} + F^\alpha_{d+m\,a} \omega^\alpha + F^\mu_{d+m\,a} \omega^\mu) + (\omega^a + \omega^{a+m})\delta_{bc} \\
&+ (L^a_{bcd} + \delta_{ab}\delta_{cd} - \delta_{ad}\delta_{bc} +  L^c_{ae}L^e_{bd} - F^\alpha_{a+m\,c} F^\alpha_{b+m\,d} - F^\mu_{a+m\,c} F^\mu_{b+m\,d}) \omega^d \\
&+ (L^c_{ab\,d+m} + \delta_{ac} \delta_{bd} - \delta_{ad}\delta_{bc} + L^c_{eb}L^e_{ad} + L^c_{ae}L^e_{bd} \\
&\qquad\quad - F^\alpha_{b+m\,a}F^\alpha_{d+m\,c} - F^\mu_{b+m\,a} F^\mu_{d+m\,c}) \omega^{d+m}\\
&+ ( L^c_{ab\alpha} - L^c_{eb} F^\alpha_{e+m\,a} -2F^\mu_{b+m\,a} F^\mu_{\alpha c} - 2 F^\mu_{a+m\,c} F^\mu_{\alpha b+m} ) \omega^\alpha \\
&+ (L^c_{ab\mu} - L^c_{eb} F^\mu_{e+m\,a} + 2F^\alpha_{b+m\,a} F^\mu_{\alpha c} + 2 F^\alpha_{a+m\,c} F^\mu_{\alpha b+m} ) \omega^\mu \\
&= L^c_{db} \nu^d_a + \nu^a \delta_{bc}
\endaligned
\end{equation*}
by~\eqref{eq:g:17} and~\eqref{eq:g:19} and the following.
Using~\eqref{eq:8:41} for $L^a_{bcd}$, the coefficient of $\omega^d$ becomes
\begin{equation*} 
\aligned
&\delta_{ab}\delta_{cd} -\frac12(\delta_{ac}\delta_{bd} + \delta_{ad}\delta_{bc}) + \frac12(L^c_{ae} L^e_{bd} - L^a_{de}L^e_{bc}) \\
&+ F^\mu_{c+m\,a}F^\mu_{b+m\,d} - F^\mu_{c+m\,b}F^\mu_{d+m\,a}
\endaligned
\end{equation*}
which is zero by~\eqref{eq:8:45} combined with the first equation of~\eqref{eq:5:6}.
In the coefficient of $\omega^{d+m}$, substitute~\eqref{eq:8:42} for $L^a_{bc\,d+m} = L^c_{ab\,d+m}$, and gather together terms using skew symmetries, to get
\begin{equation*}  
\aligned
&\quad \frac32(\delta_{ac}\delta_{bd} - \delta_{ad}\delta_{bc}) + L^a_{be} L^d_{ce} + \frac12 L^a_{de}L^e_{bc} +\frac12 L^a_{ce} L^d_{be} \\
& - F^\alpha_{a+m\,d} F^\alpha_{c+m\,b} - F^\alpha_{a+m\,b} F^\alpha_{c+m\,d} +F^\mu_{d+m\,b} F^\mu_{a+m\,c} + F^\mu_{a+m\,b} F^\mu_{d+m\,c} 
\endaligned
\end{equation*}
which is zero by using the first equation in~\eqref{eq:5:6} and then~\eqref{eq:8:45}.
The coefficient of $\omega^\alpha$ is zero by~\eqref{eq:8:43}.  The coefficient of $\omega^\mu$ is zero by~\eqref{eq:8:44}.

The coefficient of $e_{d+m}$ in $(dQ_a)e_{b+m}$ in~\eqref{eq:g:42} is, using~\eqref{eq:4:13}
\begin{equation*} 
\aligned
&(-F^\alpha_{b+m\,a} F^\alpha_{d+m\,c} + F^\mu_{b+m\,a} F^\mu_{d+m\,c}
-F^\alpha_{d+m\,a} F^\alpha_{b+m\,c} + F^\mu_{d+m\,a} F^\mu_{b+m\,c})
\omega^c \\ 
&+ (L^c_{ba} F^\alpha_{c+m\,d} - 2F^\mu_{a+m\,b} F^\mu_{\alpha d} +
L^c_{da} F^\alpha_{c+m\,b} - 2F^\mu_{a+m\,d} F^\mu_{\alpha\, b+m})
\omega^\alpha \\ 
&+ (-L^c_{ba} F^\mu_{c+m\,d} + 2F^\alpha_{b+m\,a} F^\mu_{\alpha\, d+m} -
L^c_{da} F^\mu_{c+m\,b} + 2F^\alpha_{d+m\,a} F^\mu_{\alpha\,
b+m})\omega^\mu \\ 
&= 0
\endaligned
\end{equation*}
because the coefficient of $\omega^c$ is 0 by the first equation
in~\eqref{eq:5:6}, the coefficient of $\omega^\alpha$ is 0
by~\eqref{eq:8:33} and~\eqref{eq:8:1}, and the coefficient of
$\omega^\mu$ is 0 
by~\eqref{eq:8:34} and~\eqref{eq:8:1}. 

The coefficient of $e_\alpha$ in $(dQ_a)e_{b+m}$ in~\eqref{eq:g:42}
is, using~\eqref{eq:4:13} and~\eqref{eq:5:2} 
\begin{equation*}  
\aligned
&F^\alpha_{b+m\,c}(\theta^c_a + L^c_{da} \omega^{d+m} +
F^\beta_{c+m\,a} \omega^\beta + F^\mu_{c+m\,a} \omega^\mu) \\ 
&+  (F^\alpha_{b+m\,ac} + F^\mu_{b+m\,a} F^\mu_{\alpha c} +
2F^\mu_{\alpha a} F^\mu_{b+m\,c}) \omega^c \\ 
&+ (F^\alpha_{b+m\,a \, d+m} + L^c_{ba}F^\alpha_{c+m\,d}
+ L^c_{da}F^\alpha_{c+m\,b} 
 + F^\mu_{b+m\,a} F^\mu_{\alpha\, d+m}) \omega^{d+m}\\
& + (-\delta_{ab} \delta_{\alpha\b} + F^\alpha_{b+m\,a\b} +
F^\alpha_{a+m\,c} F^\b_{b+m\,c} - F^\alpha_{b+m\,c} F^\beta_{c+m\,a}  
+ 4F^\mu_{\alpha a} F^\mu_{\b\, b+m} ) \omega^\b \\
&+ (-2L^c_{ab} F^\mu_{\alpha c} + F^\alpha_{b+m\, a\mu} -
F^\alpha_{a+m\,c} F^\mu_{b+m\,c} - F^\alpha_{b+m\,c} F^\mu_{c+m\,a}
)\omega^\mu \\ 
&= F^\alpha_{b+m\,c} \nu^c_a
\endaligned
\end{equation*}
by~\eqref{eq:g:19}, because the other terms are zero as follows.
The coefficient of $\omega^c$ is zero by the first equation
in~\eqref{eq:5:7}.  The coefficient of $\omega^{d+m}$ is zero by the
second equation in~\eqref{eq:5:7} and~\eqref{eq:8:33}.  The
coefficient of $\omega^\b$ is zero by the third equation
of~\eqref{eq:5:7} and the third equation of~\eqref{eq:5:6}.  The
coefficient of $\omega^\mu$ is zero by the third equation in~\eqref{eq:5:11}
and~\eqref{eq:8:35}. 

Finally, the coefficient of $e_\mu$ in $(dQ_a)e_{b+m}$ in~\eqref{eq:g:42} is, using~\eqref{eq:4:13} and~\eqref{eq:5:2}
\begin{equation*}\label{eq:g:53}
\aligned
& -F^\mu_{b+m\,c}(\theta^c_a - L^c_{ad} \omega^{d+m} + F^\alpha_{c+m\,a}\omega^\alpha + F^\nu_{c+m\,a} \omega^\nu)\\
&+(F^\alpha_{b+m\,a} F^\mu_{\alpha c} - F^\mu_{b+m\,ac} + 2F^\mu_{\alpha a} F^\alpha_{b+m\,c}) \omega^c \\
&+ ( -L^c_{ab} F^\mu_{d+m\,c} -L^c_{ad} F^\mu_{b+m\,c}
+ F^\alpha_{b+m\,a} F^\mu_{\alpha\,d+m} - F^\mu_{b+m\,a\,d+m}) \omega^{d+m} \\
&+ ( -2L^c_{ab} F^\mu_{\alpha c} - F^\mu_{b+m\,a\alpha} + F^\mu_{a+m\,c} F^\alpha_{b+m\,c} + F^\mu_{b+m\,c} F^\alpha_{c+m\,a} ) \omega^\alpha \\
&+ ( \delta_{ab} \delta_{\mu\nu} - F^\mu_{b+m\,a\nu} - F^\mu_{a+m\,c}
F^\nu_{b+m\,c} + F^\mu_{b+m\,c} F^\nu_{c+m\,a} - 4F^\mu_{\alpha a}
F^\nu_{\alpha\,b+m}) \omega^\mu \\ 
&= -F^\mu_{b+m\,c}\nu^c_a
\endaligned
\end{equation*}
by~\eqref{eq:g:19}, because the other terms are zero as follows.
The coefficient of $\omega^c$ is zero by the first equation
in~\eqref{eq:5:8}.  The coefficient of $\omega^{d+m}$ is zero by the
second equation in~\eqref{eq:5:8} and by~\eqref{eq:8:34}.  The
coefficient of $\omega^\alpha$ is zero by~\eqref{eq:5:11}
and~\eqref{eq:8:35}.  The coefficient of $\omega^\nu$ is zero by the
fourth equation in~\eqref{eq:5:6}. 

That concludes the verification of~\eqref{eq:g:55a}.

The next case is to verify~\eqref{eq:g:21} when both sides are applied to $e_\alpha$.  We must verify that
\begin{equation}\label{eq:g:55}
\aligned
&d(Q_a e_\alpha) - Q_a d e_\alpha = (dQ_a) e_\alpha = \nu^0_a Q_0 e_\alpha + \nu^c_a Q_c e_\alpha \\
&=\nu^a e_\alpha + F^\alpha_{b+m\,c}\nu^b_a e_c + F^\alpha_{c+m\,b} \nu^b_a e_{c+m} -2F^\mu_{\alpha b}\nu^b_a e_\mu
\endaligned
\end{equation}
Using~\eqref{eq:g:4} and~\eqref{eq:4:15}, and gathering together the coefficients of each basis vector, we get
\begin{equation}\label{eq:g:dQaealpha}
\aligned
&(dQ_a)e_\alpha = (-F^\alpha_{b+m\, a} \omega^{b+m} + 2F^\mu_{\alpha a} \omega^\mu - \theta^a_\alpha)\x \\
&+ (-F^\alpha_{a+m\,b} \omega^b + 2F^\mu_{\alpha a} \omega^\mu - \theta^{a+m}_\alpha) e_0 \\
&+ (dF^\alpha_{a+m\,c} + F^\alpha_{a+m\,b}\theta^c_b + F^\alpha_{b+m\,a} \theta^c_{b+m} - 2F^\mu_{\alpha a} \theta^c_\mu \\
&\qquad + \delta_{ac} \omega^\alpha - L^c_{ab} \theta^{b+m}_\alpha -F^\beta_{a+m\,c} \theta^\beta_\alpha - F^\mu_{a+m\,c} \theta^\mu_\alpha) e_c \\
&+ (F^\alpha_{a+m\,b} \theta^{c+m}_b + dF^\alpha_{c+m\,a} + F^\alpha_{b+m\,a} \theta^{c+m}_{b+m} -2F^\mu_{\alpha a} \theta^{c+m}_\mu \\
&\qquad -\delta_{ac} \omega^\alpha + L^c_{ab} \theta^b_\alpha - F^\beta_{c+m\,a} \theta^\beta_\alpha + F^\mu_{c+m\,a} \theta^\mu_\alpha) e_{c+m} \\
&+(F^\alpha_{a+m\,b} \theta^\beta_b + F^\alpha_{b+m\,a} \theta^\beta_{b+m} - 2 F^\mu_{\alpha a} \theta^\beta_\mu - F^\beta_{a+m\,b} \theta^b_\alpha \\
&\qquad - F^\beta_{b+m\,a} \theta^{b+m}_\alpha + 2 F^\mu_{\beta a} \theta^\mu_\alpha)e_\beta \\
&+( F^\alpha_{a+m\,b} \theta^\mu_b + F^\alpha_{b+m\,a} \theta^\mu_{b+m} - 2 F^\nu_{\alpha a} \theta^\mu_\nu -2 dF^\mu_{\alpha a} \\
&\qquad - F^\mu_{a+m\,b} \theta^b_\alpha + F^\mu_{b+m\,a} \theta^{b+m}_\alpha + 2 F^\mu_{\beta a} \theta^\beta_\alpha ) e_\mu
\endaligned
\end{equation}
The coefficient of $\x$ is zero and the coefficient of $e_0$ is zero, both by~\eqref{eq:4:13}.  
For the coefficient of $e_c$, use~\eqref{eq:4:13}, \eqref{eq:5:2} and~\eqref{eq:8:4}, and add and subtract appropriate terms, to rewrite it as
\begin{equation*}
\aligned
&F^\alpha_{b+m\,c}(\theta^b_a - L^b_{ad} \omega^{d+m} + F^\beta_{b+m\,a} \omega^\beta + F^\mu_{b+m\,a} \omega^\mu) \\
&+(F^\alpha_{a+m\,cd} - F^\alpha_{b+m\,c} L^b_{ad} + F^\alpha_{b+m\,d} L^c_{ab} - F^\mu_{a+m\,c} F^\mu_{\alpha d}) \omega^d \\
&+( F^\alpha_{a+m\,c\,d+m} - 2F^\mu_{\alpha a} F^\mu_{d+m\,a} - F^\mu_{a+m\,c} F^\mu_{\alpha\,d+m}) \omega^{d+m} \\
&+(-F^\alpha_{b+m\,c} F^\beta_{b+m\,a} - F^\alpha_{b+m\,a} F^\beta_{b+m\,c} + F^\alpha_{a+m\,c\beta} -4 F^\mu_{\alpha a} F^\mu_{\beta c} + \delta_{ac} \delta_{\alpha\beta}) \omega^\beta \\
&+( -F^\alpha_{b+m\,c} F^\mu_{b+m\,a} + F^\alpha_{b+m\,a} F^\mu_{b+m\,c} + F^\alpha_{a+m\,c\mu} - 2L^c_{ab} F^\mu_{\alpha\,b+m}) \omega^\mu \\
&= F^\alpha_{b+m\,c} \nu^b_a
\endaligned
\end{equation*}
by~\eqref{eq:g:19}, because the other terms are zero as follows.  The coefficient of $\omega^d$ is zero by the first equation of~\eqref{eq:5:7} and~\eqref{eq:8:33}.  The coefficient of $\omega^{d+m}$ is zero by the second equation of~\eqref{eq:5:7}.  The coefficient of $\omega^\beta$ is zero by the third equation of~\eqref{eq:5:7} and then the second equation of~\eqref{eq:5:6}.  The coefficient of $\omega^\mu$ is zero by~\eqref{eq:5:11} and then~\eqref{eq:8:35}.

The coefficient of $e_{c+m}$ in $(dQ_a)e_\alpha$ in~\eqref{eq:g:dQaealpha} is, using~\eqref{eq:4:13}
\begin{equation*}
\aligned
&F^\alpha_{c+m\,b} ( \theta^b_a + L^b_{da} \omega^{d+m} + F^\beta_{b+m\,a} \omega^\beta + F^\mu_{b+m\,a}\omega^\mu) \\
&+( F^\alpha_{c+m\,ad} + 2F^\mu_{\alpha a} F^\mu_{c+m\,d} + F^\mu_{c+m\,a} F^\mu_{\alpha d}) \omega^d \\
&+( F^\alpha_{c+m\,a\,d+m} - L^c_{ab} F^\alpha_{d+m\,b} + F^\mu_{c+m\,a} F^\mu_{\alpha\,d+m} - L^b_{da} F^\alpha_{c+m\,b}) \omega^{d+m} \\
&+(F^\alpha_{c+m\,a\beta} + F^\alpha_{a+m\,b} F^\beta_{c+m\,b} + 4 F^\mu_{\alpha a} F^\mu_{\beta\,c+m} - \delta_{ac} \delta_{\alpha \beta} - F^\alpha_{c+m\,b} F^\beta_{b+m\,a}) \omega^\beta \\
&+(F^\alpha_{c+m\,a\mu} - F^\alpha_{a+m\,b} F^\mu_{c+m\,b} + 2 L^c_{ab} F^\mu_{\alpha b} - F^\alpha_{c+m\,b} F^\mu_{b+m\,a}) \omega^\mu \\
&= F^\alpha_{c+m\,b} \nu^b_a
\endaligned
\end{equation*}
by~\eqref{eq:g:19}, because the other terms are zero as follows.  The coefficient of $\omega^d$ is zero by the first equation in~\eqref{eq:5:7}.  The coefficient of $\omega^{d+m}$ is zero by the second equation in~\eqref{eq:5:7} and then~\eqref{eq:8:33}.  The coefficient of $\omega^\beta$ is zero by the third equation in~\eqref{eq:5:7} and then the second equation in~\eqref{eq:5:6}.  The coefficient of $\omega^\mu$ is zero by~\eqref{eq:5:11} and then~\eqref{eq:8:35}.

The coefficient of $e_\beta$ in $(dQ_a)e_\alpha$ in~\eqref{eq:g:dQaealpha} is, using~\eqref{eq:4:13}
\begin{equation*}
\aligned
&(F^\alpha_{b+m\, a} F^\beta_{b+m\,c} + 2F^\mu_{\alpha a} F^\mu_{\beta c} + 2F^\mu_{\beta a} F^\mu_{\alpha c} + F^\beta_{b+m\,a} F^\alpha_{b+m\,c}) \omega^c \\
&+( F^\alpha_{a+m\,b} F^\beta_{c+m\,b} + 2F^\mu_{\alpha a} F^\mu_{\beta\,c+m} + 2F^\mu_{\beta a} F^\mu_{\alpha\,c+m} + F^\beta_{a+m\,b} F^\alpha_{c+m\,b}) \omega^{c+m} \\
&-2(F^\alpha_{a+m\,b}F^\mu_{\beta b} + F^\alpha_{b+m\,a} F^\mu_{\beta\,b+m} + F^\beta_{a+m\,b} F^\mu_{\alpha b} + F^\beta_{b+m\,a} F^\mu_{\alpha\,b+m}) \omega^\mu \\
&= \delta_{\alpha \beta}(\omega^a + \omega^{a+m}) = \delta_{\alpha \beta} \nu^a
\endaligned
\end{equation*}
by~\eqref{eq:g:17}, because the coefficient of $\omega^c$ is $\delta_{\alpha\beta}\delta_{ac}$ by the second equation of~\eqref{eq:5:6}, and the coefficient of $\omega^{c+m}$ is also $\delta_{\alpha\beta}\delta_{ac}$ by~\eqref{eq:8:1}, \eqref{eq:8:2} and the second equation of~\eqref{eq:5:6}; and the coefficient of $\omega^\mu$ is zero by~\eqref{eq:8:1} and~\eqref{eq:8:2}.

The coefficient of $e_\mu$ in $(dQ_a)e_\alpha$ in~\eqref{eq:g:dQaealpha} is, using~\eqref{eq:4:13} and~\eqref{eq:5:2},
\begin{equation*}
\aligned
&-2F^\mu_{\alpha b} ( \theta^b_a + L^b_{ca} \omega^{c+m} + F^\beta_{b+m\,a} \omega^\beta + F^\nu_{b+m\,a} \omega^\nu) \\
&+( -2F^\mu_{\alpha a c} + F^\alpha_{b+m\,a} F^\mu_{b+m\,c} - F^\mu_{b+m\,a} F^\alpha_{b+m\,c})\omega^c \\
&+ ( -2F^\mu_{\alpha a\,c+m} + 2F^\mu_{\alpha b} L^b_{ca} - F^\alpha_{a+m\,b} F^\mu_{c+m\,b} + F^\mu_{a+m\,b} F^\alpha_{c+m\,b}) \omega^{c+m} \\
&+2(-F^\mu_{\alpha a \beta} + F^\beta_{b+m\,a} F^\mu_{\alpha b} - F^\alpha_{a+m\,b} F^\mu_{\beta b} + F^\alpha_{b+m\,a} F^\mu_{\beta\,b+m}) \omega^\beta \\
&+2( -F^\mu_{\alpha a \nu} + F^\mu_{\alpha b} F^\nu_{b+m\,a} - F^\mu_{a+m\,b} F^\nu_{\alpha b} + F^\mu_{b+m\,a} F^\nu_{\alpha\,b+m}) \omega^\nu \\
&= -2F^\mu_{\alpha b} \nu^b_a
\endaligned
\end{equation*}
by~\eqref{eq:g:17}, because the coefficient of $\omega^c$ is zero by the first equation in~\eqref{eq:5:9}, the coefficient of $\omega^{c+m}$ is zero by~\eqref{eq:5:11} and then~\eqref{eq:8:35}, the coefficient of $\omega^\beta$ is zero by~\eqref{eq:8:1}, \eqref{eq:8:2} and the second equation in~\eqref{eq:5:9}; and the coefficient of $\omega^\nu$ is zero by~\eqref{eq:8:1}, \eqref{eq:8:3} and the third equation in~\eqref{eq:5:9}.

This completes the verification of~\eqref{eq:g:55}, which verifies that~\eqref{eq:g:21} holds when both sides are applied to $e_\alpha$.

The final case is to verify~\eqref{eq:g:21} when both sides are
applied to $e_\mu$.  We must verify that 
\begin{equation}\label{eq:g:56}
\aligned
&d(Q_ae_\mu) - Q_a(de_\mu) = (dQ_a)e_\mu =\nu^0_a Q_0 e_\mu + \nu^b_a Q_b e_\mu \\
&= -\nu^a e_\mu + F^\mu_{b+m\,c} \nu^b_a e_c - F^\mu_{c+m\,b} \nu^b_a e_{c+m} - 2 F^\mu_{\alpha b} \nu^b_a e_\alpha
\endaligned
\end{equation}
Using~\eqref{eq:g:4} and~\eqref{eq:4:15}, and gathering together the coefficients of each basis vector, we get for the left hand side
\begin{equation}\label{eq:g:dQaemu}
\aligned
&(dQ_a)e_\mu = (F^\mu_{b+m\,a} \omega^{b+m} + 2F^\mu_{\alpha a} \omega^\alpha - \theta^a_\mu)\x \\
&- (F^\mu_{a+m\,b} \omega^b + 2F^\mu_{\alpha a}\omega^\alpha + \theta^{a+m}_\mu)e_0 \\
&+ (dF^\mu_{a+m\,c} + F^\mu_{a+m\,b} \theta^c_b - F^\mu_{b+m\,a} \theta^c_{b+m} - 2F^\mu_{\alpha a} \theta^c_\alpha + \delta_{ac} \omega^\mu \\
&\qquad -L^c_{ab} \theta^{\b+m}_\mu - F^\alpha_{a+m\,c} \theta^\alpha_\mu - F^\nu_{a+m\,c} \theta^\nu_\mu) e_c \\
&+(F^\mu_{a+m\,b}\theta^{c+m}_b - dF^\mu_{c+m\,a} - F^\mu_{b+m\,a} \theta^{c+m}_{b+m} -2 F^\mu_{\alpha a} \theta^{c+m}_\alpha + \delta_{ac}\omega^\mu \\
&\qquad +L^c_{ab}\theta^b_\mu - F^\alpha_{c+m\,a} \theta^\alpha_\mu + F^\nu_{c+m\,a} \theta^\nu_\mu) e_{c+m} \\
&+( F^\mu_{a+m\,b} \theta^\alpha_b - F^\mu_{b+m\,a} \theta^\alpha_{b+m} -2 dF^\mu_{\alpha a} -2F^\mu_{\beta a} \theta^\alpha_\beta - F^\alpha_{a+m\,b} \theta^b_\mu \\
&\qquad -F^\alpha_{b+m\,a} \theta^{b+m}_\mu + 2F^\nu_{\alpha a} \theta^\nu_\mu)e_\alpha \\
&+(F^\mu_{a+m\,b} \theta^\nu_b - F^\mu_{b+m\,a} \theta^\nu_{b+m} -2 F^\mu_{\alpha a} \theta^\nu_\alpha - F^\nu_{a+m\,b} \theta^b_\mu + F^\nu_{b+m\,a} \theta^{b+m}_\mu + 2F^\nu_{\alpha a} \theta^\alpha_\mu)e_\nu
\endaligned
\end{equation}
The coefficient of $\x$ is zero by~\eqref{eq:4:13}.  The coefficient of $e_0$ is zero by~\eqref{eq:4:13} and~\eqref{eq:8:1}.

After applying~\eqref{eq:5:2} and~\eqref{eq:4:13} and adding and then subtracting some terms in the definition of $\nu^b_a$ in~\eqref{eq:g:19}, we can rewrite the coefficient of $e_c$ as
\begin{equation*}
\aligned
&F^\mu_{b+m\,c}(\theta^b_a + L^b_{da} \omega^{d+m} + F^\beta_{b+m\,a} \omega^\beta + F^\nu_{b+m\,a} \omega^\nu) \\
&+ (F^\mu_{a+m\,cd} - F^\mu_{b+m\,c} L^b_{ad} + L^c_{ab} F^\mu_{b+m\,d} + F^\alpha_{a+m\,c} F^\mu_{\alpha d}) \omega^d \\
&+( F^\mu_{a+m\,c\,d+m} + 2 F^\mu_{\alpha a} F^\alpha_{d+m\,c} +  F^\alpha_{a+m\,c} F^\mu_{\alpha\,d+m}) \omega^{d+m} \\
&+( F^\mu_{a+m\,c\alpha} + F^\mu_{b+m\,a} F^\alpha_{b+m\,c} +2L^c_{ab} F^\mu_{\alpha\,b+m} - F^\mu_{b+m\,c} F^\alpha_{b+m\,a})\omega^\alpha \\
&+(F^\mu_{a+m\,c\nu} - F^\mu_{b+m\,a} F^\nu_{b+m\,c} -4F^\mu_{\alpha a} F^\nu_{\alpha c} + \delta_{ac} \delta_{\mu\nu} - F^\mu_{b+m\,c} F^\nu_{b+m\,a}) \omega^\nu \\
&= F^\mu_{b+m\,c} \nu^b_a
\endaligned
\end{equation*}
by~\eqref{eq:g:19} and the following.  The coefficient of $\omega^d$ is zero by the first equation in~\eqref{eq:5:8} and then~\eqref{eq:8:34}.  The coefficient of $\omega^{d+m}$ is zero by~\eqref{eq:8:1} and~\eqref{eq:5:8}.  The coefficient of $\omega^\alpha$ is zero by~\eqref{eq:5:11} and then~\eqref{eq:8:35}.  The coefficient of $\omega^\mu$ is zero by the third equation in~\eqref{eq:5:8} and then~\eqref{eq:8:1} and the fourth equation in~\eqref{eq:5:6}.

Using~\eqref{eq:5:2} and~\eqref{eq:4:13}, we can rewrite the coefficient of $e_{c+m}$ in $(dQ_a)e_\mu$ in~\eqref{eq:g:dQaemu} as
\begin{equation*}
\aligned
&-F^\mu_{c+m\,b} (\theta^b_a + L^b_{da} \omega^{d+m} + F^\alpha_{b+m\a} \omega^\alpha + F^\nu_{b+m\,a} \omega^\nu) \\
&+(-F^\mu_{c+m\,ab} + 2F^\mu_{\alpha a} F^\alpha_{c+m\,b} + F^\alpha_{c+m\,a} F^\mu_{\alpha b})\omega^b \\
&+(L^b_{da} F^\mu_{c+m\,b} - F^\mu_{c+m\,a\,d+m} + L^c_{ab} F^\mu_{d+m\,b} + F^\alpha_{c+m\,a} F^\mu_{\alpha\,d+m}) \omega^{d+m} \\
&+( F^\mu_{c+m\,b} F^\alpha_{b+m\,a} - F^\mu_{c+m\,a\alpha} + F^\mu_{a+m\,b} F^\alpha_{c+m\,b} + 2F^\mu_{\alpha b} L^c_{ab}) \omega^\alpha \\
&+( F^\mu_{c+m\,b} F^\nu_{b+m\,a} - F^\mu_{c+m\,a\nu} - F^\mu_{a+m\,b} F^\nu_{c+m\,b} -4 F^\mu_{\alpha a} F^\nu_{\alpha \,c+m} + \delta_{ac} \delta_{\mu\nu}) \omega^\nu\\
&= - F^\mu_{c+m\,b} \nu^b_a
\endaligned
\end{equation*}
by~\eqref{eq:g:19} and the following.  The coefficient of $\omega^b$ is zero by the first equation in~\eqref{eq:5:8}.  The coefficient of $\omega^{d+m}$ is zero by the second equation in~\eqref{eq:5:8} and then~\eqref{eq:8:1} and~\eqref{eq:8:34}.  The coefficient of $\omega^\alpha$ is zero by~\eqref{eq:5:11}, then~\eqref{eq:8:2} and~\eqref{eq:8:3} and~\eqref{eq:8:35}.  The coefficient of $\omega^\nu$ is zero by the third equation in~\eqref{eq:5:8}, then~\eqref{eq:8:1} and~\eqref{eq:8:3} and then the fourth equation in~\eqref{eq:5:6}.

Using~\eqref{eq:5:2} and~\eqref{eq:4:13}, we can rewrite the coefficient of $e_\alpha$ in $(dQ_a)e_\mu$ in~\eqref{eq:g:dQaemu} as
\begin{equation*}
\aligned
&-2F^\mu_{\alpha b}( \theta^b_a + L^b_{ca} \omega^{c+m} + F^\beta_{b+m\,a} \omega^\beta + F^\nu_{b+m\,a} \omega^\nu) \\
&+(-2F^\mu_{\alpha a c} - F^\mu_{b+m\,a} F^\alpha_{b+m\,c} + F^\alpha_{b+m\,a} F^\mu_{b+m\,c}) \omega^c \\
&+(-2F^\mu_{\alpha a\,c+m} + 2F^\mu_{\alpha b} L^b_{ca} + F^\mu_{a+m\,b} F^\alpha_{c+m\,b} - F^\mu_{c+m\,b} F^\alpha_{a+m\,b}) \omega^{c+m} \\
&+2( -F^\mu_{\alpha a \beta} + F^\mu_{\alpha b} F^\beta_{b+m\,a} - F^\mu_{\beta b} F^\alpha_{a+m\,b} + F^\mu_{\beta\, b+m} F^\alpha_{b+m\,a}) \omega^\beta \\
&+2( -F^\mu_{\alpha a \nu} + F^\mu_{\alpha b} F^\nu_{b+m\,a} - F^\mu_{a+m\,b} F^\nu_{\alpha b} + F^\mu_{b+m\,a} F^\nu_{\alpha\,b+m}) \omega^\nu \\
&= -2F^\mu_{\alpha b} \nu^b_a
\endaligned
\end{equation*}
by~\eqref{eq:g:19} and the following.  The coefficient of $\omega^c$ is zero by the first equation in~\eqref{eq:5:9}.  The coefficient of $\omega^{c+m}$ is zero by~\eqref{eq:5:11}, then~\eqref{eq:8:2} and~\eqref{eq:8:3} and then~\eqref{eq:8:35}.  The coefficient of $\omega^\beta$ is zero by the second equation in~\eqref{eq:5:9} and then~\eqref{eq:8:1}-\eqref{eq:8:3}.  The coefficient of $\omega^\nu$ is zero by the third equation in~\eqref{eq:5:9}, then~\eqref{eq:8:1} and~\eqref{eq:8:3}.

Using~\eqref{eq:4:13}, we can rewrite the coefficient of $e_\mu$ in $(dQ_a)e_\mu$ in~\eqref{eq:g:dQaemu} as
\begin{equation*}
\aligned
&(-F^\mu_{b+m\,a} F^\nu_{b+m\,c} -2F^\mu_{\alpha a} F^\nu_{\alpha c} - F^\nu_{b+m\,a} F^\mu_{b+m\,c} -2F^\nu_{\alpha a} F^\mu_{\alpha c}) \omega^c \\
&+( -F^\mu_{a+m\,b} F^\nu_{c+m\,b} -2F^\mu_{\alpha a} F^\nu_{\alpha\,c+m} - F^\nu_{a+m\,b} F^\mu_{c+m\,b} - 2F^\nu_{\alpha a} F^\mu_{\alpha\,c+m}) \omega^{c+m} \\
&+2(-F^\mu_{a+m\,b} F^\nu_{\alpha b} -F^\mu_{b+m\,a} F^\nu_{\alpha\,b+m} - F^\nu_{a+m\,b} F^\mu_{\alpha b} - F^\nu_{b+m\,a} F^\mu_{\alpha\,b+m}) \omega^\alpha \\
&= -\delta_{\mu\nu}(\omega^a + \omega^{a+m}) = - \delta_{\mu\nu} \nu^a
\endaligned
\end{equation*}
by~\eqref{eq:8:1}, \eqref{eq:8:3} and the fourth equation in~\eqref{eq:5:6}.  This completes the verification of~\eqref{eq:g:21} when both sides are applied to $e_\mu$, and therefore also completes the verification of~\eqref{eq:g:21}.
\end{proof}


\section{The quadratic forms}
For the remainder of the paper, we will again refer to the two multiplicities
as $m_{1}$ and $m_{2}$, rather than $m$ and $N$, respectively, and we will no longer use the Einstein summation convention.
Our task now is to solve~\eqref{eq:8:1} through~\eqref{eq:8:4}. It is
known that
$m_{1}=m_{2}$ only when $m_{1}=m_{2}=1$, which is of FKM-type, or
$m_{1}=m_{2}=2$, which is not of FKM-type~\cite{Ab}.
Therefore we assume
$m_{1}\neq m_{2}$ henceforth. Our convention is that $m_{1}<m_{2}$ and we
denote by
$M_{+}$ (respectively, $M_{-}$) the focal submanifold whose co-dimension is
$m_{1}+1$ (respectively, $m_{2}+1$) in the ambient sphere. We change the
Cartan-M\"{u}nzner polynomial $F$ to $-F$ if necessary so that always
$M_{+}=f^{-1}(1)$ with respect to the isoparametric function $f$. Now that
$m_{1}<m_{2}$ we only have to show the validity
of~\eqref{eq:8:1} and the spanning property in view of Theorem~\ref{th:g:1}.

As in Section~\ref{section4} let $\x\in M_{+}$ and let $e_{0}$ be a unit normal vector to
$M_{+}$ at $\x$ for which the shape operator $S_{e_{0}}$ assumes the
eigenspaces $V_{0},V_{+}$ and $V_{-}$ with eigenvalues $0,1,$ and $-1$,
respectively. For an orthonormal basis $e_{0},\dots,e_{m_{1}}$ of the
normal space to $M_{+}$ at $\x$ we introduce the quadratic homogeneous
polynomials
$$
\check{p}_{i}(\mathbf z):=S_{e_{i}}\mathbf z \cdot \mathbf z
$$
for $0\leq i\leq m_{1}$, where $\mathbf z$ is tangent to $M_+$ at $\x$. Consider the set
$$
\check{\mathcal D}:=\{\mathbf z\in V_{+}\oplus V_{-}:|\mathbf z|=1,\check{p}_{i}(\mathbf z)=0,
0\leq i\leq m_{1}\}.
$$

\begin{proposition}\label{pro9.1}
$\check{\mathcal D}=(V_{+}\oplus V_{-})\cap M_{+}$.
\end{proposition}

\begin{proof}
This follows from the formula of~\cite[I, pp.524-526]{OT},
that reads
\begin{eqnarray}\label{eq9.1}
\begin{aligned}
F(t\x+\y+\mathbf w)
&=t^{4}+(2|\y|^2-6|\mathbf w|^2)t^2
+8(\sum_{i=0}^{m_{1}}\check{p}_{i}w_{i})t\\
&+|\y|^4-2\sum_{i=0}^{m_{1}}(\check{p}_{i})^{2}
+8\sum_{i=0}^{m_{1}}q_{i}w_{i}\\
&+2\sum_{i,j=0}^{m_{1}} (\nabla \check{p}_{i} \cdot \nabla \check{p}_{j}) w_{i}w_{j}
-6|\y|^2|\mathbf w|^2+|\mathbf w|^4.
\end{aligned}
\end{eqnarray}
Here, $q_{i}$ are the components of the third fundamental form of $M_{+}$,
$\mathbf w=\sum_{i=0}^{m_{1}}w_{i}e_{i}$  and $\y$ is tangent to
$M_{+}$. For the 
convenience of the reader, let us briefly recall that Ozeki and Takeuchi
expanded $F(t\x+\y+\mathbf w)$ in terms of $t$ and substituted it into its 
governing partial differential equations mentioned in Section~\ref{section2} to get
$$
F(t\x+\y+\mathbf w)=t^{4}+At^{2}+Bt+C,
$$
where $A$ is derived on p525, $B$ is on p526, and 
$C=C_{0}+\dots+C_{4}$, in which $C_{s}$, given on p526, is
the homogeneous part of $C$ of degree $s$
in the normal coordinates $w_{0},\dots,w_{m_{1}}$.
When one sets $t=0$, $\mathbf w = 0$ and $\y\in V_{+}\oplus V_{-}$ in
the formula~\eqref{eq9.1} one gets 
$$
F(\y)-{|\y|}^{4}=-2\sum_{i=0}^{m_{1}}(\check{p}_{i}(\y))^{2}.
$$
Hence when $|\y|=1$, we have $F(\y)=1$ if and only if $\check{p}_{i}(\y)=0$ for
$0\leq i\leq m_{1}$.
\end{proof}

In view of Proposition~\ref{pro9.1} we set $p_{i}$ to be the restriction of
$\frac14 \check{p}_{i}$ to the space
$V_{+}\oplus V_{-}$ for $1\leq i\leq m_{1}$, and set $p_0$ to be the
restriction of $\check p_0$ to this space. These are the quadratic
polynomials $p_0, p_a$ defined in~\eqref{eq:6:polys}.
Recall from~\eqref{eq:4:27} and~\eqref{eq:6:polys}, that relative to a
second order Darboux frame we have variables $x = (x_\alpha)$ and $y =
(y_\mu)$ in terms of which these polynomials are 
\begin{eqnarray}
p_{0}(x,y)&=&\sum_{\alpha=1}^{m_{2}}(x_{\alpha})^2
-\sum_{\mu=1}^{m_{2}}(y_{\mu})^{2} \label{eq9.2}\\
p_{a}(x,y)&=&\sum_{\alpha,\mu=1}^{m_{2}}F^{\mu}_{\alpha a}x_{\alpha}y_{\mu} \label{eq9.3}
\end{eqnarray}
For notational ease without any possibility of
confusion, we will stick to the range
$1\leq \alpha,\mu\leq m_{2}$ for $x_{\alpha}$ and $y_{\mu}$ from now on even
though $\alpha$ and $\mu$ live in the designated
ranges as given in~\eqref{eq:4:convention}.

As mentioned in Section~\ref{section4}, we know $e_{0}$ also lies in $M_{+}$
with the normal space span$(\x,e_{m_{1}+1},\dots,e_{2m_{1}})$.
The $0,+1,-1$ eigenspaces of the shape operator $S_{\x}$ at $e_{0}$
are, respectively, span$(e_{1},\dots,e_{m_{1}})$, $V_{+}$ and $V_{-}$. 
With respect to the normal
basis $\x,e_{p},m_{1}+1\leq p\leq 2m_{1},$ at $e_{0}$, we let
$\overline{p}_{0},\dots,\overline{p}_{m_{1}}$ be the counterparts of
$p_{0},\dots,p_{m_{1}}$, respectively, as in~\eqref{eq:6:13}. Then Proposition~\ref{pro9.1}
immediately gives the following simple but crucial observation.

\begin{proposition}\label{pro9.2}
$\check{D}:=\{\mathbf z\in V_{+}\oplus V_{-}:|\mathbf z|=1, \overline{p}_{i}(\mathbf z)=0,0\leq i\leq m_{1}\}.$
\end{proposition}
Now $\check{\mathcal D}$ can be viewed from a different angle.
Observe that all
$\z=(x_{1},\dots,x_{m_{2}},y_{1},\dots,y_{m_{2}})\in \check{D}$ must
satisfy $\sum_{\alpha=1}^{m_{2}}(x_{\alpha})^2
+\sum_{\mu=1}^{m_{2}}(y_{\mu})^2=1$. It follows that
$\z \in S^{m_{2}-1}\times S^{m_{2}-1}$ due to the fact that $p_{0}(\z)=0$,
where $S^{m_{2}-1}$ is the standard sphere of radius $1/\sqrt{2}$. The real
projective variety out of $S^{m_{2}-1}\times S^{m_{2}-1}$ is
${\mathbf R}P^{m_{2}-1}\times {\mathbf R}P^{m_{2}-1}$. Note that the solution
to $p_{a}=0$, $1\leq a\leq m_{1}$, lives naturally in
${\mathbf R}P^{m_{2}-1}\times {\mathbf R}P^{m_{2}-1}$, which is parametrized
by $[x_{1}:\dots:x_{m_{2}}]\times [y_{1}:\dots:y_{m_{2}}]$. As a consequence
the projectivized $\check{D}$ in
${\mathbf R}P^{m_{2}-1}\times {\mathbf R}P^{m_{2}-1}$ via the map
$S^{m_{2}-1}\times S^{m_{2}-1}\longrightarrow
{\mathbf R}P^{m_{2}-1}\times {\mathbf R}P^{m_{2}-1}$ is exactly
\begin{equation}\label{eq9.4}
\aligned
{\mathcal D}:=\{[\z]\in {\mathbf R}P^{m_{2}-1}\times {\mathbf R}P^{m_{2}-1}
:p_{a}(\z)=0,1\leq a\leq m_{1}\}.
\endaligned
\end{equation}
Since the $+1$ and $-1$ eigenspaces of the shape operator $S_\x$ at $e_0$ are $V_+$ and $V_-$, respectively, it follows from~\eqref{eq9.2} that 
$\overline{p}_{0}=p_{0}$. Hence,
Proposition~\ref{pro9.2} can be rephrased as follows.

\begin{proposition}\label{pro9.3}
The zero locus of $p_{1},\dots,p_{m_{1}}$ in
${\mathbf R}P^{m_{2}-1}\times {\mathbf R}P^{m_{2}-1}$ is identical with that
of $\overline{p}_{1},\dots,\overline{p}_{m_{1}}$.
\end{proposition}

\begin{lemma}\label{le9.1}
If $m_2 \geq m_1+2$, then the quadratic forms $p_1,\dots, p_{m_1}$ are linearly independent and irreducible, both over the real numbers $\R$ and over the complex numbers $\C$.
\end{lemma}

\begin{proof}
The quadratic form $p_a(x,y)$ is given by
\begin{equation*}
4p_a(x,y) = \begin{pmatrix} 0 & A_a \\ \tose A_a & 0 \end{pmatrix}
\begin{pmatrix} x \\ y \end{pmatrix} \cdot \begin{pmatrix} x \\ y
\end{pmatrix} = 2A_a y \cdot x 
\end{equation*}
where $A_a$ is the matrix with respect to $e_\alpha, e_\mu$ of the
operator defined in the first equation of~\eqref{eq:6:5}. 
Recall that the rank of $p_a$ is defined to be the rank of the matrix
of the associated bilinear form, 
\begin{equation*}
\mbox{rank\,}(p_a) = \mbox{rank}\begin{pmatrix} 0 & A_a \\ \tose A_a &
0 \end{pmatrix} = 2\,\mbox{rank\,}A_a 
\end{equation*}
Let $S_a = U+V$ be the shape operator given in~\eqref{eq:6:7}, where 
\begin{equation*}
U = \begin{pmatrix} 0 & A_a & 0 \\ \tose A_a & 0 & 0 \\ 0&0&0
\end{pmatrix}, \qquad V = \begin{pmatrix} 0&0& B_a \\ 0&0& C_a \\ \tose
B_a & \tose C_a &0 \end{pmatrix}
\end{equation*}
Then $\text{rank\,}S_a \leq \text{rank\,}U + \text{rank\,}V$.  Since
$\text{rank\,}S_a = 2m_2$, $\text{rank\,}U = 2\text{\,rank\,}A_a$ and
$\text{rank\,}V \leq 2m_1$, we get
\begin{equation*}
2(m_2-m_1) \leq 2\text{\,rank\,}A_a = \text{\,rank\,}(p_a)
\end{equation*}
for all $a$, as proved by Ozeki and Takeuchi~\cite[II, p45]{OT}.
If $p_a$ is reducible, then $p_a = f g$ is a product of linear forms
$f = a_\alpha x_\alpha + a_\mu y_\mu$ and $g = b_\alpha x_\alpha +
b_\mu y_\mu$.  If we let 
$\mathbf a = \tose(a_\alpha\, a_\mu)$ and $\mathbf b = \tose(b_\alpha\, b_\mu) \in \R^{2m_2}$
then the symmetric matrix of the bilinear form must be
\begin{equation*}
\begin{pmatrix} 0 & A_a \\ \tose A_a & 0 \end{pmatrix} =
\frac12(\mathbf a \tose\mathbf b + \mathbf b \tose \mathbf a) 
\end{equation*}
which has rank $\leq 2$, as each column is a linear combination of $\mathbf a$ and $\mathbf b$.  
In particular, if $m_2-m_1 \geq 2$, then $\mbox{rank}\, (p_a) \geq 4 >2$ and hence, $p_a$ is irreducible over $\R$.
Notice that this discussion is unchanged if we work over the complex numbers, which shows that they are irreducible over $\C$ as well.
Linear independence of $p_1,\dots,p_{m_1}$ over $\R$ is equivalent to linear independence of $A_1,\dots,A_{m_1}$, which follows under our hypotheses from Proposition~\ref{pr:8:4}.  Being real polynomials, they are also linearly independent over $\C$.
\end{proof}

\section{Some commutative algebra and algebraic geometry}
We will explore in more depth the fact that $p_{a},1\leq a\leq m_{1}$, are
irreducible when $m_{2}\geq m_{1}+2$ and are bi\-homo\-geneous, i.e., are
homogeneous in $x_{1},\dots,
x_{m_{2}}$ and in $y_{1},\dots,y_{m_{2}}$, of bi-degree $(1,1)$ in this
section. We shall pursue commutative algebra only to the extent that serves
our need, and shall stress the geometry behind the algebra. A few ad hoc
proofs and examples will be given to convey to the reader, who might be
unfamiliar with the
subject, some intuition about the concepts encountered.
Henceforth, $n$ is just an index that has nothing to do with the dimension
of the ambient sphere in which the isoparametric hypersurface sits.

\begin{definition}\label{reducedness}
 Let ${\mathbf F}$ be either ${\mathbf R}$ or ${\mathbf C}$
and let ${\mathbf F}[x_{1},\dots,x_{s},y_{1},\dots,y_{s}]$ be the polynomial
ring in variables $x_{1},\dots,x_{s},y_{1},\dots,y_{s}$ over ${\mathbf F}$. 
Given bi\-homo\-geneous polynomials $p_{1},\dots,p_{n}$, we say that the ideal
$I:=(p_{1},\dots,p_{n})$ in ${\mathbf F}[x_{1},\dots,x_{s},y_{1},\dots,
y_{s}]$ is {\em reduced} if

(i). The bi-projective variety 
$${\mathbf P}_{b}V_{I}:=
\{([x], [y]) \in {\mathbf F}P^{s-1}\times {\mathbf F}P^{s-1}:p_{a}(x,y)=0,1
\leq a\leq n\}
$$
is not empty, and

(ii). whenever
$f\in {\mathbf F}[x_{1},\dots,x_{s},y_{1},\dots,y_{s}]$ satisfies
$f|_{{\mathbf P}_{b}V_{I}}\equiv 0$ then we have
$$
f=p_{1}f_{1}+\dots+p_{n}f_{n}
$$
for some $f_{1},\dots,f_{n}\in {\mathbf F}[x_{1},\dots,x_{s},y_{1},
\dots,y_{s}]$.

We call the affine variety
$V_{I}:=\{(x,y)\in {\mathbf C}^{s}\times {\mathbf C}^{s}: p_{a}(x,y)=0,
1\leq a\leq n\}$ a \textit{bi-affine cone}.
\end{definition}

For instance, when ${\mathbf F}={\mathbf C}$, the radical of $I$, denoted by
${\rm rad}(I)$, is always reduced. This is Hilbert's Nullstellensatz
indeed~\cite{Ei}.
In particular, since a prime ideal equals its radical, the ideal
$I$ will be
reduced if $I$ is a prime ideal. ${\mathbf P}_{b}V_{I}$ is not empty
automatically in this case, because otherwise
$V_{I}=({\mathbf C}^{s}\times\{0\})\cup(\{0\}\times {\mathbf C}^{s})$ would
not
be irreducible. We will extensively probe the primeness of $I$ subsequently.
(See~\cite{Fu} and~\cite{Mum} for bi-projective geometry.)
Before we proceed, let us introduce a notation. When $p$ is a real
polynomial, we denote by $p^{\mathbf C}$ the same polynomial whose variables
are over the complex numbers. We call $p^{\mathbf C}$ the \textit{complexification} of
$p$. Likewise, when $p_{1},\cdots,p_{n}$ are bi-homogeneous in 
${\mathbf R}[x_{1},\dots,x_{s},y_{1},\dots,y_{s}]$, we denote by $V$ the
resulting real bi-affine cone
and by $V^{\mathbf C}$ the complex bi-affine cone defined by the
complexifications of $p_{1},\cdots,p_{n}$.

\begin{lemma}\label{le10.1}
Suppose $V$ is a bi-affine cone in
${\mathbf R}^{s}\times {\mathbf R}^{s}$ defined by the real polynomials
$p_{1},\cdots,p_{n}$, such that its complex counterpart
$V^{\mathbf C}$ is irreducible and such that
$\dim_{\mathbf R}(V)=\dim_{\mathbf C}(V^{\mathbf C})$. If a real polynomial
$p(x_{1},\dots,x_{s},y_{1},\dots,y_{s})$ satisfies  
$p|_{V}\equiv 0$, then $p^{\mathbf C}|_{V^{\mathbf C}}\equiv 0$.
Here, by the dimension of\; $V$ we mean the maximal dimension of all the
irreducible components of\; $V$.
\end{lemma}

\begin{proof}
Suppose $p^{\mathbf C}|_{V^{\mathbf C}}$ is not identically zero on
$V^{\mathbf C}$. Then $p^{\mathbf C}$ cuts out a subvariety $X$, all of
whose irreducible components are of co-dimension
1 in $V^{\mathbf C}$ \cite[p59]{Sh1}. Clearly, $V\subset X$. Then we have
$$
\dim_{\mathbf R}(V)\leq\dim_{\mathbf C}(X)=\dim_{\mathbf C}(V^{\mathbf C})-1,
$$
in contradiction to the assumption that
$\dim_{\mathbf R}(V)=\dim_{\mathbf C}(V^{\mathbf C})$.
The inequality holds true because any real analytic parametrization
$\sigma:t=(t_{1},\cdots,t_{k})\longmapsto
(x_{1},\cdots,x_{s},y_{1},\cdots,y_{s})\in V$ around a smooth point, at
$t=0$, of $V$, satisfies
$p_{1}(\sigma(t))=\cdots=p_{n}(\sigma(t))=p(\sigma(t))=0$. 
The convergent power series defining $\sigma$ remain so when 
$t_{1},\cdots,t_{k}$ are allowed to be
complex variables, and then $\sigma(t)$ is a holomorphic map,
nonsingular at $t=0$, such that 
$p_{1}^{\mathbf C}(\sigma(t))=\cdots=p_{n}^{\mathbf C}(\sigma(t))=
p^{\mathbf C}(\sigma(t))=0$ because
a holomorphic function vanishing on the real part is identically
zero. That is, $\sigma(t)$, with $t$ complex, is a 
holomorphic map, nonsingular at $t=0$, into $X$. Therefore, we conclude that 
$\dim_{\mathbf C}(X)\geq \dim_{\mathbf R}(V)$.
\end{proof}

\begin{proposition}\label{pro10.1}
If $p_{1},\dots,p_{n}\in {\mathbf R}[x_{1},\dots,x_{s},y_{1},\dots,
y_{s}]$ are bi\-homo\-geneous polynomials of positive degree in each set of variables, and if $p_{1}^{\mathbf C},
\dots,p_{n}^{\mathbf C}$, their complexifications, are such that
\begin{description}
\item[(1)] $V^{\mathbf C}:=\{z\in {\mathbf C}^{s}\times
{\mathbf C}^{s}:p_{a}^{\mathbf C}(z)=0,1\leq a\leq n\}$ is irreducible,
\item[(2)] ${\rm rad}(I)=I$, where $I:=(p_{1}^{\mathbf C},\dots,p_{n}^{\mathbf
C})$, and
\item[(3)] $\dim_{\mathbf R}(V)=\dim_{\mathbf C}(V^{\mathbf C})$, where
$V:=\{z\in {\mathbf R}^{s}\times
{\mathbf R}^{s}:p_{a}(z)=0,1\leq a\leq n\}$.
\end{description}
then the real ideal $(p_{1},\dots,p_{n})$ is reduced.
\end{proposition}

\begin{proof} $I$ is a prime ideal by the first two assumptions.
Therefore, the remark immediately after Definition~\ref{reducedness} ensures
that ${\mathbf P}_{b}V^{\mathbf C}$ is not empty. Moreover,
$\dim_{\mathbf C}(V^{\mathbf C})>s$ by the first assumption and the fact that the
reducible
$({\mathbf C}^{s}\times \{0\})\cup (\{0\}\times {\mathbf C}^{s})$ is
contained in $V^{\mathbf C}$. Hence 
${\mathbf P}_{b}V$ is not empty either by the third assumption. So the first
condition in
Definition~\ref{reducedness} holds. Let $f$ be a real polynomial vanishing
on ${\mathbf P}_{b}V$ so that $f$ vanishes on $V$ as well;
by Lemma~\ref{le10.1} its complexification $f^{\mathbf C}$ vanishes on
$V^{\mathbf C}$.
It follows from the reducedness of $I$ that there are complex bi-homogeneous
polynomials $h_{1},\dots,h_{n}$ such that
$$
f^{\mathbf C}=p_{1}^{\mathbf C}h_{1}+\dots+p_{n}^{\mathbf C}h_{n}.
$$
Let $f_{1},\dots,f_{n}$ be, respectively, the real parts of $h_{1},
\dots,h_{n}$ when they are restricted to the real variables. We have, by
the realness of $f$ and $p_{1},\dots,p_{n}$, that
$$
f=p_{1}f_{1}+\dots+p_{n}f_{n}.
$$                                      
\end{proof}

We now review some important notions and properties from commutative algebra,
leaving detailed expositions to~\cite{Ei} and~\cite{Ma}.

\begin{definition}
Let $R$ be a commutative ring with identity. We say that $n$ elements
$x_{1},\dots,x_{n}\in R$ form a \textit{regular sequence} if $(x_{1},\dots,x_{n})
\neq R$, $x_{1}$ is not a zero
divisor in $R$ and $x_{i+1}$ is not a zero divisor in the quotient ring
$R/I_{i}$, where $I_{i}$ is the ideal $(x_{1},\dots,x_{i})$, for
$1\leq i\leq n-1$.
\end{definition}

\begin{example} A single nonconstant
$p\in {\mathbf C}[z_{1},\dots,z_{L}]$
clearly forms a regular sequence.
\end{example}

\begin{example}\label{example2} Let $p_{1}$ and $p_{2}$ in 
${\mathbf C}[z_{1},\dots,z_{L}]$ be relatively prime homogeneous polynomials
of degree $\geq 1$. Then $p_{1}$ and
$p_{2}$ form a regular sequence. This follows simply from the fact that
$p_{2}f=p_{1}g$ implies $f=p_1h$ for some $h$. Moreover, $(p_{1},p_{2})$ is
not the entire polynomial ring due to the homogeneity of $p_{1}$ and $p_{2}$.
\end{example}

\begin{definition} Let ${\mathcal P}$ be a prime ideal in a commutative ring
$R$ with identity. We define the {\em co-dimension} of ${\mathcal P}$ to be
$$
{\rm codim}({\mathcal P})
=\sup\{s:{\rm there\; is\; a\; prime\; chain}\; {\mathcal P}_{s}\subset
\dots\subset {\mathcal P}_{1}\subset {\mathcal P}_{0}={\mathcal P}\},
$$
where the set inclusions are all proper. For an arbitrary ideal $I$ we define
$$
{\rm codim}(I)=\inf_{I\subset {\mathcal P}}\{{\rm codim}({\mathcal P})\},
$$
and define the {\em depth} of $I$ to be
$$
{\rm depth}(I)=\sup\{n:{\rm there\; is\; a\; regular\; sequence}\; x_{1},
\dots,x_{n} \in I\}.
$$
We define the {\em dimension} of $R$ to be
$$
\dim(R)=\sup\{s:{\rm there\; is\; a\; prime\; chain}\; {\mathcal P}_{s}
\subset
\dots\subset {\mathcal P}_{1}\subset {\mathcal P}_{0}\;\subset R\}.
$$
Lastly, $R$ is {\em Cohen-Macaulay} if, for every maximal
ideal ${\mathcal M}$ of $R$ (and such ideals are necessarily prime), we have
$$
{\rm depth}({\mathcal M})={\rm codim}({\mathcal M}).
$$
\end{definition}

\begin{example}
Consider $R:={\mathbf C}[x,y,z]$ with $p_{1}=xz$ and $p_{2}=yz$. 
The ideal $I:=(p_{1},p_{2})$ has the property ${\rm rad}(I)=I$ so that $R/I$
is the coordinate ring
of the zero locus of $p_{1}$ and $p_{2}$, which is made up of the
$(x,y)$-plane
and the $z$-axis. It is not hard to see that $\dim(R/I)=2\neq 1$, the ambient
dimension minus the number of equations. So the ring $R/I$ is not
Cohen-Macaulay. In
fact, at the origin the maximal ideal ${\mathcal M}=(x,y,z)/I$ is the first
term in a maximal descending prime chain $(x,y,z)/I$, $(y,z)/I$ and $(z)/I$
so that
${\rm codim}({\mathcal M})=2$. However, ${\rm depth}({\mathcal M})=1$, since
$x+z\ {\rm mod} (I)$, for instance, forms a maximal regular sequence in
${\mathcal M}$.
\end{example}

The following ingredient, on the other hand, generates many
Cohen-Macaulay rings.

\vspace{3mm}
\noindent {\bf FACT}(\cite[p455]{Ei}). If $p_{1},\dots,p_{n}$ form a regular sequence
in the ring $R:={\mathbf C}[z_{1},\dots,z_{L}]$ with ideal
$I=(p_{1},\dots, p_{n})$, then ${\rm codim}(I)=n$, the ring $R/I$ is Cohen-Macaulay, 
and $\dim(R/I)=L-n$.

\begin{remark} The {\bf FACT} can be interpreted geometrically. In the case
when ${\rm rad}(I)=I$, for instance, the
quotient ring $R/I$ is the coordinate
ring of an affine variety. This quotient ring being Cohen-Macaulay says that
every point of the variety is
cut out by $L-n$ functions (technically, in a maximal regular sequence
vanishing at the
point) in such a way that the co-dimension of the point is the expected
value $L-n$. The variety is
called a {\em complete intersection}, which is of equal dimension $L-n$ on
all of its irreducible components.
\end{remark}

We now come to the major recipe for inductively constructing 
Cohen-Macaulay rings in this paper.

\begin{proposition}\label{pro10.2}
If $p_{1},\dots,p_{n}$ are linearly independent homogeneous
polynomials of equal degree
$\geq 1$ in the
ring ${\mathbf C}[z_{1},
\dots,z_{L}]$ such that the ideal $(p_{1},\dots,p_{n-1})$ is prime and
such that $p_1,\dots,p_{n-1}$ form a regular sequence, then 
$p_{1},\cdots,p_{n}$ form a regular sequence. In
particular, the {\bf FACT} above implies that the quotient ring
\begin{equation*}
{\mathbf C}[z_{1},\dots,
z_{L}]/(p_{1},\dots,p_{n})
\end{equation*}
is Cohen-Macaulay.
\end{proposition}

\begin{proof}
We know $V_{n-1}^{\mathbf C}$ is irreducible since
$I_{n-1}:=(p_{1},\cdots,p_{n-1})$ is prime.  
Thus $p_{n}^{\mathbf C}$ cannot vanish identically on
$V_{n-1}^{\mathbf C}$. Otherwise the Nullstellensatz applied to
$p_{n}^{\mathbf C}$ on the prime $I_{n-1}$ would imply
$$
p_{n}^{\mathbf C}=p_{1}^{\mathbf C}f_{1}+\dots+p_{n-1}^{\mathbf C}f_{n-1}
$$
for some $f_{1},\dots,f_{n-1}\in {\mathbf C}[z_{1},\dots,
z_{L}]$. As shown in Proposition~\ref{pr:6:10}, we may assume
that $f_{1},\dots,f_{n-1}$ are
constant polynomials, because all of $p_1^\C,\dots,p_n^\C$ are homogeneous of the same 
degree $\geq 1$.
But this would imply that $p_{1}^{\mathbf C},\dots,
p_{n}^{\mathbf C}$
were linearly dependent, which is not the case by assumption.

Suppose there are $f,f_{1},\dots,f_{n-1}\in
{\mathbf C}[z_{1},\dots,z_{L}]$ such that
$$
p_{n}^{\mathbf C}f=p_{1}^{\mathbf C}f_{1}+\dots+p_{n-1}^{\mathbf C}f_{n-1}.
$$
Then $f|_{V_{n-1}^{\mathbf C}}\equiv 0$ since $p_{n}^{\mathbf C}$ does not
vanish identically on the irreducible $V_{n-1}^{\mathbf C}$. So once
more the Nullstellensatz applied to $f$ on $I_{n-1}$ implies
that
$$
f=p_{1}^{\mathbf C}g_{1}+\dots+p_{n-1}^{\mathbf C}g_{n-1}
$$
for some $g_{1},\dots,g_{n-1}\in
{\mathbf C}[z_{1},\dots,z_{L}]$.

Lastly,
$(p_{1}^{\mathbf C},\dots,p_{n}^{\mathbf C})\neq
{\mathbf C}[z_{1},\dots,z_{L}]$ since
$p_{1}^{\mathbf C},\dots,p_{n}^{\mathbf C}$ are all homogeneous of the same degree $\geq 1$.
This confirms
that $p_{1}^{\mathbf C},\dots,p_{n}^{\mathbf C}$ form a regular sequence.
\end{proof}

For our later applications on the \textit{variety} level,
Proposition~\ref{pro10.2} is not quite sufficient, because the ring
${\mathbf C}[z_{1},\dots,z_{L}]/(p_{1},\dots,p_{n})$ in the proposition,
though being Cohen-Macaulay, may have nilpotent elements, in which case
the ring is not the coordinate ring of an affine variety.  If the ring
contains no nilpotent elements, then it is called \textit{reduced}.

\begin{example}
Let $p_{1}=y-x^{2}$ and $p_{2}=y$ in ${\mathbf C}[x,y]$. The zero locus of
$p_{1}$ and $p_{2}$ is $\{(0,0)\}$. However, the Cohen-Macaulay ring
${\mathbf C}[x,y]/(p_{1},p_{2})$ has a nilpotent, namely,
$x\; {\rm mod}((p_{1},p_{2}))$. Geometrically, the parabola $y=x^{2}$
intersects $y=0$ with multiplicity $2$.
\end{example}

What we must
do now is to find conditions under which the quotient ring in
Proposition~\ref{pro10.2} is 
reduced, in which case the variety associated with
the ring is called a \textit{Cohen-Macaulay variety}.

\begin{proposition}\label{pro10.3}
Assume the hypotheses of Proposition~{\rm \ref{pro10.2}}. Let $J_{n}$ be the
subvariety of the
variety $V_{n}:=\{z \in \C^L :p_{1}(z)=0,\dots,p_{n}(z)=0\}$ where the Jacobian matrix
of $p_{1},\dots,p_{n}$ is not of rank $n$. If\; ${\rm codim}(J_{n})\geq 1$
in $V_{n}$, then the ring
${\mathbf C}[z_{1},\dots,z_{L}]/(p_{1},\dots,p_{n})$ is reduced.
\end{proposition}

\begin{proof}
This is just Serre's criterion of reducedness~\cite[p457]{Ei}.
\end{proof}

\begin{remark}
If we assume in Proposition~\ref{pro10.2} that
$J_{n-1}$, the
subvariety of
$V_{n-1}=\{z:p_{1}(z)=\dots=p_{n-1}(z)=0\}$ where the Jacobian of
$p_{1},\dots,p_{n-1}$ is not of rank $n-1$, is of co-dimension $\geq 2$
in $V_{n-1}$, then we can give a somewhat more geometric account
of
Proposition~\ref{pro10.3} as follows. (In fact, in our applications to follow,
${\rm codim}(J_{n-1})\geq 2$ always holds true.)
Let $R = \mathbf C[z_1,\dots,z_L]$, let $I = (p_1,\dots,p_{n-1})$ and
let $J = (p_n)$.
We must show $R/(I+J)$ has no
nilpotents. That is, whenever $f\in R$ satisfies
$$
f^{k}=p_{1}f_{1}+\dots+p_{n}f_{n}\in I+J
$$
for some $k$ and $f_{1},\dots,f_{n}$, we must have $f\in I+J$. We may assume
$f^{k}$ is not in $I$, or else we are done since then $f\in I$ by the primeness
of $I$. It follows that $f$ is nonzero on $V_{n-1}$ and is zero on $V_{n}$.

Let $V_{n}=W_{1}\cup\dots\cup W_{s}$ be the irreducible decomposition of
$V_{n}$ in $V_{n-1}$. We know ${\rm codim}(W_{i})=1$ in $V_{n-1}$ for all $i$.
Then by ${\rm codim}(J_{n})\geq 1$ in $V_{n}$ the polynomial $p_{n}$ cuts out
$W_{i}$ with multiplicity 1 for each $i$ (it comes down to the implicit
function theorem in calculus). That is, $p_{n}=0$ defines the
{\em divisor} $W_{1}+\dots+W_{s}$ in $V_{n-1}$.

Now since $f$ vanishes on $V_{n}$, the divisor defined by $f=0$ assumes
multiplicity $\geq 1$ on each
$W_{i}$. At this point the principle that says that the poles get cancelled
by the zeros seems to suggest that the rational function $f/p_{n}$ is regular
everywhere on $V_{n-1}$. This is certainly true if $V_{n-1}$ is
smooth~\cite[p129]{Sh1},
because the germs of local regular functions on $V_{n-1}$ then form
a unique factorization domain; more generally, the
{\em normality} of the variety suffices for the conclusion~\cite[p111]{Sh1}.
From this it
follows that $(f/p_{n})|_{V_{n-1}}=g$ for some regular $g$ on $V_{n-1}$. In
other words, $(f-p_{n}g)|_{V_{n-1}}\equiv 0$. Therefore,
$$
f-p_{n}g=p_{1}g_{1}+\dots+p_{n-1}g_{n-1}\in I
$$
by the primeness of $I$. We conclude that $f\in I+J$, proving the reducedness
of $R/(I+J)$.

It remains to ensure the normality of $V_{n-1}$, which is true if the
co-dimension of 
$J_{n-1}$ is at least $2$. This is a consequence of Serre's criterion of
primeness~\cite[p457]{Ei}, because $V_{n-1}$ is a Cohen-Macaulay variety
due to {\rm codim}$(I)=n-1$.
In any event we resort to Serre's criterion one way or another.
\end{remark}
The next proposition plays a vital role in the applications to follow.

\begin{proposition}\label{pro10.4}
We assume the hypotheses of Proposition~{\rm \ref{pro10.2}} and the notation
in
Proposition{\rm ~\ref{pro10.3}}. If\; {\rm codim}$(J_{n})\geq 2$ in $V_{n}$
and 
$V_{n}$ is connected, then $(p_{1},\dots,p_{n})$ is a prime ideal.
\end{proposition}

\begin{proof}
Proposition~\ref{pro10.3} asserts that $V_{n}$ is a
connected Cohen-Macaulay variety. Now $X_{n}$, the complement of $J_{n}$ in
$V_{n}$, is smooth on the one hand.
On the other hand, $X_{n}$ is also connected on account of Hartshorne's
connectedness theorem~\cite[p454]{Ei}, that says that a connected
Cohen-Macaulay
variety remains connected when a subvariety of co-dimension $\geq 2$ is
removed. Being both smooth and connected, $X_{n}$ must be irreducible.
However, since {\rm codim}$(J_{n})\geq 2$, $J_{n}$ cannot be an irreducible
component of $V_{n}$ due to the fact that a Cohen-Macaulay variety is of
equal dimension
on all of its irreducible components. $V_{n}$ is then irreducible. As a
consequence $(p_{1},\dots,p_{n})$ is a prime ideal because
Proposition~\ref{pro10.3} establishes the reducedness of $(p_{1},\dots,p_{n})$.
\end{proof}

\begin{example}
This example shows that ${\rm codim}(J_{n})\geq 2$ in $V_{n}$ is a must in
Proposition~\ref{pro10.4}.
Let $p_{1}=z$ and $p_{2}=x^{2}-y^{2}+z^{2}$ in ${\mathbf C}[x,y,z]$. Then
$V_{2}=\{(x,\pm x,0)\}$ and $J_{2}=\{(0,0,0)\}$, which is of co-dimension 1
in $V_{2}$. But $V_{2}$ is reducible albeit connected. It also illustrates
that the co-dimension $2$ condition in Hartshorne's connectedness theorem
cannot be improved to co-dimension $1$.
\end{example}

\section{The classification theorem}
We now return to the isoparametric case. For a given second order Darboux frame field~\eqref{eq:4:12a} along $\x$ on $U \subset M$, recall that we have,
for $1\leq a\leq m_{1}$, bi\-homo\-geneous polynomials 
$$
p_{a}=\sum_{\alpha,\mu=1}^{m_{2}} F^{\mu}_{\alpha a}x_{\alpha}y_{\mu}
$$
of bi-degree $(1,1)$ in the polynomial ring ${\mathbf R}[x_{1},\dots,
x_{m_{2}},y_{1},\dots,y_{m_{2}}]$, irreducible and linearly
independent if $m_{2}\geq m_{1}+2$ by Lemma~\ref{le9.1}.
Before proving the theorem, we first introduce a generalized spanning
property. For $n = 1,\dots,m_1$, we define the linear map
$S_{n}^{x}:{\mathbf R}^{m_{2}}\to {\mathbf R}^{n}$ by
\begin{equation}\label{eq11.1}
S_{n}^{x}(y) = \begin{pmatrix}
\sum_{\alpha}F^{1}_{\alpha 1}x_{\alpha}&
\dots&\sum_{\alpha}F^{m_{2}}_{\alpha 1}x_{\alpha} \\
\vdots&\dots&\vdots\\
\sum_{\alpha}F^{1}_{\alpha n}x_{\alpha}&
\dots&\sum_{\alpha}F^{m_{2}}_{\alpha n}x_{\alpha}\end{pmatrix}
\begin{pmatrix}y_{1}\\\vdots\\y_{m_{2}}\end{pmatrix} = \begin{pmatrix}
p_1(x,y) \\ \vdots \\ p_n(x,y) \end{pmatrix}
\end{equation}
and the linear map $S_{n}^{y}:{\mathbf R}^{m_{2}}\to
{\mathbf R}^{n}$ by
\begin{equation}\label{eq11.2}
S_{n}^{y}(x) = \begin{pmatrix}
\sum_{\mu}F^{\mu}_{1 1}y_{\mu}&
\dots&\sum_{\mu}F^{\mu}_{m_{2} 1}y_{\mu}\\
\vdots&\dots&\vdots\\
\sum_{\mu}F^{\mu}_{1 n}y_{\mu}&
\dots&\sum_{\mu}F^{\mu}_{m_{2} n}y_{\mu}\end{pmatrix}
\begin{pmatrix}x_{1}\\\vdots\\x_{m_{2}}\end{pmatrix} = \begin{pmatrix}
p_1(x,y) \\ \vdots \\ p_n(x,y) \end{pmatrix}
\end{equation}

\begin{definition}
We say that the $n$-\textit{spanning property} holds if
there is an $x\in {\mathbf R}^{m_{2}}$ such that $S_{n}^{x}$ is surjective and
there is a $y\in {\mathbf R}^{m_{2}}$ such that $S_{n}^{y}$ is surjective.
\end{definition}

Note that when $n=m_{1}$, this definition agrees with that of the spanning
property in Definition~\ref{spanning} for the second fundamental form (see Remark~\ref{re:6:8}).  As for the spanning property, the $n$-spanning property is an open condition.

We now set up an induction procedure toward our solution to~\eqref{eq:8:1}
and the spanning property.

\vspace{3mm}
\noindent {\bf Induction hypothesis ${\mathcal S}(n)$}
\begin{description}
\item[(I)] $p_{1},\dots,p_{n},n\leq m_{1},$ being irreducible and linearly
independent imply that $p_{1}^{\mathbf C},\dots,p_{n}^{\mathbf C}$ form
a regular sequence.
\item[(II)] $V_{n}:=\{z = (x,y) \in \R^{m_2} \times \R^{m_2} :p_{a}(z)=0, a = 1,\dots, n\}$ and
$V_{n}^{\mathbf C}:=\{z = (x,y) \in \C^{m_2} \times \C^{m_2} :p_{a}^{\mathbf C}(z)=0, a = 1, \dots, n\}$ satisfy
$\dim_{\mathbf R}(V_{n})=\dim_{\mathbf C} (V_{n}^{\mathbf C})=2m_{2}-n$, where 
$\dim_{\mathbf R} V_{n}$ is the maximal dimension of all the irreducible
components of $V_{n}$.
\item[(III)] $I_{n}:=(p_{1}^{\mathbf C},\dots,p_{n}^{\mathbf C})$ is a prime
ideal.
\item[(IV)] The $n$-spanning property is true.
\end{description}
Let $J_{n}$ be the subvariety of
$V_{n}^{\mathbf C}$ where the Jacobian matrix of
$p_{1}^{\mathbf C},\dots,p_{n}^{\mathbf C}$ is of rank $< n$.
Proposition~\ref{pro10.4} points out that {\rm codim}$(J_{n})\geq 2$ plays a
decisive role in determining the primeness of $I_{n}$. We will establish
in the next section the following estimate.
\begin{proposition}\label{pro11.1}
Assume $m_{1}\geq 2$. If $m_{2}\geq 3m_{1}-1$ , then ${\rm codim}(J_{n})\geq 2$
for all $n\leq m_{1}$.
\end{proposition}
Assuming this proposition for the time being, let us prove the classification
theorem of this paper.

\begin{theorem}[Classification]\label{th11.1} If $m_1 = 1$ or if $m_1
\geq 2$ and $m_2 \geq 3m_1-1$, then
the isoparametric hypersurface is of FKM-type.
\end{theorem}

\begin{proof}  When $m_1 = 1$, then $a=1$, $p=2$ and equations~\eqref{eq:5:6}
through~\eqref{eq:5:10} simplify sufficiently that one easily shows
that there exists a second order frame field for which
\[
\aligned
F^\mu_{\a\, a+m_1} &= \delta_{\a+m_2\, \mu} = F^\mu_{\a a} \\
F^\a_{pa} &= 0 = F^\mu_{pa}
\endaligned
\]
for all $\a$, $\mu$.  The first line of these equations
implies~\eqref{eq:8:1} and the 
spanning property.
Hence, Theorem~\ref{th:g:1} implies Takagi's
result~\cite{Ta} that all such isoparametric hypersurfaces are of
FKM-type.  

Suppose $m_{1}\geq 2$. Our strategy is to show that the induction
procedure can be completed 
for $n\leq m_{1}$. When $n=m_{1}$ what we achieve out of the induction is
that~\eqref{eq:8:1} and the spanning property hold true. It follows from
Theorem~\ref{th:g:1} that the isoparametric hypersurface is of FKM-type.

${\mathcal S}(1)$ is true. (I) holds because $p_{1}^{\mathbf C}$ is
irreducible by Lemma~\ref{le9.1}, and $p_{1}^{\mathbf C}$ cannot
generate the polynomial ring since it is of degree 2.
(II) is valid because $p_{1}$ is bi\-homo\-geneous of bi-degree (1,1), and so
one can easily solve for one variable in terms of the remaining ones regardless of
whether the variables are real or complex.
(III) is verified because $(p_{1}^{\mathbf C})$ is
a prime ideal due to the irreducibility of $p_1^{\mathbf C}$. (IV) is also clear since
$p_{1}\neq 0$.

Suppose ${\mathcal S}(n-1)$ is true for $n-1\leq m_{1}$. We show
${\mathcal S}(n)$ is true if $n\leq m_{1}$.
Now, (I) comes from Proposition~\ref{pro10.2}, so that 
the same proposition allows us to conclude that
${\mathbf C}[x_{1},\dots,x_{m_{2}},y_{1},\dots,y_{m_{2}}]/
(p_{1}^{\mathbf C},\dots,p_{n}^{\mathbf C})$ is Cohen-Macaulay.

We wish to
establish (II) next. To this end, note first that $V_{n}^{\mathbf C}$ is of
equal dimension $2m_{2}-n$
on all irreducible components, because
$V_{n}^{\mathbf C}$ is the intersection of the irreducible
$V_{n-1}^{\mathbf C}$ and the irreducible
hypersurface defined by $p_{n}^{\mathbf C}=0$. It follows that the real variety $V_{n}$
has the property
$$
\dim_{\mathbf R}(V_{n})\leq \dim_{\mathbf C}(V_{n}^{\mathbf C})=2m_2-n,
$$
because as established in Lemma~\ref{le10.1}, $V_{n}$ is a real subvariety
of $V_{n}^{\mathbf C}$ and any real subvariety is of dimension at most half
the (real) dimension of $V_{n}^{\mathbf C}$. We claim that there is a
component of $V_{n}$ having dimension $2m_{2}-n$ so that
$$
\dim_{\mathbf R}(V_{n})=\dim_{\mathbf C}(V_{n}^{\mathbf C}),
$$
which will establish (II). To prove the claim, consider 
$$
V_{n}\stackrel{\iota}{\longrightarrow}
{\mathbf R}^{m_{2}}\times {\mathbf R}^{m_{2}}\stackrel{\pi_{1}}
{\longrightarrow}{\mathbf R}^{m_{2}},
$$
where $\iota$ is the natural embedding and $\pi_{1}$ is the projection onto the
first summand.
Now $\pi_{1}\circ\iota$ is
surjective. This is because $(x,y)\in (\pi_{1}\circ\iota)^{-1}(x)$ precisely
when $y$ belongs to the kernel of the linear map $S_{n}^{x}$,
which has dimension $\geq m_{2}-n>0$; the set ${\mathcal L}$ of $x$ where this dimension
achieves the minimum value $t$ is Zariski open.
Since $\pi_{1}\circ\iota$ is
surjective, one of the irreducible components $W$ of $V_{n}$ must be mapped
onto an open subset of ${\mathcal L}$; or else Sard's theorem would imply
that $\pi_{1}\circ\iota$ is not surjective. Around a regular value $x$ of
$\pi_{1}\circ\iota$ in ${\mathcal L}$ we know $V_{n}$ is a product with fiber
${\mathbf R}^{t}$, which is therefore contained in the irreducible $W$.
Then since $t\geq m_{2}-n$, we have
$$
\dim(W)= m_{2}+t\geq m_{2}+m_{2}-n=2m_{2}-n.
$$
Therefore
$$
\dim_{\mathbf R}(V_{n})=2m_{2}-n=\dim_{\mathbf C}(V_{n}^{\mathbf C}),
$$
which proves (II).

Now that $\dim W=2m_{2}-n$, the fact that $V_{n}$ is a product with fiber
${\mathbf R}^{t}$ around the
regular value $x$ gives that
$$
\dim((\pi_{1}\circ\iota)^{-1}(x))=m_{2}-n.
$$
That is, $S_{n}^{x}$ spans ${\mathbf R}^{n}$. Likewise, there is some
$y\neq 0$ in ${\mathbf R}^{m_{2}}$ such
that $S_{n}^{y}$ spans ${\mathbf R}^{n}$ if we consider the projection
$\pi_{2}:{\mathbf R}^{m_{2}}\times{\mathbf R}^{m_{2}}
\longrightarrow {\mathbf R}^{m_{2}}$ onto the second summand.
In conclusion, we have shown that (IV) is true.

To finish the induction, we must show that $I_{n}$ is a prime ideal so that (III)
holds.
Proposition~\ref{pro10.4} and Proposition~\ref{pro11.1} tell us that this
is true if $V_{n}^{\mathbf C}$ is connected, which is the case because $V^\C_n$ is a cone.  In fact, if $z$ and $w$ are any two points in $V_n^\C$, then the real lines from $z$ to the origin and from the origin to $w$ are in $V^\C_n$, thus showing that $V^\C_n$ is path connected.

Thus, by Propositions~\ref{pro10.4} and~\ref{pro11.1}, the induction procedure is completed.

Setting $n=m_{1}$ in the induction, we obtain the spanning property in Definition~\ref{spanning} by
induction item
(IV). Note also that $V_{{m_{1}}}$ is exactly ${\mathcal D}$ defined
in~\eqref{eq9.4}.

We are only left with handling~\eqref{eq:8:1}.
By Proposition~\ref{pro9.3} we know $\overline{p}_{a},1\leq a\leq m_{1},$
vanish on ${\mathbf P}_{b}V_{m_{1}}$ so that
$\overline{p}_{a}|_{V_{m_{1}}}\equiv 0$, which warrants that
$\overline{p}_{a}^{\mathbf C}|_{V_{m_{1}}^{\mathbf C}}\equiv 0$ in view of the
induction item (II)
and Lemma~\ref{le10.1}, so that $\overline{p}_{a}^{\mathbf C}\in I_{m_{1}}$
by the induction item (III). 
Hence there are complex polynomials
$\tau_{ab},1\leq a,b\leq m_{1},$ such that
$$
\overline{p}_{a}^{\mathbf C}=\sum_{b=1}^{m_{1}}\tau_{ab}p_{b}^{\mathbf C}.
$$
As shown in the proof of Proposition~\ref{pr:6:10}, we may assume that the
$\tau_{ab}$ are constant polynomials, since each of the polynomials $\overline p^\C_a$ and $p_b^\C$ is of bi-degree $(1,1)$.  Restricting to the real variables we obtain
$$
\overline{p}_{a}=\sum_{b=1}^{m_{1}}f_{ab}p_{b}
$$
for some real constants $f_{ab}$.  The above argument establishes this at every point of the open set $U$ on which the frame is defined.  By Proposition~\ref{pr:6:10}, after a possible change of second order frame field along $\x$ on $U$, equation~\eqref{eq:8:1} holds on $U$.  
Theorem~\ref{th:g:1} then finishes the proof in the case $m_{1}\geq 2$.
\end{proof}

\begin{remark}\label{re11.1}
In contrast,
for $m_{1}=m_{2}=2$ of non-FKM-type, we have two pairs of
$(p_{1},p_{2})$ depending
on which one of the two focal submanifolds is referred to as $M_{+}$.
One pair of
$(p_{1},p_{2})=(0,0)$. The other pair is
$(2x_{2}y_{1}-2x_{1}y_{2},-2x_{1}y_{1}-2x_{2}y_{2})$, out of which the
real bi-projective variety ${\mathbf P}_{b}V_{2}$ is empty whereas the complex bi-projective
variety ${\mathbf P}_{b}V_{2}^{\mathbf C}$ consists of four points
$[1:\pm\sqrt{-1}]\times [1:\pm\sqrt{-1}]$. This case fails to satisfy
Proposition~\ref{pro10.1} miserably.
\end{remark}

\section{The estimate}
We now prove Proposition~\ref{pro11.1} to complete the classification
theorem in the preceding section. Recall for $V_{n}^{\mathbf C}$, its
subvariety $J_{n}$ is where the Jacobian matrix of $p_{1}^{\mathbf C},\dots,
p_{n}^{\mathbf C}$ fails to be
of rank $n$. From now on $S_{n}^{x}$ and $S_{n}^{y}$ in~\eqref{eq11.1}
and~\eqref{eq11.2} will be set in the complex category.

\begin{lemma}\label{le12.1}
Notation is as in~\eqref{eq:6:7}. For any choice of $a \in \{1,\dots,m_{1}\}$, there is an
orthonormal basis in $V_{+}$ and an orthonormal basis in $V_{-}$ such that
relative to these bases we have
\begin{description}
\item[(1)] $B_{a}=C_{a}$ with rank $=r\leq m_{1}$, and
\item[(2)] $A_{a}=\begin{pmatrix}I&0\\0&\Delta\end{pmatrix}$, where $\Delta$
is an $r\times r$ matrix in block form
$$
\Delta = \begin{pmatrix}\Delta_{1}&0&0&0&\dots\\0&\Delta_{2}&0&0&\dots
\\0&0&\Delta_{3}&0&\dots\\\vdots&\vdots&\vdots&\vdots&\vdots\end{pmatrix}
$$
with $\Delta_{1}=0$ and $\Delta_{i},i\geq 2,$ nonzero skew-symmetric
matrices in block form
$$
\Delta_i = \begin{pmatrix}0&f_{i}&0&0&\dots\\-f_{i}&0&0&0&\dots\\0&0&0&f_{i}&\dots\\
0&0&-f_{i}&0&\dots\\\vdots&\vdots&\vdots&\vdots&\vdots\end{pmatrix}.
$$
\end{description}
\end{lemma}

\begin{proof}
We know $B_{a}:V_{0}\longrightarrow V_{+}$, so that
$B_{a}{\tose B_{a}}:V_{+}\longrightarrow V_{+}$. Pick an orthonormal basis
$X_{1},\dots,X_{m_{2}-r},Y_{1},\dots,Y_{r}$ of $V_{+}$ for some $r$
such that
\begin{eqnarray}\label{eq12.0}
B_{a}{\tose B_{a}}&:&X_{t}\longmapsto 0,\\\nonumber
                  &:&Y_{s}\longmapsto (\sigma_{s})^{2}Y_{s},\nonumber
\end{eqnarray}
where $1\leq t\leq m_{2}-r$, $1\leq s\leq r$ and $\sigma_{s}>0$. Now
${\tose B_{a}}(X_{t})=0$ because ${\rm Ker}(B_{a})\cap {\rm Im}
({\tose B_{a}})=0$; hence $X_{t}\in {\rm Ker}({\tose B_{a}})$. That is,
${\rm Ker}({\tose B_{a}})$ is the eigenspace of $B_{a}{\tose B_{a}}$ with
eigenvalue zero. On the other hand, we know
$({\rm Ker}({\tose B_{a}}))^{\perp}={\rm Im}(B_{a})$. So the eigenspace
decomposition of $B_{a}{\tose B_{a}}$ is
$$
V_{+}={\rm Ker}({\tose B_{a}})\oplus {\rm Im}(B_{a})
$$
with $X_{1},\dots,X_{m_{2}-r}$ spanning the first summand and
$Y_{1},\dots,
Y_{r}$ spanning the second. As a result, it follows that $r={\rm rank}(B_{a})$.
Likewise,
$$
V_{0}={\rm Ker}(B_{a})\oplus {\rm Im}({\tose B_{a}}).
$$
We know from above that ${\tose B_{a}}(X_{t})=0$ and we set
\begin{eqnarray}\label{eq12.1}
{\tose B_{a}}:Y_{s}\longmapsto \sigma_{s}W_{s}
\end{eqnarray}
for some $W_{s}$. An easy calculation shows
$W_i \cdot W_j=\delta_{ij}$ so that $W_{1},\dots,W_{r}$
form an orthonormal basis of ${\rm Im}({\tose B_{a}})$. In conclusion,
$$
V_{0}={\rm Ker}(B_{a})\oplus {\rm Im}({\tose B_{a}}),
$$
where $W_{1},\dots,W_{r}$ span the second summand and we let $Z_{1},
\dots,Z_{m_{1}-r}$  be an orthonormal basis generating the first. We find
by~\eqref{eq12.0} that
\begin{eqnarray}\label{eq12.2}
B_{a}&:&Z_{t}\longmapsto 0,\\\nonumber
     &:&W_{s}\longmapsto \sigma_{s} Y_{s}.\nonumber
\end{eqnarray}
We calculate to see that ${\tose B_{a}}B_{a}:V_{0}\longrightarrow V_{0}$
satisfies
\begin{eqnarray}\label{eq12.3}
{\tose B_{a}}B_{a}&:&Z_{t}\longmapsto 0,\\\nonumber
                  &:&W_{s}\longmapsto (\sigma_{s})^{2}W_{s}.\nonumber
\end{eqnarray}
Now consider
$$
C_{a}:V_{0}\longrightarrow V_{-}.
$$
In the same manner as above for $B_{a}$, we get
$V_{0}={\rm Ker}(C_{a})\oplus {\rm Im}({\tose C_{a}})$ with
\begin{eqnarray}\label{eq12.4}
C_{a}&:&Z_{t}^{*}\longmapsto 0,\\\nonumber
     &:&W_{s}^{*}\longmapsto \sigma^{*}_{s}Y_{s}^{*},\nonumber
\end{eqnarray}
where $Z_{1}^{*},\dots,Z_{m_{1}-p}^{*}$ span ${\rm Ker}(C_{a})$ and
$W_{1}^{*},
\dots,W_{p}^{*}$ span ${\rm Im}({\tose C_{a}})$ for some $p$. However,
$$
{\tose C_{a}}C_{a}={\tose B_{a}}B_{a}
$$
by the first equation of~\eqref{eq:5:6}, we thus obtain
${\rm Ker}(B_{a})={\rm Ker}(C_{a})$ and 
${\rm Im}({\tose B_{a}})={\rm Im}({\tose C_{a})}$. 
In particular, $p=r$ and we may take $Z_{1},\dots,Z_{m_{1}-r}$ to be
identical with $Z_{1}^{*},\dots,Z_{m_{1}-r}^{*}$, and
$W_{1},\dots,W_{r}$ to be identical with $W_{1}^{*},\dots,W_{r}^{*}$.
Therefore~\eqref{eq12.2} and~\eqref{eq12.4} imply that we can pick a
basis of $V_{+}$ and a basis of $V_{-}$ relative to which the matrices of these operators, denoted by the same letters as the operators, satisfy
\begin{eqnarray}\label{eq12.45}
B_{a}=C_{a},
\end{eqnarray}
because from
\begin{eqnarray}\nonumber
{\tose C_{a}}C_{a}&:&Z_{t}^{*}\longmapsto 0,\\\nonumber
                  &:&W_{s}^{*}\longmapsto (\sigma_s^{*})^{2}W_{s}^{*}\nonumber
\end{eqnarray}
and $W_{s}=W_{s}^{*}$, we know $(\sigma_{s})^{2}=(\sigma^{*}_{s})^{2}$, and
hence we may assume $\sigma_{s}=\sigma_{s}^{*}$ by adjusting the basis in
$V_{-}$.

The second and the fourth equations of~\eqref{eq:5:6} together
with~\eqref{eq12.45} yield
\begin{eqnarray}\label{eq12.5}
A_{a}{\tosea A_{a}}={\tosea A_{a}}A_{a}=I-2B_{a}{\tose B_{a}}.
\end{eqnarray}
We have three more equations
\begin{eqnarray}
B_{a}{\tose B_{a}}{\tose A_{a}}+A_{a}B_{a}{\tose B_{a}}=0,\label{eq12.6}\\
B_{a}{\tose B_{a}}A_{a}+{\tose A_{a}}B_{a}{\tose B_{a}}=0,\label{eq12.7}\\
{\tose B_{a}}{\tose A_{a}}B_{a}+{\tose B_{a}}A_{a}B_{a}=0,\label{eq12.8}
\end{eqnarray}
which can be derived from~\eqref{eq12.45} and the three diagonal blocks of
the first equation of~\eqref{eq:6:12}. Let
$$
A_{a}=\begin{pmatrix}\alpha&\beta\\\gamma&\mu\end{pmatrix}
$$
where $\alpha$ is of size $(m_{2}-r)\times (m_{2}-r)$ and $\mu$ is of
size $r\times r$. Let $\sigma={\rm diag}(\sigma_{1},\dots,\sigma_{r})$ be
the diagonal matrix with the indicated diagonal entries so that
by~\eqref{eq12.1} and~\eqref{eq12.2}, $B_a$ and ${\T B_a}$ are of the
same form
\begin{equation}\label{block}
\begin{pmatrix}0&0\\0&\sigma\end{pmatrix},
\end{equation}
with
$$
B_{a}{\tose B_{a}}=\begin{pmatrix}0&0\\0&\sigma^{2}\end{pmatrix}
$$
of the same block sizes as $A_{a}$. From~\eqref{eq12.6} we obtain
\begin{eqnarray}
\beta=\gamma=0\label{eq12.9},\\
\sigma^{2}({\tose \mu})=-\mu \sigma^{2}.\label{eq12.10}
\end{eqnarray}
Moreover from~\eqref{eq12.5} we see
\begin{eqnarray}
\alpha{\tose \alpha}=I,\label{eq12.11}\\
\mu{\tose \mu}={\tose \mu}\mu=I-2\sigma^{2}\label{eq12.115}.
\end{eqnarray}
Similarly,~\eqref{eq12.7} yields
\begin{equation}
\sigma^{2}\mu=-{\tose \mu}\sigma^{2},\label{eq12.12}
\end{equation}
and~\eqref{eq12.8} gives
\begin{equation}
\sigma{\tose \mu}\sigma=-\sigma\mu\sigma.\label{eq12.13}
\end{equation}
With~\eqref{eq12.10} and~\eqref{eq12.12} we deduce
$$
\mu_{ij}=-(\sigma_{i}/\sigma_{j})^{2}\mu_{ji},
$$
and
$$
\mu_{ji}=-(\sigma_{i}/\sigma_{j})^{2}\mu_{ij}.
$$
We therefore conclude
$$
\mu_{ij}=0\;\; {\rm if}\; \sigma_{i}\neq \sigma_{j},
$$
and
$$
\mu_{ij}=-\mu_{ji}\;\; {\rm if}\; \sigma_{i}=\sigma_{j}.
$$
In other words,
$$
A_{a}=\begin{pmatrix}\alpha&0\\0&\mu\end{pmatrix}
$$
with
$$
\alpha{\tose \alpha}=I
$$
and $\mu$ is in blocked form
$$
\mu=\begin{pmatrix}\Delta_{1}&0&0&0&\dots\\0&\Delta_{2}&0&0&\dots
\\0&0&\Delta_{3}&0&\dots\\\vdots&\vdots&\vdots&\vdots&\vdots\end{pmatrix},
$$
where all the $\Delta_{i}$ are skew-symmetric such that the number of
$\Delta_{i}$ is the number of different non-zero eigenvalues of $B_{a}{\tose B_{a}}$.
Then~\eqref{eq12.13} is automatically satisfied. Now by the skew-symmetry of
$\mu$ and~\eqref{eq12.115} we derive
\begin{equation}
\Delta_{i}^{2}=-(1-2\sigma_{i}^{2})I.\label{eq12.14}
\end{equation}
In view of~\eqref{eq12.11} and the skew-symmetry of $\mu$ 
we can perform an orthonormal basis change so that
$$
\alpha=I
$$
and
$$
\Delta_{i}=
\begin{pmatrix}0&r_{1}&0&0&\dots\\-r_{1}&0&0&0&\dots\\0&0&0&r_{2}&\dots\\
0&0&-r_{2}&0&\dots\\\vdots&\vdots&\vdots&\vdots&\vdots\end{pmatrix}.
$$
Thus~\eqref{eq12.14} implies $r_{1}^{2}=r_{2}^{2}=\dots=1-2\sigma_{i}^{2}$,
and so
$$
\Delta_{i}=\sqrt{1-2\sigma_{i}^{2}}
\begin{pmatrix}0&1&0&0&\dots\\-1&0&0&0&\dots\\0&0&0&1&\dots\\
0&0&-1&0&\dots\\\vdots&\vdots&\vdots&\vdots&\vdots\end{pmatrix},
$$
if $1-2\sigma_{i}^{2}>0$. We set $\Delta_{1}\equiv 0$ so that
$\sigma_{1}=1/\sqrt{2}$. We are done.
\end{proof}

\begin{corollary}\label{co12.1}
$\dim(Ker(A_{a}))=\dim(\Delta_{1})\leq r={\rm rank}(B_{a})\leq m_{1}$.
\end{corollary}

\begin{remark}\label{re12.1}
When $(m_{1},m_{2})=(2,m),m\geq 3$, Ozeki and Takeuchi
showed~\cite[II, p49]{OT}, that $r$ given in Lemma~\ref{le12.1} is $1$,
essentially by exploring the fact that $p_{1}$
and $p_{2}$ form a regular sequence in the spirit of Example~\ref{example2}
above.
It follows immediately from Lemma~\ref{le12.1} that we have $\Delta=0$
and so
$$
A_{1}=\begin{pmatrix}I&0\\0&0\end{pmatrix}
$$
as given in~\cite{OT}. With this it is not hard to see~\cite[II, p51]{OT},
that
$$
A_{2}=\begin{pmatrix}B&0\\0&0\end{pmatrix}
$$
of the same block sizes as $A_{1}$ with
$$
B=\begin{pmatrix}0&-I\\I&0\end{pmatrix},
$$
where $I$ in $B$ is of size $l\times l$ and $m_{2}=2l+1$.
\end{remark}

\begin{proof}[Proof of Proposition~{\rm \ref{pro11.1}}]  We must
estimate the codimension in $V^\C_n$ of
\[
J_n = \{(x,y)\in V_n^\C : dp_1^\C \wedge \dots \wedge dp_n^\C = 0\}
\]
We first
estimate the dimension of the 
subvariety $Z_{n}$ of  ${\mathbf C}^{m_{2}}\times {\mathbf C}^{m_{2}}$ at
each point
of which the Jacobian matrix of $p_{1}^{\mathbf C},\dots,p_{n}^{\mathbf C}$
is of rank $<n$.
At $(x,y)\in Z_{n}$, the differentials
$dp_{1}^{\mathbf C},\dots,dp_{n}^{\mathbf C}$ are linearly
dependent, i.e., there are $c_{1},\dots,c_{n}\in {\mathbf C}$, depending on
$(x,y)$, such that
$$
0=\sum_{a=1}^{n}c_{a}dp_{a}^{\mathbf C}=\sum_{\alpha} (\sum_{a,\mu}
c_{a}F^{\mu}_{\alpha a}y_{\mu})dx_{\alpha}
+\sum_{\mu}(\sum_{a,\alpha}
c_{a}F^{\mu}_{\alpha a}x_{\alpha})dy_{\mu},
$$
which requires that the coefficients of $dx_{\alpha}$ be zero and the
coefficients of $dy_{\mu}$ be zero.
Thus
$$
Z_{n}=\{(x,y)\in
{\mathbf C}^{m_{2}}\times {\mathbf C}^{m_{2}}
:\exists (c_{1},\dots,c_{n}),
\sum_{a}c_{a}{\tose A_{a}}x=\sum_{a}c_{a}A_{a}y=0\}.
$$
Accordingly, for a fixed $(c_{1},\dots,c_{n})$ let us define
$$
Z_{(c_{1},\dots,c_{n})}:=\{(x,y)\in
{\mathbf C}^{m_{2}}\times {\mathbf C}^{m_{2}}
:\sum_{a}c_{a}{\tose A_{a}}x=\sum_{a}c_{a}A_{a}y=0\}.
$$
Consider the incidence space $Y_{n}$ in
${\mathbf C}P^{n-1}\times {\mathbf C}^{m_{2}}\times {\mathbf C}^{m_{2}}$
given by
\begin{equation}\label{incidence}
Y_n = \{([c_{1}:\cdots:c_{n}],x,y): (x,y)\in Z_{(c_{1},\dots,c_{n})}\}.
\end{equation}
The standard projection of $Y_{n}$ to
${\mathbf C}^{m_{2}}\times {\mathbf C}^{m_{2}}$ maps $Y_{n}$ onto $Z_{n}$.
Let $\pi$ be the standard projection of $Y_{n}$ to ${\mathbf C}P^{n-1}$.
Then with respect to $\pi$ we have
\begin{equation}\label{eq12.145}
\dim(Z_{n})\leq \dim(Y_{n})\leq \dim({\rm base})+\dim({\rm fiber}),
\end{equation}
where $\dim({\rm fiber})$ is the
maximal dimension of all fibers.
We first estimate the dimension
of the fibers $\pi^{-1}\{[c_1:\cdots :c_n]\} =
Z_{(c_{1},\dots,c_{n})}$. In fact, it comes down to estimating the
dimension of 
$$
T_{(c_{1},\dots,c_{n})}:=\{y\in {\mathbf C}^{m_{2}}:\sum_{a}c_{a}A_{a}y=0\}
$$
for a fixed $(c_{1},\dots,c_{n})$, because 
\begin{equation}\label{equaldims}
\dim (\mbox{ker}\,(\sum_a c_a \T A_a)) = \dim (\mbox{ker}\,(\sum_a c_a
A_a))
\end{equation}
thus giving us the estimate 
$$
\dim(Z_{(c_{1},\dots,c_{n})})\leq 2\dim(T_{(c_{1},\dots,c_{n})}).
$$
\begin{remark}\label{re12.2}
Let us examine the case $(m_{1},m_{2})=(2,m_2),m_2\geq 3$, before we proceed.
By the above standard matrix form
of $A_{1}$ and of $A_{2}$ in Remark~\ref{re12.1} we see that for $\tose y = (\tose z, t) \in \mathbf C^{m_2}$, where $t \in \mathbf C$,
$$
A_{1}\begin{pmatrix}z\\t\end{pmatrix}=\begin{pmatrix}z\\0\end{pmatrix}, \quad
A_{2}\begin{pmatrix}z\\t\end{pmatrix}=\begin{pmatrix}Bz\\0\end{pmatrix}
$$
Hence 
$\sum_{a=1}^{n=2}c_{a}A_{a}y=0$ precisely when $z$ is an eigenvector of $B$,
with eigenvalue $\pm\sqrt{-1}$. In other words, when $[c_{1}:c_{2}]
=[\pm \sqrt{-1}:1]$ in ${\mathbf C}P^{1}$, $Z_{(c_{1},c_{2})}$ is made
up of vectors of the form
$$(\begin{pmatrix}z\\t\end{pmatrix},\begin{pmatrix}w\\s\end{pmatrix})
$$
where $z$ and $w$ both belong to the $\sqrt{-1}$-eigenspace or to the
$-\sqrt{-1}$-eigenspace of $B$ and $Z_{(c_{1},c_{2})}=\{(0,0)\}$ for other
values of $[c_{1}:c_{2}]$. Thus
\begin{equation}\label{sharp}
\dim(Z_{2})=m_{2}+1.
\end{equation}
\end{remark}
We continue on now to estimate the dimension
of $Z_{(c_{1},\dots,c_{n})}$.

\noindent Case (1). $c_{1},\dots,c_{n}$ are either all real or all purely
imaginary. Say it is the latter, so that $c_{k}=\sqrt{-1}d_{k}$ with
$d_{k}$ real. Then for $y\in T_{(c_{1},\dots,c_{n})}$, we have
$$
\sum_{k=1}^{n}d_{k}A_{k}y=0.
$$
However, the second fundamental form $S$ has the property
$$
d_{1}S_{e_{1}}+\dots+d_{n}S_{e_{n}}=\sqrt{d_{1}^2+\dots+d_{n}^{2}}S_{e},
$$
where
$$
e=(d_{1}e_{1}+\dots+d_{n}e_{n})/\sqrt{d_{1}^2+\dots+d_{n}^{2}}.
$$
We may therefore rename $e$ to be $e_{1}$ in the normal basis, and so by
restricting to the $A$-block in the matrix of $S$ we see that $S_{e}y=0$
comes down to, after the renaming, $A_{1}y=0$. Corollary~\ref{co12.1} then
establishes that
$$
\dim(T_{(c_{1},\dots,c_{n})})\leq r\leq m_{1}
$$
and
\begin{equation*}
\dim (Z_{(c_1,\dots, c_n)}) \leq 2 \dim T_{(c_1,\dots, c_n)} \leq 2m_1
\end{equation*}
Case (2). $c_{1},\dots,c_{n}$ are not all real and not all purely imaginary.
Write
$$
c_{k}=\alpha_{k}+\sqrt{-1}\beta_{k},
$$
where not all $\alpha_{k}$ and not all $\beta_{k}$ are zero.
Then
$$
c_{1}S_{e_{1}}+\dots+c_{n}S_{e_{n}}=(\alpha_{1}S_{e_{1}}
+\dots+\alpha_{n}S_{e_{n}})+\sqrt{-1}(\beta_{1}S_{e_{1}}
+\dots+\beta_{n}S_{e_{n}}).
$$
As in Case (1), we know $\alpha_{1}S_{e_{1}}+\dots+\alpha_{n}S_{e_{n}}$
is a multiple of $S_{e}$ for some unit vector $e$. Hence without loss of generality
we may assume, after renaming $e$ to be $e_{1}$, that
$$
c_{1}S_{e_{1}}+\dots+c_{n}S_{e_{n}}=\alpha_{1}S_{e_{1}}
+\sqrt{-1}(\beta_{1}S_{e_{1}}+\dots+\beta_{n}S_{e_{n}}).
$$
On the other hand $\beta_{2}S_{e_{2}}+\dots+\beta_{n}S_{e_{n}}$ is a
multiple of $S_{f}$ for some unit vector $f$ perpendicular to $e_{1}$. We rename $f$
to be $e_{2}$ so that we may assume without loss of generality that
$$
c_{1}S_{e_{1}}+\dots+c_{n}S_{e_{n}}=(\alpha_{1}
+\sqrt{-1}\beta_{1})S_{e_{1}}+\sqrt{-1}\beta_{2}S_{e_{2}}.
$$
By restricting to the $A$-block in $S$ again we see that 
$(\sum_{a}c_{a}A_{a})y=0$ is reduced to
$$
\beta_{2}A_{2}y
=\sqrt{-1}(\alpha_{1}+\sqrt{-1}\beta_{1})A_{1}y.
$$
We may assume both coefficients are nonzero, or else
we would be back to Case (1). Hence we are now handling
\begin{equation}\label{eq12.15}
(A_{2}-zA_{1})y=0
\end{equation}
for some nonzero $z\in {\mathbf C}$. By Lemma~\ref{le12.1}, we may assume
$$
A_{1}=\begin{pmatrix}I&0\\0&\Delta\end{pmatrix}.
$$
Write
$$
A_{2}=\begin{pmatrix}\Theta&\Lambda\\\Omega&\Gamma\end{pmatrix}
$$
of the same block sizes as $A_{1}$. By the second equation of~\eqref{eq:5:6},
which is
$$
A_{2}{\T A_{1}}+A_{1}{\T A_{2}}+2(B_{2}{\T B_{1}}
+B_{1}{\T B_{2}})=0,
$$
we obtain
\begin{equation}\label{eq12.16}
\Theta+{\tose \Theta}=0
\end{equation}
when we invoke~\eqref{block}.
If we write
$$
y=\begin{pmatrix}u\\v\end{pmatrix}, \quad u \in \mathbf C^{m_2-r}, \quad v \in \mathbf C^r
$$
then part of~\eqref{eq12.15} reads, 
\begin{equation}\label{eq12.17}
(zI-\Theta)u=\Lambda v.
\end{equation}
Consider the map $G:{\mathbf C}^{m_{2}}\longrightarrow {\mathbf C}^{m_{2}-r}$
given by
$$
G:(u,v)\longmapsto (zI-\Theta)u-\Lambda v.
$$
The kernel of $G$ consists of all $y=\tose (u,v)$ satisfying~\eqref{eq12.17}.
If $z$ is not an eigenvalue of $\Theta$, then the rank of $G$ is at least
the rank of $zI-\Theta$, which is $m_{2}-r$.  Thus, the rank of $G$ is $m_{2}-r$, so that
the kernel of $G$ has dimension $r$. On the other hand if $z$ is an eigenvalue of
$\Theta$, then because $\Theta$ is skew-symmetric by~\eqref{eq12.16}, the
rank of $zI-\Theta$ is at least $(m_{2}-r)/2$ due to the fact that a nonzero
eigenvalue of $\Theta$ is purely imaginary, and its conjugate is also an
eigenvalue of $\Theta$. It follows that the rank of $G$ is no less than
$(m_{2}-r)/2$, so that its kernel is of dimension $\leq (m_{2}+r)/2$.
The upshot is that, since $r\leq m_{1}$ and since
$\dim(T_{c_1,\dots,c_n})$ is an integer, we have arrived at the estimate
$$
\dim(T_{(c_{1},\dots,c_{n})}) \leq [(m_{2}+r)/2] \leq [(m_2+m_1)/2]
= (m_{2}+m_{1}-1)/2,
$$
where $[p]$ is the greatest integer in the number $p$, and the last
equality is true
because $m_{2}+m_{1}$ is an odd number when $2\leq m_{1}<m_{2}$
by a result of M\"{u}nzner~\cite[II]{Mu}. Therefore,
\begin{equation}\label{extra}
\dim(\mbox{fiber})=\dim(Z_{(c_{1},\dots,c_{n})})\leq
2\dim(T_{(c_{1},\dots,c_{n})})\leq 
m_{2}+m_{1}-1.
\end{equation}
This estimate is sharp in light of~\eqref{sharp}. Note that $m_2+m_{1}-1$ is greater than 
the upper bound $2m_1$ for $\dim (Z_{(c_1,\dots, c_n)})$ in Case~(1),
since $m_2 \geq 3m_1-1$ and $m_1 \geq 2$ by assumption.

We next stratify the incidence space $Y_{n}$ of~\eqref{incidence} in
another way as follows. 
We let $s\leq m_{2}$ be the largest
integer for which $\sum_{i=1}^{n}c_{i}A_{i}$ is of rank $s$ for some, and
hence for generic, $[c_{1}:\cdots : c_{n}]$, the set of which
constitute a Zariski 
open set $U$ of 
${\mathbf C}P^{n-1}$. 
A look at Corollary~\ref{co12.1} shows that
$$
s\geq m_{2}-m_{1},
$$
so that for $(c_{1},\cdots,c_{n})$ in $U$, 
$$
{\rm rank}(\sum_{i=1}^{n}c_{i}A_{i})=s\geq m_{2}-m_{1},
$$
and thus, by~\eqref{equaldims},
$$
\aligned
\dim(\mbox{fiber}) &= \dim(Z_{(c_1,\dots,c_n)}) = \dim(\ker(\sum_1^n
c_iA_i)) + \dim(\ker(\sum_1^n c_i \T A_i)) \\
&= 2( m_2 - \mbox{rank}(\sum_1^n c_iA_i)) = 2(m_2-s) \leq 2m_{1}.
\endaligned
$$
It follows that over $U$, \eqref{eq12.145} extends to 
\begin{equation}\label{eq200}
\dim(\mbox{fiber})+\dim(\mbox{base})\leq 2m_{1} + (n-1),
\end{equation}
On the other hand, over a subvariety $W$, contained in ${\mathbf C}P^{n-1}$,
of dimension $\leq n-2$, the rank of $\sum_{i=1}^{n}c_{i}A_{i}$ is less than
$s$. 
In view of~\eqref{extra}, we have that over $W$
\begin{equation}\label{eq201}
\aligned
\dim(\mbox{fiber})+\dim(\mbox{base})&\leq \dim(\mbox{fiber})+n-2\\&\leq m_{1}+m_{2}-1+n-2\\
&=m_{1}+m_{2}+n-3.
\endaligned
\end{equation}
Now the part of $Y_{n}$ over $U$, call it $A$, is irreducible because each fiber over 
$U$ is a Euclidean space of a fixed dimension, whereas the part over $W$, call it $B$, is 
Zariski closed in $Y_{n}$. It follows that the closure of $A$, call it $\overline{A}$, in 
$Y_{n}$ is an irreducible 
component of $Y_{n}$, and the closure of $B$ not in ${\overline A}$
constitutes the
remaining irreducible components in $Y_{n}$. Therefore, the larger of
the two upper bounds
in~\eqref{eq200} and~\eqref{eq201} will be an upper bound for the
dimension of $Y_{n}$, and
hence of $Z_{n}$. However, 
$2m_{1} + n-1 < m_{1}+m_{2}+n-3$, due to
$m_{2}\geq 3m_{1}-1$ and $m_1 \geq 2$. We conclude that over 
${\mathbf C}P^{n-1}$
$$
\dim(Z_{n})\leq m_{1}+m_{2}+n-3
$$
if $m_{2}\geq 3m_{1}-1$ and $m_1 \geq 2$.

Now $J_{n}$, the subvariety of $V_{n}^{\mathbf C}$ where
$dp_{1}^{\mathbf C},\dots,p_{n}^{\mathbf C}$ are dependent, is clearly a
subvariety of $Z_{n}$. Hence
$$
\dim(J_{n})\leq \dim(Z_{n})\leq m_{1}+m_{2}+n-3.
$$
On the other hand, $\dim(V_{n}^{\mathbf C})\geq 2m_{2}-n$ on all of its
irreducible components because $V_{n}^{\mathbf C}$ is cut out from
${\mathbf C}^{m_{2}}\times {\mathbf C}^{m_{2}}$ by $n$
equations~\cite[p59]{Sh1}.
It follows that ${\rm codim}(J_{n})\geq 2$ in $V_{n}^{\mathbf C}$ of dimension
at least $2m_{2}-n$ if $m_{1}+m_{2}+n-3\leq 2m_{2}-n-2$, i.e., if
$m_{2}\geq m_{1}+2n-1$, which is true if $m_{2}\geq 3m_{1}-1$
since $n\leq m_{1}$. This finishes the proof of Proposition~\ref{pro11.1}.
\end{proof}

The classification result Theorem~\ref{th11.1} is therefore established.

\begin{remark}\label{re:51}
The standard matrix form of $A_{1}$ and of $A_{2}$ in the case
$(m_{1},m_{2})=(2,m_2),m_2 \geq 3$, give~\cite[II, p51]{OT} that
\begin{eqnarray}\nonumber
p_{1}&=&-2\sum_{j=1}^{l}(x_{j}y_{j}+x_{l+j}y_{l+j}),\nonumber\\
p_{2}&=&2\sum_{j=1}^{l}(x_{j}y_{l+j}-x_{l+j}y_{j}),\nonumber
\end{eqnarray}
where as before $m_{2}=2l+1$. It turns out that $J_{2}=Z_{2}$. This is
because any element $X$ of $Z_{2}$ is either $(0,0)$ or is of the form
$$
(\begin{pmatrix}u\\\pm\sqrt{-1} u\\t\end{pmatrix},
\begin{pmatrix}v\\\pm\sqrt{-1} v\\s\end{pmatrix}),
$$
where $u$ is of size $l\times 1$. It is immediate to verify that $X$ is
annihilated by both $p_{1}^{\mathbf C}$ and $p_{2}^{\mathbf C}$.
It follows that $\dim(J_{2})= \dim(Z_2) = m_{2}+1$, as shown in~\eqref{sharp}.  Thus, 
$\mbox{codim}\,(J_{2})\geq 2$ in $V_{2}^{\mathbf C}$ (which is of dimension
$2m_{2}-2$), provided that $m_{2}\geq 5$, which is exactly
equal to $3m_{1}-1$ given in our classification theorem.

When $m_1 = 2$, our approach misses only the case $m_{2}=3$, so $l=1$
(which is of FKM-type by Ozeki-Takeuchi~\cite[II]{OT}).  This is not
surprising 
in view of the fact that in this case $m_2 = 3$, the bi-projective variety
${\mathbf P}_{b}V_{2}^{\mathbf C}$ defined by
$p_{1}^{\mathbf C}=p_{2}^{\mathbf C}=0$ in
${\mathbf C}P^{2}\times {\mathbf C}P^{2}$ is made up of six irreducible components
$$
\aligned
&\{[1:\pm\sqrt{-1}:z]\times [1:\pm \sqrt{-1}:w] : z,w \in \C \},\\
&\{[0:0:1]\}\times {\mathbf C}P^{2},\\
&{\mathbf C}P^{2}\times \{[0:0:1]\},
\endaligned
$$
so that $(p_1^{\mathbf C}, p_2^{\mathbf C})$
is not a prime ideal and thus Proposition~\ref{pro10.4} says then that 
$\mbox{codim}\,(J_2) \leq 1$ in $V_2^{\mathbf C}$.
\end{remark}

In view of the known classification of Takagi~\cite{Ta} for $m_{1}=1$,
Ozeki-Takeuchi~\cite[II]{OT} for
$m_{1}=2$,
and Stolz's result~\cite{St} on the multiplicities $m_{1}\leq m_{2}$ that
states that $(m_{1},m_{2})\neq (2,2)$ or $(4,5)$ must be that of an
isoparametric hypersurface of FKM-type, we obtain from
Theorem~\ref{th11.1} that all isoparametric hypersurfaces with four
principal curvatures in spheres, whose multiplicities are not
$(2,2)$ or $(4,5)$, are of FKM-type, except possibly for
those whose multiplicities are one of the following 9 pairs
$(3,4)$, $(4,7)$, $(5,10)$,
$(6,9)$, $(7,8)$, $(7,16)$, $(8,15)$, $(9,22)$, $(10,21)$. The $(4,5)$
case also remains open.


\end{document}